\documentclass[showpacs,preprint,superscriptaddress,amsmath,amssymb]{revtex4-1}
\allowdisplaybreaks[4]
\UseRawInputEncoding
\usepackage{txfonts}
\usepackage{dcolumn}
\usepackage{bm}
\usepackage{booktabs}
\usepackage{amssymb}
\usepackage{amsmath}
\usepackage{array}
\usepackage{color}
\usepackage{tabularx}
\usepackage{multirow} 
\usepackage{graphicx}
\usepackage{epstopdf}
\usepackage{subfig} 

\usepackage{xcolor}
\usepackage{algorithm}
\usepackage{algpseudocode}

\begin{document} 
	
	\title{Regularized lattice Boltzmann method based maximum principle and energy stability preserving finite-difference scheme for the Allen-Cahn equation}
 
	\author{Ying Chen}
	\affiliation{School of Mathematics and Statistics, Huazhong
		University of Science and Technology, Wuhan 430074, China}%
			\author{Xi Liu}
		\affiliation{School of Mathematics and Statistics, Huazhong
			University of Science and Technology, Wuhan 430074, China}%
	
	\author{Zhenhua Chai}
	\email[Corresponding author]{(hustczh@hust.edu.cn)}
	\affiliation{School of Mathematics and Statistics, Huazhong
		University of Science and Technology, Wuhan 430074, China}%
	\affiliation{Institute of Interdisciplinary Research for Mathematics and Applied Science,
		Huazhong University of Science and Technology, Wuhan 430074, China}%
	\affiliation{Hubei Key Laboratory of Engineering Modeling and
		Scientific Computing, Huazhong University of Science and Technology,
		Wuhan 430074, China}
	
	\author{Baochang Shi}
	\affiliation{School of Mathematics and Statistics, Huazhong
		University of Science and Technology, Wuhan 430074, China}%
	\affiliation{Institute of Interdisciplinary Research for Mathematics and Applied Science,
		Huazhong University of Science and Technology, Wuhan 430074, China}%
	\affiliation{Hubei Key Laboratory of Engineering Modeling and
		Scientific Computing, Huazhong University of Science and Technology,
		Wuhan 430074, China}
	
	\date{\today} 
	
	\begin{abstract}
	 The Allen-Cahn equation (ACE) inherently possesses two crucial properties: the maximum principle and the energy dissipation law.  Preserving these two properties at the discrete level is also necessary in the numerical methods for the ACE.  In this paper, unlike the traditionally top-down macroscopic numerical schemes which discretize the ACE directly, we first propose a novel bottom-up mesoscopic regularized lattice Boltzmann method based  macroscopic numerical scheme for $d$ (=1, 2, 3)-dimensional ACE, where the D$d$Q$(2d+1)$ [($2d+1$) discrete velocities in $d$-dimensional space] lattice structure is adopted. In particular, the proposed macroscopic numerical scheme has a second-order accuracy in space, and can also be viewed as an implicit-explicit finite-difference scheme for the ACE, in which the nonlinear term is discretized semi-implicitly, the temporal derivative and dissipation term of the ACE are discretized by using the explicit Euler method and second-order central difference method, respectively.  Then we  demonstrate that the proposed scheme can preserve the maximum bound principle and the original energy dissipation law at the discrete level under some conditions.  Finally, some numerical experiments are conducted to validate our theoretical analysis. 
	\end{abstract}
	
	\maketitle
	\section{Introduction}\label{sec:introduction}
	The Allen-Cahn equation (ACE), as one of the popular phase-field type models, has been widely used to describe some important physical phenomena in nature and engineering, for instance, the mean curvature flows \cite{ilmanen1993convergence, feng2003numerical, lee2015mean}, crystal growth \cite{li2023second}, image segmentation \cite{benevs2004geometrical}, to name but a few. In particular, it has also become a fundamental model in the diffuse interface approaches for phase transitions and interface dynamics \cite{kim2012phase, Wang2019}. However, it should be noted that the ACE is a nonlinear partial differential equation (PDE), and  it is  challenging to obtain its exact solution. Therefore, it is desirable and necessary to develop some numerical methods for the ACE.  In this work, we will focus on  the following $d$-dimensional ACE:
	 \begin{align}\label{ACE}
	 	\begin{cases}
	 		\frac{\partial \phi}{\partial t}=\varepsilon^2\nabla^2\phi-f(\phi),&\mathbf{x}\in\Omega,t\in[0,T_t],\\
	 		\phi(\mathbf{x},0)=\phi_0(\mathbf{x}),&\mathbf{x}\in\overline{\Omega},
	 	\end{cases}
	 \end{align}
	 where $\phi(\mathbf{x},t)$ represents the unknown function, $\Omega=\Pi_{i=1}^d(L_i,R_i)$ is a hypercube domain in $\mathbb{R}^d$, the interfacial  parameter $\varepsilon>0$ and the nonlinear term $f(\phi)=\phi^3-\phi$ is a cubic polynomial. It should be noted that to make Eq. (\ref{ACE}) well-posed, some proper boundary conditions should be imposed. For simplicity, here we only consider the homogeneous Neumann boundary condition (HN-BC), homogeneous Dirichlet boundary condition (HD-BC) and periodic boundary condition (P-BC). 
	 
	It is well known that the ACE (\ref{ACE}) has two important properties: the one is the maximum principle, which means that if the absolute values of the initial and boundary data are bounded by the constant 1, then the absolute value of the solution  at any time $t$ is also bounded by 1 \cite{yang2018uniform}, i.e.,
	\begin{align}\label{maximum}
	\big\|\phi(\cdot,t)\big\|_{\infty}\leq 1\quad \forall t>0.
	\end{align} 
	The other one is the energy dissipation law. Actually, as a phase-field type model, the ACE (\ref{ACE}) can be viewed as an $L^2$ gradient flow of a total free energy functional, i.e.,
	 \begin{align}
	 \partial_t\phi=-\frac{\delta E}{\delta \phi}=\varepsilon^2\nabla^2\phi-f(\phi),
	 \end{align}
	where the integration by parts and proper boundary condition, such as the HN-BC ($\mathbf{n}\cdot\nabla\phi=0$ with $\mathbf{n}$ as the outer unit normal vector on the boundary $\partial \Omega$), have been used. The free energy $E(\phi)$ is given by
	\begin{align}
		E(\phi)=\int_{\Omega} \Big(F(\phi)+\frac{\varepsilon^2}{2}|\nabla\phi|^2\Big)d\mathbf{x},
	\end{align}
  $F(\phi)$ is the bulk energy and has a double-well form:
	 \begin{align}
	 	F(\phi)=\frac{1}{4}\big(\phi^2-1\big)^2.
	 \end{align}
	 Taking the inner product between Eq. (\ref{ACE}) and the variational derivative of energy $E(\phi)$ in $L^2$ space, we can obtain the energy dissipation law as
	 \begin{align}\label{energy-diss}
	 	\frac{\partial E(\phi)}{\partial t}=\bigg<\phi_t,\frac{\delta E(\phi)}{\delta \phi}\bigg>=-\int_{\Omega}|\phi_t|^2d\mathbf{x}\leq 0.
	 \end{align} 
	From Eq. (\ref{energy-diss}), one can see that the total energy is non-increasing with time, i.e.,
	\begin{align}\label{energy}
		 E\big(\phi(\cdot,t_n)\big)\leq E\big(\phi(\cdot,t_m)\big) \quad \forall t_n>t_m.
	\end{align}
	These two properties, i.e., Eqs. (\ref{maximum}) and (\ref{energy}), are crucial for the ACE (\ref{ACE}). Whether the maximum principle and energy dissipation law could be preserved at the discrete level is also an important issue in the numerical methods for ACE (\ref{ACE}).

In the past years, some traditionally macroscopic numerical schemes have been developed to solve the ACE while maintaining the discrete maximum principle and energy dissipation law. Tang and Yang \cite{tang2016implicit} adopted the first-order and central difference methods for the temporal  derivative  and diffusion term, and obtained the unconditionally energy stable and maximum principle preserving linear stabilized schemes for the ACE, and the generalized ACE  with an advection term in the work \cite{shen2016maximum}. Ham and Kim \cite{ham2023stability}  analyzed the explicit time step constraints to ensure that the fully explicit finite-difference (FEX-FD) scheme for the ACE with the HN-BC can preserve the discrete maximum principle and energy dissipation law. Hou et al. \cite{hou2017numerical} proposed a second-order convex splitting scheme based on the Crank-Nicolson (CN) approach for the fractional-in-space ACE, and further discussed the discrete maximum principle and energy stability, while a nonlinear system has to be solved at each step. To avoid this problem, Hou and Leng \cite{hou2020numerical} constructed a stabilized second-order CN/Adams–Bashforth scheme for the ACE, and presented the conditions that preserve the discrete maximum principle and energy stability.  By using a modified leap-frog finite-difference scheme, Feng et al. \cite{feng2021maximum} developed a novel linear second-order finite-difference scheme for the ACE, and the maximum principle and energy dissipation law can be ensured under some constraints on the time step and coefficient of the stabilized term. In addition to above maximum principle and energy stability preserving schemes based on the single time-stepping method, Liao et al. \cite{liao2020energy} investigated the two-step backward differentiation formula with non-uniform grids for the ACE, and also discussed the discrete maximum principle and energy stability. More recently, the exponential time differencing (ETD) method \cite{du2019maximum, fu2022energy1,fu2022energy2} was proposed to naturally construct  maximum principle and energy stability preserving schemes. For instance, Du et al. \cite{du2019maximum} considered the first-order ETD and second-order ETD Runge-Kutta (ETDRK2) schemes for the non-local ACE, and demonstrated that the discrete maximum principle can be preserved unconditionally. Then Ju et al. \cite{fu2022energy1,fu2022energy2} integrated the scalar auxiliary variable (SAV) technique into stabilized first-order ETD and ETDRK2 methods, and further constructed the discrete maximum principle and energy stability preserving schemes for a class of Allen-Cahn type equations. To obtain a higher-order scheme for  the ACE, Li et al. \cite{li2020arbitrarily} developed an arbitrarily high-order multistep exponential integrator method by employing the cut-off operation \cite{cheng2022new1,cheng2022new2}  to preserve the maximum principle. Yang et al. \cite{yang2022arbitrarily} adopted the cut-off approach  and the SAV method \cite{akrivis2019energy,shen2018scalar} to develop a class of arbitrarily high-order energy stability and maximum preserving schemes for the ACE.  Subsequently, Zhang et al.  \cite{zhang2021numerical,zhang2021preserving,li2021stabilized} proposed a family of stabilized integrating factor Runge-Kutta schemes up to the third- and fourth-order in time for the ACE, which preserve the discrete maximum principle unconditionally.

 In fact, apart from the traditionally top-down macroscopic numerical methods, the mesoscopic lattice Boltzmann (LB) method, as a highly efficient second-order kinetic theory-based approach for the fluid flow problems governed by the Navier-Stokes equations (NSEs) \cite{Guo2013,Succi2008,Kruger2017,Wang2019,liu2022diffuse,liu2023improved}, can also be extended to solve the ACE (\ref{ACE}). However, it should be noted that although the LB method has a distinct advantage in terms of parallelism due to the combination of discretization and relaxation processes, the relaxation principle of the LB method comes at the cost of introducing bottom-up idea that differs from traditionally top-down discretized  methods, such as the finite-difference method \cite{tang2016implicit,shen2016maximum,hou2017numerical,hou2020numerical,feng2021maximum,ham2023stability} and the finite-element method \cite{ xiao2022second}. This fundamental feature also brings a challenge to the   numerical analysis of the LB method \cite{simonis2023lattice}.  Recently, some works have been performed to derive the macroscopic finite-difference schemes of LB methods for some specific PDEs, such as the diffusion equation \cite{Suga2010, Lin2022, silva2023discrete, Chen2023}, the convection-diffusion equation (CDE)  \cite{straka2020accuracy, chen2023fourth}, the Burgers equation \cite{li2012multilevel, chen2023cole} and also the NSEs \cite{junk2001finite}, which aim to bridge the gap between LB method and macroscopic numerical scheme, and simultaneously, the consistency and stability analysis of the LB method can be performed. However, it is important to note that the macroscopic finite-difference scheme derived from the LB method is inherently a multilevel scheme unless all diagonal entries in the diagonal relaxation matrix are fixed to be unity \cite{Be2023, Chen2023-GP}, and this would bring two problems to the LB method in terms of memory storage and initialization implementation. In order to reduce the memory requirement and improve the computational efficiency, Liu et al. \cite{Liu2023rlb} proposed an efficient second-order two-level macroscopic finite-difference scheme based on the  regularized lattice Boltzmann (RLB) method for the NSEs and CDE.  In particular, it should be noted that the macroscopic finite-difference scheme developed from the LB method in Refs. \cite{Be2023, Chen2023-GP} is at least a three-level scheme, while that in Ref. \cite{Liu2023rlb} is only a two-level scheme, which avoids the problem in the initialization implementation and also simplifies the numerical analysis. In this work, we will focus on how to develop  the RLB method based
  two-level macroscopic finite-difference scheme for the ACE (\ref{ACE}), and simultaneously,  the discrete maximum principle and energy dissipation law can be perserved.

	The rest of this paper is organized as follows. In Sec. \ref{IE-FD-RLB}, we develop a novel  RLB method based
	 macroscopic implicit-explicit finite-difference (RLB-MIE-FD)  scheme for the ACE (\ref{ACE}). In Sec. \ref{stability}, we first demonstrate that  the developed RLB-MIE-FD scheme has the first-order accuracy in time and second-order accuracy in space, then present the  conditions to ensure that the RLB-MIE-FD scheme preserves the discrete maximum principle and energy dissipation law. In Sec. \ref{Numer}, several numerical experiments are performed to confirm the theoretical analysis. Finally, some conclusions are summarized in Sec. \ref{conclusion}.

	\section{the RLB method based RLB-MIE-FD scheme for the ACE}\label{IE-FD-RLB} 
	For the sake of brevity and to simplify the following analysis, we first define the grid function space $\bm{\Phi}_{\Delta x}$ in the $d$-dimensional domain $\Omega$ as 
		\begin{align}
		\bm{\Phi}_{\Delta x}:=\bigg\{{\phi^n_\mathbf{i}}=\phi(\mathbf{x}_l+\mathbf{x}_i,t_n)\:\makebox{with}\:\mathbf{x}_l=\{L_j\}_{j=1}^d\: \makebox{is the grid function at the position} \:  \mathbf{x}_l+\mathbf{x}_i\: \makebox{and time}\: t_n\bigg\},
		\end{align} 
	where $\mathbf{x}_i$ represents $\Delta x\big(\mathbf{i}-0.5\mathbf{I}_d\big)$ with $\mathbf{I}_d=\{\underbrace{1,1,\ldots,1}_d\}$ for the HN-BC, and $\Delta x\mathbf{i}$ for the HD-BC and P-BC. It is also the node in the lattice space $\mathcal{L}=\Delta x[0,\tilde{N}]^d$ ($\tilde{N}=N$, $N$ or $N+1$ corresponds to the HD-BC, P-BC or HN-BC) with the lattice spacing $\Delta x=X/N>0$ and multi-index $\mathbf{i}\in[0,\tilde{N}]^d$ [here we only consider the square computational domain $\Omega$ with $X=R_j-L_j$ ($j=1,2,\ldots,d$) denoting the characteristic length in each direction and $N$ being a positive integer], $t_n=n\Delta t$ represents the $n_{th}$ time level with the time step $\Delta t>0$. In the LB method, the lattice velocity, defined as $c = \Delta x/\Delta t$, is also introduced. 
 
\begin{figure} [H]  
		\vspace{-0.8cm}  
		\setlength{\belowcaptionskip}{-1cm}   
		\begin{center} 
			\subfloat[D1Q3]
			{
				\includegraphics[width=0.3\textwidth]{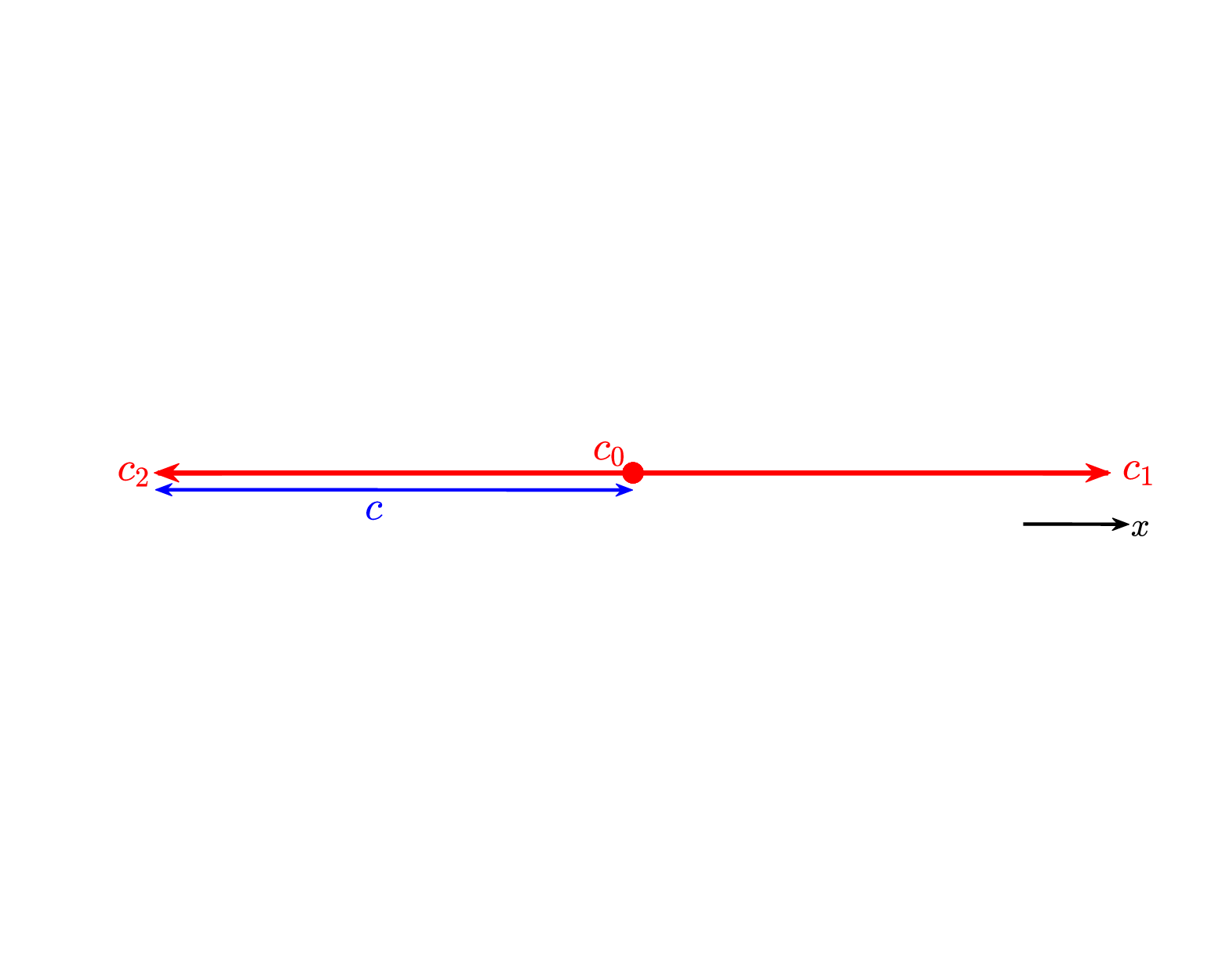}
			} 
			\subfloat[D2Q5]
			{
				\includegraphics[width=0.3\textwidth]{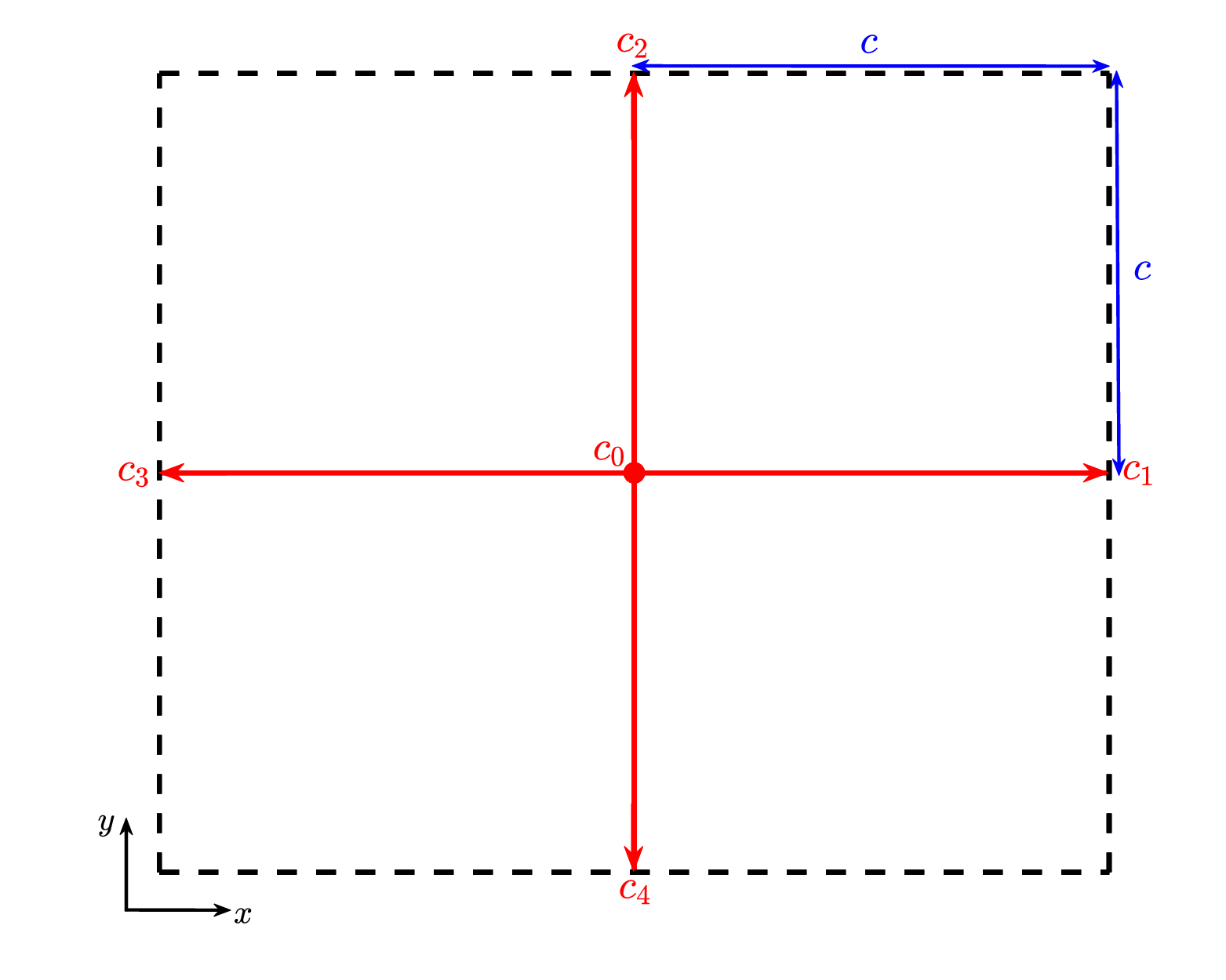}
			}
			\subfloat[D3Q7]
			{
				\includegraphics[width=0.3\textwidth]{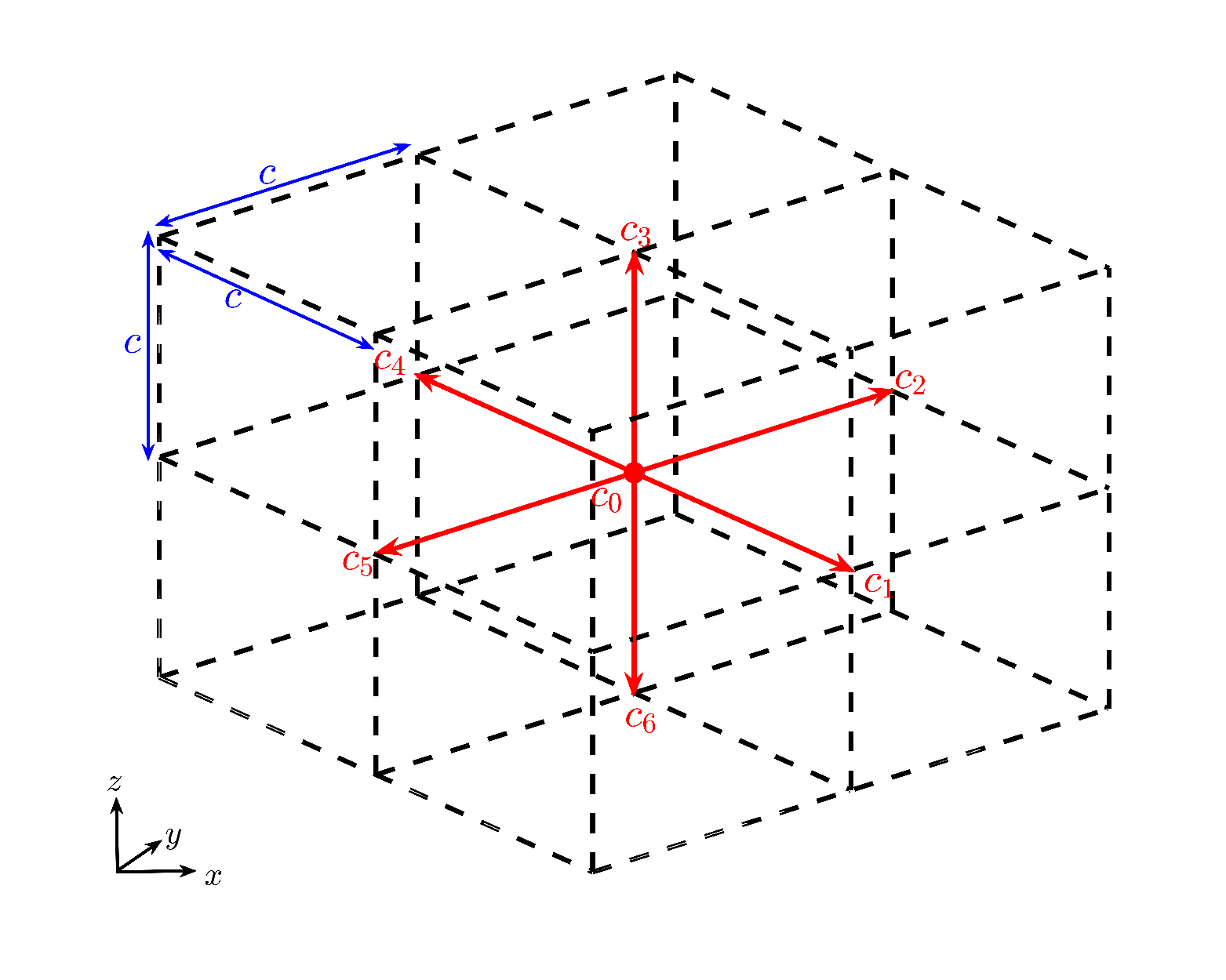}
			} 		 
			\caption{Schematic  of the D$d$Q$(2d+1)$ lattice structure.}  
			\label{Fig-DdQq}  
		\end{center}  
	\end{figure}

It is known that the discrete-velocity-Boltzmann equation (DVBE) with the D$d$Q($2d+1$) [($2d+1$) discrete velocities in $d$-dimensional space] lattice structure can be written as
	\cite{Chai2020square} 
	\begin{align}\label{dvbe}
		\partial_tf_k(\mathbf{x},t)+\mathbf{c}_k\cdot\nabla f_k(\mathbf{x},t)=-\tilde{s}f_k^{ne}(\mathbf{x},t)+F_k(\mathbf{x},t),k=0,1,2\ldots,2d,
	\end{align}
	where $f_k(\mathbf{x},t)$ is the distribution function at the position $\mathbf{x}$ and time $t$, $\tilde{s}$ is the relaxation parameter, $f_k^{ne}(\mathbf{x},t)=f_k(\mathbf{x},t)-f_k^{eq}(\mathbf{x},t)$ is the non-equilibrium distribution function, and  $f_k^{eq}(\mathbf{x},t)$ is the equilibrium distribution function. $F_k(\mathbf{x},t)$ is the distribution function of the source term. In the D$d$Q$(2d+1)$ lattice structure (see Fig. \ref{Fig-DdQq}), the discrete velocity $\mathbf{c}_i$ is given by
	\begin{align} \label{ci}
		&\mathbf{c}_0=\big(\underbrace{0,0,\ldots,0}_{d}),\mathbf{c}_i=c\mathbf{e}_i,\mathbf{c}_{i+d}=-c\mathbf{e}_i, i=1,2,\ldots,d,
	\end{align} 
	where $\mathbf{e}_i$ denotes the $i_{th}$ row of the identity matrix $\mathbf{I}_{d\times d}$.
	
	Integrating Eq. (\ref{dvbe}) along the characteristic line $\mathbf{x}'=\mathbf{x}+\mathbf{c}_kt'$ with $t'\in[0,\Delta t]$, we have
	\begin{align}\label{dvbe2}
		f_k(\mathbf{x}+\mathbf{c}_k\Delta t,t+\Delta t)=f_k(\mathbf{x},t)+\int_{0}^{\Delta t}\Big(-\tilde{s}f_k^{ne}+F_k\Big)(\mathbf{x}+\mathbf{c}_kt',t+t')dt'.
	\end{align}
	Unlike the treatment for the last term on the right-hand side of Eq. (\ref{dvbe2}) in Ref. \cite{Chai2020square}, here we consider  the following approximations: 
	\begin{subequations}\label{fk-Fk}
		\begin{align}
			&\int_0^{\Delta t}-\tilde{s}f_k^{ne}(\mathbf{x}+\mathbf{c}_kt',t+t')dt'=\frac{\Delta t}{2}\bigg[\Big(-\tilde{s}f_k^{ne}(\mathbf{x},t)\Big)+\Big(-\tilde{s}f_k^{ne}(\mathbf{x}+\mathbf{c}_k\Delta t,t+\Delta t)\Big)\bigg]+O(\Delta t^3),\\
			&\int_0^{\Delta t}F_k(\mathbf{x}+\mathbf{c}_kt',t+t')dt'= \int_0^{\Delta t}F_k\big(\mathbf{x}+\mathbf{c}_k(\Delta t+t'-\Delta t),t+t'\big)dt'\notag\\
				&\qquad\qquad\qquad\qquad\qquad \:= \int_0^{\Delta t}\Big[F_k(\mathbf{x}+\mathbf{c}_k\Delta t,t'+ t)\Big]dt' +  O(\Delta t \Delta x)\notag\\
			&\qquad\qquad\qquad\qquad\qquad \:=
			\frac{\Delta t}{2}\Big[F_k(\mathbf{x}+\mathbf{c}_k\Delta t,t)+F_k(\mathbf{x}+\mathbf{c}_k\Delta t,t+\Delta t)\Big]+O(\Delta t \Delta x+\Delta t^3).
		\end{align}
		\end{subequations}

		Based on Eq. (\ref{fk-Fk}), we can write Eq. (\ref{dvbe2}) as 
		\begin{align}\label{dvbe-1}
			\overline{f}_k(\mathbf{x}+\mathbf{c}_k\Delta t,t+\Delta t)=\tilde{f}_k(\mathbf{x},t)+\frac{\Delta t}{2}\Big[F_k(\mathbf{x}+\mathbf{c}_k\Delta t,t)+F_k(\mathbf{x}+\mathbf{c}_k\Delta t,t+\Delta t)\Big]+O(\Delta t \Delta x+\Delta t^3),
		\end{align}
		where
		\begin{subequations}\label{dvbe3}
			\begin{align}
				&\overline{f}_k=f_k-\frac{\Delta t}{2}\Big(-\tilde{s}f_k^{ne}\Big) =\Big(1+\frac{\Delta t}{2}\tilde{s}\Big)f_k-\frac{\Delta t}{2}\tilde{s}f_k^{eq},\\
				&\tilde{f}_k=f_k+\frac{\Delta t}{2}\Big(-\tilde{s}f_k^{ne}\Big) =\Big(1-\frac{\Delta t}{2}\tilde{s}\Big)f_k+\frac{\Delta t}{2}\tilde{s}f_k^{eq},
			\end{align}
		\end{subequations} 
		from above equation, we can derive
		 				\begin{subequations}\label{dvbe5}
		 	\begin{align}
		 		&f_k=\Big(1+\frac{\Delta t}{2}\tilde{s}\Big)^{-1}\Big[\overline{f}_k+\frac{\Delta t}{2}\big(\tilde{s}f_k^{eq}\big)\Big],\\
		 		&\tilde{f}_k=\Big(1-\frac{\Delta t}{2}\tilde{s}\Big) \Big(1+\frac{\Delta t}{2}\tilde{s}\Big)^{-1} \Big[\overline{f}_k+\frac{\Delta t}{2}\big(\tilde{s}f_k^{eq}\big)\Big] +\frac{\Delta t}{2}\tilde{s}f_k^{eq}.
		 	\end{align}
		 \end{subequations} 
		 Let 
		 \begin{align}
		 	s=1-\Big(1-\frac{\Delta t}{2}\tilde{s}\Big) \Big(1+\frac{\Delta t}{2}\tilde{s}\Big)^{-1}=\Big(\frac{1}{2}+\frac{1}{\Delta t\tilde{s}}\Big)^{-1},
		 \end{align}
		 then 
		 	\begin{align}\label{dvbe6}
		 		 &\tilde{f}_k=(1-s)  \overline{f}_k+(2-s)\frac{\Delta t}{2}\tilde{s}f_k^{eq} =(1-s)  \overline{f}_k+sf_k^{eq}  ,
		 	\end{align} 
		 	where $(2-s)\tilde{s}\Delta t/2=s$ has been used.

	Substituting above Eq. (\ref{dvbe6}) into Eq. (\ref{dvbe-1}) and removing the term $O(\Delta t \Delta x+\Delta t^3)$, a novel LB (N-LB) equation of the DVBE (\ref{dvbe}) can be obtained as follows:
	\begin{align}\label{lbe}
	\overline{f}_k(\mathbf{x}+\mathbf{c}_k\Delta t,t+\Delta t)&-\frac{\Delta t}{2}F_k(\mathbf{x}+\mathbf{c}_k\Delta t,t+\Delta t)=\overline{f}_k(\mathbf{x},t)-s\big[\overline{f}_k(\mathbf{x},t)-f_k^{eq}(\mathbf{x},t)\big]\notag\\
	&+\frac{\Delta t}{2}F_k(\mathbf{x}+\mathbf{c}_k\Delta t,t),k=0,1,2\ldots,2d.
	\end{align}
	to keep the notation compact, hereafter we  denote $\overline{f}_i$ as $f_i$. It should be noted that the commonly used DVBE of the SRT-LB method in the previous work \cite{Chai2020square} is written as
	\begin{align}\label{dvbe-general}
		f_k(\mathbf{x}+\mathbf{c}_k,t+\Delta t)=f_k(\mathbf{x},t)-sf_k^{ne}(\mathbf{x},t)+\Delta t\Bigg[F_k(\mathbf{x},t)+\frac{\Delta t}{2}\big(\partial_t+\gamma\mathbf{c}_k\cdot\nabla\big)F_k(\mathbf{x},t)\Bigg],
	\end{align}
	where the parameter $\gamma=0$ or $1$. In particular, if we consider $\gamma=2$ in Eq. (\ref{dvbe-general}), we have
	\begin{align}\label{gamma2}
		f_k(\mathbf{x}+\mathbf{c}_k,t+\Delta t)=f_k(\mathbf{x},t)-sf_k^{ne}(\mathbf{x},t)+\Delta t\Bigg[F_k(\mathbf{x},t)+\frac{\Delta t}{2}\big(\partial_t+2\mathbf{c}_k\cdot\nabla\big)F_k(\mathbf{x},t)\Bigg],
	\end{align}
	actually, Eq. (\ref{gamma2}) becomes the N-LB Eq. (\ref{lbe}) when the following  approximation with the truncation error   $O(\Delta t\Delta x^2)$ removed is adopted for the last term on the right-hand side of Eq. (\ref{gamma2}): 
	\begin{align}
		&\Delta t\Bigg[F_k(\mathbf{x},t)+\frac{\Delta t}{2}\big(\partial_t+2\mathbf{c}_k\cdot\nabla\big)F_k(\mathbf{x},t)\Bigg]\notag\\
		&=\frac{\Delta t}{2}\Bigg[F_k(\mathbf{x},t)+\Delta t \mathbf{c}_k\cdot\nabla F_k(\mathbf{x},t)\Bigg]+\frac{\Delta t}{2}\Bigg[F_k(\mathbf{x},t)+\Delta t \big(\partial_t+\mathbf{c}_k\cdot\nabla\big) F_k(\mathbf{x},t)\Bigg] \notag\\
		&\approx\frac{\Delta t}{2}F_k(\mathbf{x}+\mathbf{c}_k\Delta t,t)  +\frac{\Delta t}{2}F_k(\mathbf{x}+\mathbf{c}_k\Delta t,t+\Delta t).
	\end{align}

To recover the ACE (\ref{ACE}) correctly, some proper requirements on the equilibrium and source distribution functions should be satisfied. Based on the developed N-LB Eq. (\ref{lbe}) and the D$d$Q$(2d+1)$ lattice structure, the expressions of $f_k^{eq}$ and $F_k$ can be given by \cite{Chai2020square}
\begin{subequations}\label{fieq-Fi}
	\begin{align}
		&f_k^{eq}=\omega_k\phi,\\
		&F_k=\omega_k\big[-f(\phi)\big],
	\end{align}
\end{subequations}
which satisfy the following moment conditions:
\begin{subequations}\label{moment-condition}
	\begin{align}
		&\sum_{k=0}^{2d}f_k^{eq}=\phi,\sum_{k=0}^{2d}\mathbf{c}_kf_k^{eq}=\mathbf{0},\sum_{k=0}^{2d}\mathbf{c}_k\mathbf{c}_kf_k^{eq}=c_s^2\phi\mathbf{I}_{2d},\\
		&\sum_{k=0}^{2d}F_k=-f(\phi),
	\end{align}
\end{subequations}
where $c_s^2=2\omega_1c^2$	and the coefficient weights $\omega_1=\omega_2=\cdots=\omega_{2d}\in\big(0,1/(2d)\big],\omega_0=1-2d\omega_1\in[0,1)$. 

Following the idea of RLB method where a regularization process is introduced to filter out  the non-equilibrium distribution function in the developed N-LB Eq. (\ref{lbe})  \cite{wang2015regularized, Liu2023rlb}, we have
\begin{align}\label{fneq}
	f_k^{ne}(\mathbf{x},t)\approx \frac{\omega_k\mathbf{c}_k\cdot \bm{{\Pi}}^{ne}(\mathbf{x},t)}{c_s^2},
\end{align}
where  $\bm{{\Pi}}^{ne}=-\Delta tc_s^2\nabla\phi/s$ is the second-order approximation of the first-order moment of the non-equilibrium distribution functions.  

Substituting Eq. (\ref{fneq}) into Eq. (\ref{lbe}), one can obtain the following semi-discrete evolution of the RLB method for the ACE (\ref{ACE}):
\begin{align}\label{RLB}
	f_k(\mathbf{x},t+\Delta t)-\frac{\Delta t}{2}F_k(\mathbf{x},t+\Delta t)
	=&f_k^{eq}(\mathbf{x}-\mathbf{c}_k\Delta t,t)+s_A\big[\omega_k\mathbf{c}_k\cdot\bm{{\Pi}}^{ne}(\mathbf{x}-\mathbf{c}_k\Delta t,t)\big]\notag\\
	&+\frac{\Delta t}{2}F_k(\mathbf{x},t),k=0,1,2,...,2d,
\end{align}
where $s_A=(1-s)/c_s^2$.

	In order to obtain a full-discrete evolution equation of the RLB method (\ref{RLB}) for the ACE (\ref{ACE}), the gradient term in $\bm{\Pi}^{ne}$ is discretized by the following first-order difference scheme:
	\begin{subequations}\label{pine}
		\begin{align}
			&\nabla_{x_k}\phi(\mathbf{x}+\mathbf{c}_k\Delta t,t)= \frac{\phi(\mathbf{x}+\mathbf{c}_{k}\Delta t,t)-\phi(\mathbf{x},t)}{\Delta x }+O(\Delta x), k=1,2,\ldots,d,\\
			&\nabla_{x_k}\phi(\mathbf{x}+\mathbf{c}_{k+d}\Delta t,t)= \frac{\phi(\mathbf{x},t)-\phi(\mathbf{x}+\mathbf{c}_{k+d}\Delta t,t)}{\Delta x}+O(\Delta x),k=1,2,\ldots,d.
		\end{align}
	\end{subequations}
	Then, summing $f_k$ in Eq. (\ref{RLB}) over the index $k$, one can derive the  full-discrete numerical scheme, i.e., the RLB-MIE-FD scheme, which only includes the macroscopic variable $\phi$,
		\begin{align}\label{IE-FD}
			\phi(\mathbf{x},t+\Delta t)
			=&\phi(\mathbf{x},t)+\frac{\omega_1}{s} \sum_{k=1}^{d}\bigg(\phi(\mathbf{x}+\mathbf{c}_k\Delta t,t)+\phi(\mathbf{x}+\mathbf{c}_{k+d}\Delta t,t)-2\phi(\mathbf{x},t)\bigg)
			\notag\\
			&-\frac{\Delta t}{2}\Big(f(\phi)(\mathbf{x},t)+ f(\phi)(\mathbf{x},t+\Delta t)\Big),
		\end{align}
	where Eqs. (\ref{fieq-Fi}), (\ref{moment-condition}) and (\ref{pine}) have been used. In addition,  the RLB-MIE-FD scheme (\ref{IE-FD}) in the grid function space $\bm{\Phi}_{\Delta x}$ can be rewritten into the following vector form:
		\begin{align}\label{IE-FD-vector-form}
		\frac{	\bm{\Phi}^{n+1}-\bm{\Phi}^n}{\Delta t}+\frac{\big(	\bm{\Phi}^{n+1}\big)^{.3}+\big(\bm{\Phi}^n\big)^{.3}}{2}-\frac{\bm{\Phi}^{n+1}+\bm{\Phi}^n}{2}=\varepsilon^2\bm{\Lambda}_{\Delta x}\bm{\Phi}^n,
		\end{align}
	where $0\leq n\leq N_t-1$ ($N_t \Delta t=T_t$) and $\big(\bm{\Phi}^n\big)^{.3}=\Big\{\big(\phi^n_{\mathbf{i}}\big)^3\Big\}_
	{	\mathbf{i}\in [0,\tilde{N}]^d}$ with $\bm{\Phi}^n$ denoting the vector of the numerical solution at the $n_{th}$ time level, i.e.,
		\begin{subequations}
			\begin{align}
				&d=1:\:\bm{\Phi}^{n} =\Bigg[\phi^n_{0},\phi^n_{1},\ldots,\phi^n_{\tilde{N}}\Bigg]^T,\\
				&d=2:\:\bm{\Phi}^{n} =\Bigg[\phi^n_{0,0},\phi^n_{0,1},\ldots,\phi^n_{0,\tilde{N}},\phi^n_{1,0},\phi^n_{1,1},\ldots,\phi^n_{1,\tilde{N}},\ldots,\phi^n_{\tilde{N},0},\phi^n_{\tilde{N},1},\ldots,\phi^n_{\tilde{N},\tilde{N}}\Bigg]^T,\\
				&d=3:\:\bm{\Phi}^{n} =\Bigg[\phi^n_{0,0,0},\phi^n_{0,0,1},\ldots,\phi^n_{0,0,\tilde{N}},
				\phi^n_{0,1,0},\phi^n_{0,1,1},\ldots,\phi^n_{0,1,\tilde{N}},\ldots,\phi^n_{0,\tilde{N},0},\phi^n_{0,\tilde{N},1},\ldots,\phi^n_{0,\tilde{N},\tilde{N}},\notag\\
				&\qquad\qquad\qquad\:\:\phi^n_{1,0,0},\phi^n_{1,0,1},...,\phi^n_{1,0,\tilde{N}},
				\phi^n_{1,1,0},\phi^n_{1,1,1},\ldots,\phi^n_{1,1,\tilde{N}},\ldots,\phi^n_{1,\tilde{N},0},\phi^n_{1,\tilde{N},1},\ldots,\phi^n_{1,\tilde{N},\tilde{N}},
				\ldots,
				\notag\\
				&\qquad\qquad\qquad\:\:\phi^n_{\tilde{N},0,0},\phi^n_{\tilde{N},0,1},\ldots,\phi^n_{\tilde{N},0,\tilde{N}},
				\phi^n_{\tilde{N},1,0},\phi^n_{\tilde{N},1,1},\ldots,\phi^n_{\tilde{N},1,\tilde{N}},\ldots,\phi^n_{\tilde{N},\tilde{N},0},\phi^n_{\tilde{N},\tilde{N},1},\ldots,\phi^n_{\tilde{N},\tilde{N},\tilde{N}}\Bigg]^T, 
		\end{align}
	\end{subequations}
	where $[\cdot]^T$ represents the transpose of  $[\cdot]$. The diffusion coefficient $\varepsilon^2=\epsilon\Delta x^2/
	\Delta t$ with the parameter $\epsilon=\omega_1/s$, and the matrix $\bm{\Lambda}_{\Delta x}$ is defined as 
		\begin{subequations}
			\begin{align}
				&d=1:\:\bm{\Lambda}_{\Delta x}=\mathbf{D}_{\Delta x},\\
				&d=2:\:\bm{\Lambda}_{\Delta x}=\mathbf{I}_{\overline{N}\times \overline{N}}\otimes\mathbf{D}_{\Delta x}+\mathbf{D}_{\Delta x}\otimes\mathbf{I}_{\overline{N}\times \overline{N}},\\
				&d=3:\:\bm{\Lambda}_{\Delta x}=\mathbf{I}_{\overline{N}\times \overline{N}}\otimes\mathbf{I}_{\overline{N}\times \overline{N}}\otimes\mathbf{D}_{\Delta x}+\mathbf{I}_{\overline{N}\times \overline{N}}\otimes\mathbf{D}_{\Delta x}\otimes\mathbf{I}_{\overline{N}\times \overline{N}}+\mathbf{I}_{\overline{N}\times \overline{N}}\otimes\mathbf{I}_{\overline{N}\times  \overline{N}}\otimes\mathbf{D}_{\Delta x},
			\end{align}
		\end{subequations}
where the matrix $\mathbf{D}_{\Delta x}$ has the following forms under different boundary conditions, 
		\begin{subequations}\label{Ddeltax}
			\begin{align}
					&{\makebox{HN-BC}}\:(\overline{N}=N+2):\:\:\mathbf{D}_{\Delta x}=\frac{1}{\Delta x^2}\left(\begin{matrix}
						-1&1&0&0&0&\cdots&0&0&0\\
						1&-2&1&0&0&\cdots&0&0&0\\
						0&1&-2&1&0&\cdots&0&0&0\\
						\vdots&\vdots&\vdots&\vdots&\vdots&\ddots&\vdots&\vdots&\vdots\\
						0&0&0&0&0&\cdots&0&1&-1\\
						\end{matrix}\right)_{\overline{N}\times \overline{N}},\label{neumann}\\
						&{\makebox{HD-BC}}\:(\overline{N}=N-1):\:\:\mathbf{D}_{\Delta x}=\frac{1}{\Delta x^2}\left(\begin{matrix}
							-2&1&0&0&0&\cdots&0&0&0\\
							1&-2&1&0&0&\cdots&0&0&0\\
							0&1&-2&1&0&\cdots&0&0&0\\
							\vdots&\vdots&\vdots&\vdots&\vdots&\ddots&\vdots&\vdots&\vdots\\
							0&0&0&0&0&\cdots&0&1&-2\\
						\end{matrix}\right)_{\overline{N}\times \overline{N}},\label{dirichelt}\\
						&{\makebox{P-BC}}\:(\overline{N}=N+1):\:\mathbf{D}_{\Delta x}=\frac{1}{\Delta x^2}\left(\begin{matrix}
							-2&1&0&0&0&\cdots&0&0&1\\
							1&-2&1&0&0&\cdots&0&0&0\\
							0&1&-2&1&0&\cdots&0&0&0\\
							\vdots&\vdots&\vdots&\vdots&\vdots&\ddots&\vdots&\vdots&\vdots\\
							1&0&0&0&0&\cdots&0&1&-2\\\end{matrix}\right)_{\overline{N}\times \overline{N}}.\label{periodic}
			\end{align}
			\end{subequations}
	 It should be noted that all the matrices presented in Eq. (\ref{Ddeltax}) are symmetric negative semidefinite \cite{ham2023stability, lee2020numerical}.
 
	Here we would also like to point out that  under the condition of $\varepsilon^2=\epsilon\Delta x^2/\Delta t$, the RLB-MIE-FD scheme (\ref{IE-FD-vector-form})  can be   constructed directly through treating the nonlinear term $f(\phi)$ semi-implicitly $\Big[f(\phi^n_{\mathbf{i}})\approx \big(f(\phi^n_{\mathbf{i}})+f(\phi^{n+1}_{\mathbf{i}})\big)/2\Big]$ and discretizing the temporal derivative $\partial_t\phi$ and dissipative term $\nabla^2\phi$ by using the explicit  first-order Euler scheme and second-order central  difference scheme, respectively.

\section{Theoretical analysis}\label{stability}
    In this section, we will discuss the accuracy, the discrete maximum principle and energy dissipation of the RLB-MIE-FD scheme   (\ref{IE-FD-vector-form}).
    \subsection{Accuracy analysis}\label{accuracy}
   $\mathbf{{ Theorem\: 1.}}$  If the exact solution $\phi(\mathbf{x},t)$ is smooth, the proposed RLB-MIE-FD scheme (\ref{IE-FD}) has the first-order accuracy in time and second-order accuracy in space as long as $\varepsilon^2=\epsilon\Delta x^2/\Delta t$.\\
    $\mathbf{{  Proof.}}$  Applying the Taylor expansion to Eq. (\ref{IE-FD}) at the position $\mathbf{x}$ and time $t$ yields the following result:
    \begin{align}
    	\Delta t\frac{\partial\phi}{\partial t}+O(\Delta t^2)=&\epsilon\sum_{k=1}^d\Big[\Delta t\big(\mathbf{c}_k+\mathbf{c}_{k+d}\big)\cdot\nabla\phi+\frac{\Delta t^2}{2}\Big[\big(\mathbf{c}_k\cdot\nabla\big)^2+\big(\mathbf{c}_{k+d}\cdot\nabla\big)^2\Big]\phi+  O(\Delta x^4)\Big]\notag\\
    	&-\frac{\Delta t}{2}\Big[2f(\phi)+O(\Delta t)\Big].
    \end{align}
   Multiplying above equation by $1/\Delta t$, one can obtain 
        \begin{align}
    	 \frac{\partial\phi}{\partial t}=\varepsilon^2  \nabla^2\phi- f(\phi) +O(\Delta t+\Delta x^2).
    \end{align}
     The proof of this theorem is completed.

		\subsection{Discrete maximum principle}
		 $\mathbf{{ Theorem\: 2.}}$ If the initial value satisfies $\max\limits_{\mathbf{x}\in\overline{\Omega}}|\phi_0(\mathbf{x})|\leq 1$, then the full-discrete  RLB-MIE-FD scheme (\ref{IE-FD-vector-form}) preserves the discrete maximum principle provided following conditions are satisfied,
				\begin{align}\label{maximum-condition}
			\Delta t\leq 1-\frac{2d\omega_1}{s},0<\omega_1\leq\frac{1}{2d},s>1.
		\end{align}
 $\mathbf{{  Proof.}}$ We prove this theorem by induction and contradiction. First it holds that $\big\|\bm{\Phi}^0\big\|_{\infty}\leq 1$  from the assumption on $\phi_0(\mathbf{x})$ $\Big[$here we define the maximum norm as $\big\|\bm{\Phi}^n\big\|_{\infty}=\max\limits_{\mathbf{i}\in [\tilde{N}_l, \tilde{N}_r]^{d}}\big|\phi^n_{\mathbf{i}}\big|$ with $[\tilde{N}_l, \tilde{N}_r]$ denoting $[0, N + 1]$ (HN-BC), $[1, N-1]$ (HD-BC) or $[0, N]$ (P-BC)$\Big]$. Now we assume that the result holds for $n=m$, i.e., $\big\|\bm{\Phi}^m\big\|_{\infty}\leq 1$, then in the following, we will prove that  this result is also true for $n=m+1$. 
	  
   To begin with, we first rewrite the RLB-MIE-FD scheme (\ref{IE-FD-vector-form}) as follows:
	\begin{align}\label{FD-vector-2}
		 	\bm{\Phi}^{m+1}+\frac{\Delta t}{2}f(\bm{\Phi}^{m+1})=\Big[\mathbf{I}_{\overline{N}\times\overline{N}}+\epsilon{\rm diag}\big(\bm{\Lambda}\big)\Big]\bm{\Phi}^{m}-\frac{\Delta t}{2}f(\bm{\Phi}^m)+\epsilon\Big[\bm{\Lambda}-{\rm diag}\big(\bm{\Lambda}\big)\Big]\bm{\Phi}^m,
	\end{align}
	where $\bm{\Lambda}=\Delta x^2\bm{\Lambda}_{\Delta x}$ and ${\rm diag}\big(\bm{\Lambda}\big)$ is the matrix composed of the diagonal entries of the matrix $\bm{\Lambda}$. For the purpose of contradiction, we suppose $\big\|\bm{\Phi}^{m+1}\big\|_{\infty}>1$, then from   Eq. (\ref{FD-vector-2})  we have
	\begin{align}\label{proof-eq1}
		1<\Bigg\|\bm{\Phi}^{m+1}+\frac{\Delta t}{2}f(\bm{\Phi}^{m+1})\Bigg\|_{\infty}\leq
		\Bigg\|\Big[\mathbf{I}_{\overline{N}\times\overline{N}}+\epsilon{\rm diag}\big(\bm{\Lambda}\big)\Big]\bm{\Phi}^{m}-\frac{\Delta t}{2}f(\bm{\Phi}^m)\Bigg\|_{\infty}+\epsilon\Bigg\|\Big[\bm{\Lambda}-{\rm diag}\big(\bm{\Lambda}\big)\Big]\bm{\Phi}^m
	\Bigg\|_{\infty}	.
	\end{align} 
	Let $g(x)=\big(1-2d\eta\epsilon\big)x-\Delta t/2f(x)$ ($\eta=1$ or $2$), one can show that $|g(x)|\leq 1-2d\eta \epsilon$ for $x\in[-1,1]$ if $1-2d\eta\epsilon> 0$ and $\min\limits_{x\in[-1,1]}g'(x)\geq 0$ ($\Delta t\leq 1-2d\eta\epsilon)$. Therefore, if $\big\|\bm{\Phi}^{m}\big\|_{\infty}\leq 1$ and Eq. (\ref{maximum-condition}) holds, we can obtain
	\begin{align}\label{proof-eq2}
		\Bigg\|\Big[\mathbf{I}_{\overline{N}\times\overline{N}}+\epsilon{\rm diag}\big(\bm{\Lambda}\big)\Big]\bm{\Phi}^{m}-\frac{\Delta t}{2}f(\bm{\Phi}^m)\Bigg\|_{\infty}+\epsilon\Bigg\|\Big[\bm{\Lambda}-{\rm diag}\big(\bm{\Lambda}\big)\Big]\bm{\Phi}^m
		\Bigg\|_{\infty}	\leq 1,
	\end{align} 
	which is in contradiction with Eq. (\ref{proof-eq1}). Thus, the proof of this theorem
	is completed.
	
We now give a remark  on the RLB-MIE-FD scheme (\ref{IE-FD}) with the maximum principle conditions of Eq.  (\ref{maximum-condition}).\\
 $\mathbf{{  Remark\: 1.}}$ The real solution of the RLB-MIE-FD scheme (\ref{IE-FD}) at any fixed time $t_n$ ($0<n\leq N_t$) exists uniquely under conditions of Eq.  (\ref{maximum-condition}). Actually, let $\theta=\phi^{n}_{\mathbf{i}}\in[-1,1]$, $\mathbf{i}\in [\tilde{N}_l, \tilde{N}_r]^d$, and define a cubic polynomial as follows: 
\begin{align}
	p(\theta)=\theta+\frac{\Delta t}{2}f(\theta)-\xi,
\end{align}
where $\xi=\big(1-2d\epsilon\big)\phi^{n}_{\mathbf{i}}-\Delta t/2f(\phi^{n}_{\mathbf{i}})+\epsilon \sum_{k=1}^d\bigg(\phi^{n}_{\mathbf{i}+\mathbf{c}_k/{c}}+\phi^{n}_{\mathbf{i}+\mathbf{c}_{k+d}/c}\bigg)\in[-1,1]$ [see Eq. (\ref{proof-eq2})].  Since the coefficients of the odd-degree polynomial $p(\theta)$ are real, $p(\theta)$ must have real roots, and also satisfies
\begin{align}
	\begin{cases}
		p(1)=1-\xi\geq0,&\\
		p(-1)=-1-\xi \leq 0.&\\
	\end{cases}
\end{align}
Now we consider following two cases, 
\begin{align}
	&{\rm Case}\: 1: \: p(1)=1-\xi>0, 	p(-1)=-1-\xi < 0 ,\\
	&{\rm Case}\: 2: \: p(1)p(-1)=0,
\end{align} 
which indicates that the polynomial $p(\theta)$ has at least one real root and at most three real roots on $(-1,1)$ in terms of Case 1, and has the real root $1$ or $-1$ for the Case 2. Without loss of generality, we assume that one of the real roots for the two cases is $\theta_0\in[-1,1]$, and derive
\begin{align}\label{qtehta}
	\frac{2}{\Delta t}\frac{p(\theta_0)}{\theta-\theta_0}=  \theta^2+\theta_0\theta+\frac{2}{\Delta t}-1+\theta_0^2,
\end{align}
it should be noted that the discriminant of the quadratic polynomial $q(\theta)=2p(\theta_0)/\big[\Delta t(\theta-\theta_0)
\big]$  satisfies
\begin{align}\label{delta} 
	\Delta \big[q(\theta)\big]=\theta_0^2-4\bigg(\frac{2}{\Delta t}-1+\theta_0^2\bigg)<0, 
\end{align}
where the maximum principle condition (\ref{maximum-condition}) has been used.  

Thus, the polynomial $p(\theta)$ has a unique real root on $[-1,1]$, i.e.,  
\begin{align}\label{root}
	\theta_0=\bigg(\zeta+ \sqrt{\zeta^2 - \overline{\Delta} t^3}\bigg)^{1/3} +\overline{\Delta} t\bigg(\zeta + \sqrt{\zeta^2 -\overline{\Delta} t^3}\bigg)^{-1/3} ,
\end{align}
where $\zeta=\xi/\Delta t$ and $\overline{\Delta} t=(\Delta t - 2)/(3\Delta t )$. This means that ,  we do not need to numerically solve the nonlinear system at each step, which also indicates that the present RLB-MIE-FD scheme (\ref{IE-FD}) is actually an explicit scheme.
\subsection{Discrete energy dissipation law}\label{proof-energy}
Before analyzing whether the RLB-MIE-FD scheme (\ref{IE-FD-vector-form}) holds the energy dissipation law, we first present a lemma.

 $\mathbf{{ Lemma\: 1.}}$ Suppose $\mathbf{A}$ with all positive diagonal entries
is a symmetric matrix that satisfies either strictly diagonally dominant or irreducible  diagonal dominant, then $\mathbf{A}$ is a positive definite matrix  \cite{varga1976m}.

 $\mathbf{{ Theorem\: 3.}}$ Define the discrete energy as
\begin{align}\label{proof2-eq0}
	E_{\Delta x}(\bm{\Phi})=\Delta x^d\Bigg[\frac{1}{4}\sum_{\mathbf{i}\in[\tilde{N}_l,\tilde{N}_r]^d}\Big(\big(\phi_{\mathbf{i}}\big)^2-1\Big)^2-\frac{\epsilon}{2\Delta t}\bm{\Phi}^T\bm{\Lambda}\bm{\Phi}\Bigg].
\end{align}
Under the conditions of Eq. (\ref{maximum-condition}), the RLB-MIE-FD scheme (\ref{IE-FD-vector-form}) satisfies the discrete energy dissipation law, i.e.,
\begin{align}\label{discrete-energy-law}
	E_{\Delta x}\big(\bm{\Phi}^{n+1}\big)-E_{\Delta x}\big(\bm{\Phi}^{n}\big)\leq0, n=0,1,\ldots,N_t-1.
\end{align}
 $\mathbf{{  Proof.}}$ On one hand, taking the difference of the discrete energy between two successive time levels $n$ and $n+1$ yields
\begin{align}\label{proof2-eq1}
	E_{\Delta x}\big(\bm{\Phi}^{n+1}\big)-E_{\Delta x}\big(\bm{\Phi}^{n}\big)=&\frac{\Delta x^d}{4}\sum_{\mathbf{i}\in[\tilde{N}_l,\tilde{N}_r]^d}\Bigg[\Big(\big(\phi^{n+1}_{\mathbf{i}}\big)^2-1\Big)^2- \Big(\big(\phi^{n}_{\mathbf{i}}\big)^2-1\Big)^2 \Bigg]-\frac{\Delta x^d\epsilon}{2\Delta t}\Big(\big(\bm{\Phi}^{n+1}\big)^T\bm{\Lambda}\bm{\Phi}^{n+1}-\big(\bm{\Phi}^n\big)^T\bm{\Lambda}\bm{\Phi}^n\Big)\notag\\
	=&\frac{\Delta x^d}{4}	  \sum_{\mathbf{i}\in[\tilde{N}_l,\tilde{N}_r]^d}\Big[\big(\phi^{n+1}_{\mathbf{i}}\big)^3+\big(\phi^{n}_{\mathbf{i}}\big)^3+\phi^{n+1}_{\mathbf{i}}\big(\phi^{n}_{\mathbf{i}}\big)^2+\big(\phi^{n+1}_{\mathbf{i}}\big)^2\phi^{n}_{\mathbf{i}}-2\big(\phi^{n}_{\mathbf{i}}+\phi^{n+1}_{\mathbf{i}}\big)\Big]\big(\phi^{n+1}_{\mathbf{i}}-\phi^{n}_{\mathbf{i}}\big)\notag\\
	&-\frac{\Delta x^d\epsilon}{2\Delta t}\big(\bm{\Phi}^{n+1}-\bm{\Phi}^n\big)^T\bm{\Lambda}\big(\bm{\Phi}^{n+1}+\bm{\Phi}^n\big),
\end{align}
where the symmetry of the matrix $\bm{\Lambda}$ has been used [see Eq. (\ref{Ddeltax})].

On the other hand, taking $L^2$ inner product of the RLB-MIE-FD scheme (\ref{IE-FD-vector-form}) with $\Delta x^d\big(\bm{\Phi}^{n+1}-\bm{\Phi}^n\big)$, we have
\begin{align}\label{proof2-eq2}
	\frac{\Delta x^d}{4}	&\sum_{\mathbf{i}\in[\tilde{N}_l,\tilde{N}_r]^d}\Big[\big(\phi^{n+1}_{\mathbf{i}}\big)^3+\big(\phi^{n}_{\mathbf{i}}\big)^3+\phi^{n+1}_{\mathbf{i}}\big(\phi^{n}_{\mathbf{i}}\big)^2+\big(\phi^{n+1}_{\mathbf{i}}\big)^2\phi^{n}_{\mathbf{i}}-2\big(\phi^{n}_{\mathbf{i}}+\phi^{n+1}_{\mathbf{i}}\big)\Big]\big(\phi^{n+1}_{\mathbf{i}}-\phi^{n}_{\mathbf{i}}\big)\notag\\
	&+\frac{\Delta x^d}{\Delta  t}\sum_{\mathbf{i}\in[\tilde{N}_l,\tilde{N}_r]^d}\big(\phi^{n+1}_{\mathbf{i}}-\phi^{n}_{\mathbf{i}}\big)^2+\frac{\Delta x^d}{4}\sum_{\mathbf{i}\in[\tilde{N}_l,\tilde{N}_r]^d}\Big[\big(\phi^{n+1}_{\mathbf{i}}\big)^3+\big(\phi^{n}_{\mathbf{i}}\big)^3-\phi^{n+1}_{\mathbf{i}}\big(\phi^{n}_{\mathbf{i}}\big)^2-\big(\phi^{n+1}_{\mathbf{i}}\big)^2\phi^{n}_{\mathbf{i}}\Big]\big(\phi^{n+1}_{\mathbf{i}}-\phi^{n}_{\mathbf{i}}\big)\notag\\
	&+\frac{\Delta x^d\epsilon}{2\Delta t}\big(\bm{\Phi}^{n+1}-\bm{\Phi}^n\big)^T\bm{\Lambda}\big(\bm{\Phi}^{n+1}-\bm{\Phi}^{n}\big)=\frac{\Delta x^d\epsilon}{2\Delta t}\big(\bm{\Phi}^{n+1}-\bm{\Phi}^n\big)^T\bm{\Lambda}\big(\bm{\Phi}^{n+1}+\bm{\Phi}^n\big).
\end{align}
Replacing the last term on the right-hand side of Eq. (\ref{proof2-eq1}) with the term on the left-hand side of Eq. (\ref{proof2-eq2}), we can obtain 
\begin{align}\label{proof2-eq3}
	E_{\Delta x}\big(\bm{\Phi}^{n+1}\big)-E_{\Delta x}\big(\bm{\Phi}^{n}\big)=&-\frac{\Delta x^d}{4}\sum_{\mathbf{i}\in[\tilde{N}_l,\tilde{N}_r]^d}\Big[\big(\phi^{n+1}_{\mathbf{i}}\big)^3+\big(\phi^{n}_{\mathbf{i}}\big)^3-\phi^{n+1}_{\mathbf{i}}\big(\phi^{n}_{\mathbf{i}}\big)^2-\big(\phi^{n+1}_{\mathbf{i}}\big)^2\phi^{n}_{\mathbf{i}}\Big]\big(\phi^{n+1}_{\mathbf{i}}-\phi^{n}_{\mathbf{i}}\big)\notag\\
	&-\frac{\Delta x^d}{\Delta t}\sum_{\mathbf{i}\in[\tilde{N}_l,\tilde{N}_r]^d}\big(\phi^{n+1}_{\mathbf{i}}-\phi^{n}_{\mathbf{i}}\big)^2-\frac{\Delta x^d\epsilon}{2\Delta t}\big(\bm{\Phi}^{n+1}-\bm{\Phi}^n\big)^T\bm{\Lambda}\big(\bm{\Phi}^{n+1}-\bm{\Phi}^{n}\big).
\end{align}
From   Eq. (\ref{maximum-condition}) in Theorem 2, we have $\big\|\bm{\Phi}^{l}\big\|_{\infty}\leq 1$ ($l=n,n+1$). Consequently, we can derive
\begin{align}\label{proof2-eq4}
	-\frac{1}{4}&\sum_{\mathbf{i}\in[\tilde{N}_l,\tilde{N}_r]^d}\Big[\big(\phi^{n+1}_{\mathbf{i}}\big)^3+\big(\phi^{n}_{\mathbf{i}}\big)^3-\phi^{n+1}_{\mathbf{i}}\big(\phi^{n}_{\mathbf{i}}\big)^2-\big(\phi^{n+1}_{\mathbf{i}}\big)^2\phi^{n}_{\mathbf{i}}\Big]\big(\phi^{n+1}_{\mathbf{i}}-\phi^{n}_{\mathbf{i}}\big)\notag\\
	&=-\frac{1}{4}\sum_{\mathbf{i}\in[\tilde{N}_l,\tilde{N}_r]^d}
	\Big[\big(\phi^{n+1}_{\mathbf{i}}\big)^2-\big(\phi^{n}_{\mathbf{i}}\big)^2\Big]  \big(\phi^{n+1}_{\mathbf{i}}-\phi^{n}_{\mathbf{i}}\big)^2\leq \frac{1}{4}\sum_{\mathbf{i}\in[\tilde{N}_l,\tilde{N}_r]^d}\big(\phi^{n+1}_{\mathbf{i}}-\phi^{n}_{\mathbf{i}}\big)^2.
\end{align}
Combining Eq. (\ref{proof2-eq4}) with Eq. (\ref{proof2-eq3}), we have 
\begin{align}\label{proof2-eq5}
	E_{\Delta x}\big(\bm{\Phi}^{n+1}\big)-E_{\Delta x}\big(\bm{\Phi}^{n}\big)&\leq \Delta x^d\big(\bm{\Phi}^{n+1}-\bm{\Phi}^n\big)^T\Bigg[\frac{\mathbf{I}_{\overline{N}^d\times \overline{N}^d}}{4}-\frac{\mathbf{I}_{\overline{N}^d\times \overline{N}^d}}{\Delta t}-\frac{\epsilon\bm{\Lambda}}{2\Delta t}\Bigg]\big(\bm{\Phi}^{n+1}-\bm{\Phi}^{n}\big)\notag\\
	&=-\frac{\Delta x^d}{4\Delta t} \big(\bm{\Phi}^{n+1}-\bm{\Phi}^n\big)^T\Bigg[-\Delta t\mathbf{I}_{\overline{N}^d\times \overline{N}^d}+4\mathbf{I}_{\overline{N}^d\times \overline{N}^d}+2\epsilon\bm{\Lambda}\Bigg]\big(\bm{\Phi}^{n+1}-\bm{\Phi}^{n}\big).
\end{align}
Actually, one can show that the matrix $\bigg(-\Delta t\mathbf{I}_{\overline{N}^d\times \overline{N}^d} +4\mathbf{I}_{\overline{N}^d\times \overline{N}^d}+2\epsilon\bm{\Lambda}\bigg)$ in Eq. (\ref{proof2-eq5}) is a positive definite matrix (see Lemma 1) once the following conditions are staisfied for a specified boundary condition,
\begin{subequations}
	\begin{align}	
		&{\makebox{HN-BC:}}\:\begin{cases}
			4-\Delta t-4d\epsilon>0\\ 
			|4-\Delta t-4d\epsilon|\geq 4d\epsilon
		\end{cases}
		\Rightarrow  
		\Delta t\leq 4\big(1-2d\epsilon\big)  ,\\
		&{\makebox{HD-BC:}}\:\begin{cases}
			4-\Delta t-4d\epsilon>0\\
			|4-\Delta t-4d\epsilon|> 4d\epsilon
		\end{cases}
		\Rightarrow  
		\Delta t< 4\big(1-2d\epsilon\big)  ,\\
		&{\makebox{P-BC:}}\:\begin{cases}
			4-\Delta t-4d\epsilon>0\\
			|4-\Delta t-4d\epsilon|\geq 4d\epsilon
		\end{cases}
		\Rightarrow  
		\Delta t\leq 4\big(1-2d\epsilon\big).
	\end{align}
\end{subequations}
Due to the fact   $\Delta t\leq 1-2d\epsilon$ with $0<\epsilon < 1/(2d)$, i.e., the maximum principle conditions of Eq. (\ref{maximum-condition}),  the result of the Eq. (\ref{discrete-energy-law}) is proved.

For the maximum principle and energy dissipation law conditions of the RLB-MIE-FD scheme (\ref{IE-FD-vector-form}) shown in Eq.  (\ref{maximum-condition}), we give a remark.\\
 $\mathbf{{ Remark\: 2.}}$ The conditions to ensure that the FEX-FD scheme developed in Ref. \cite{ham2023stability} can preserve the maximum principle and energy dissipation law, i.e., $\Delta t\leq (1-2d\epsilon)/2$ with $0<\epsilon\leq 1/(2d)$, are actually the sufficient conditions for Eq. (\ref{maximum-condition}), which this also means that the RLB-MIE-FD scheme (\ref{IE-FD}) would be more stable than the FEX-FD scheme (see Part. \ref{Part-ex1}).
\section{Numerical results and discussion}\label{Numer}
In this section, we will carry out some numerical examples to validate the theoretical results of the present RLB-MIE-FD scheme (\ref{IE-FD-vector-form}), the relaxation parameter $s$ is chosen based on the maximum principle and energy dissipation conditions [see Eq. (\ref{maximum-condition})],  and for simplicity, the weight coefficient is set to be $\omega_1=1/3$ $(d=1)$ $\big[1/5$ $(d=2)$ and $1/6$ $(d=3)\big]$. To measure the accuracy of the RLB-MIE-FD scheme (\ref{IE-FD-vector-form}), we consider the numerical error at time $t=T_t$ in both the maximum and $l^2$ norms, i.e.,
\begin{subequations}\label{error}
	\begin{align}
		& \makebox{Err}(\bm{\Psi},\infty)=\Big\|\bm{\Phi} -\bm{\Psi}(\mathbf{x}_i,t=T_t)\Big\|_{\infty},\\
		& \makebox{Err}(\bm{\Psi},2)=\sqrt{\Delta x^d\sum_{\mathbf{i}\in[0,N]^d}  \Big(\phi_{\mathbf{i}}-\psi(\mathbf{x}_i,t=T_t)\Big)^2},
	\end{align}
\end{subequations}
where $\bm{\Phi}$ and $\bm{\Psi}$ are the vectors of the numerical result and analytical solution $\phi^{\star}$ (or reference solution $\overline{\phi}$), respectively. According to  the  Eq. (\ref{error}), one can estimate the convergence rate (CR) of the numerical scheme (\ref{IE-FD-vector-form}) with the following formulas:
\begin{subequations}
	\begin{align}
		&	 	\makebox{CR}(\bm{\Psi},\infty)=\frac{\log\Big(\makebox{Err}(\bm{\Psi},\infty)_{\Delta x}/\makebox{Err}(\bm{\Psi},\infty)_{\Delta x/2}\Big)}{\log2},\\
		& 	\makebox{CR}(\bm{\Psi},2)=\frac{\log\Big(\makebox{Err}(\bm{\Psi},2)_{\Delta x}/\makebox{Err}(\bm{\Psi},2)_{\Delta x/2}\Big)}{\log2}.
	\end{align}
\end{subequations}

\subsection{The one-dimensioanl ACE}
\subsubsection{Traveling wave solution}\label{Part-ex1}
One of the exact solutions of the ACE (\ref{ACE}) is the traveling wave solution \cite{ham2023stability}:
\begin{align} 
	\phi^{\star}(x,t)=\frac{1}{2}\Bigg[1-\tanh\bigg(\frac{x-s_{\varepsilon}t}{2\sqrt{2}\varepsilon}\bigg)\Bigg],x\in \Omega,
\end{align}
where the computational  domain is $\Omega=(-0.5,2.5)$ with the HN-BC, and $s_{\varepsilon}=3\varepsilon/\sqrt{2}$ is the speed of the traveling wave.

\begin{figure}[H]
	\vspace{-0.2cm}  
	\setlength{\belowcaptionskip}{-1cm}   
	\begin{center}     
		\includegraphics[width=0.4\textwidth]{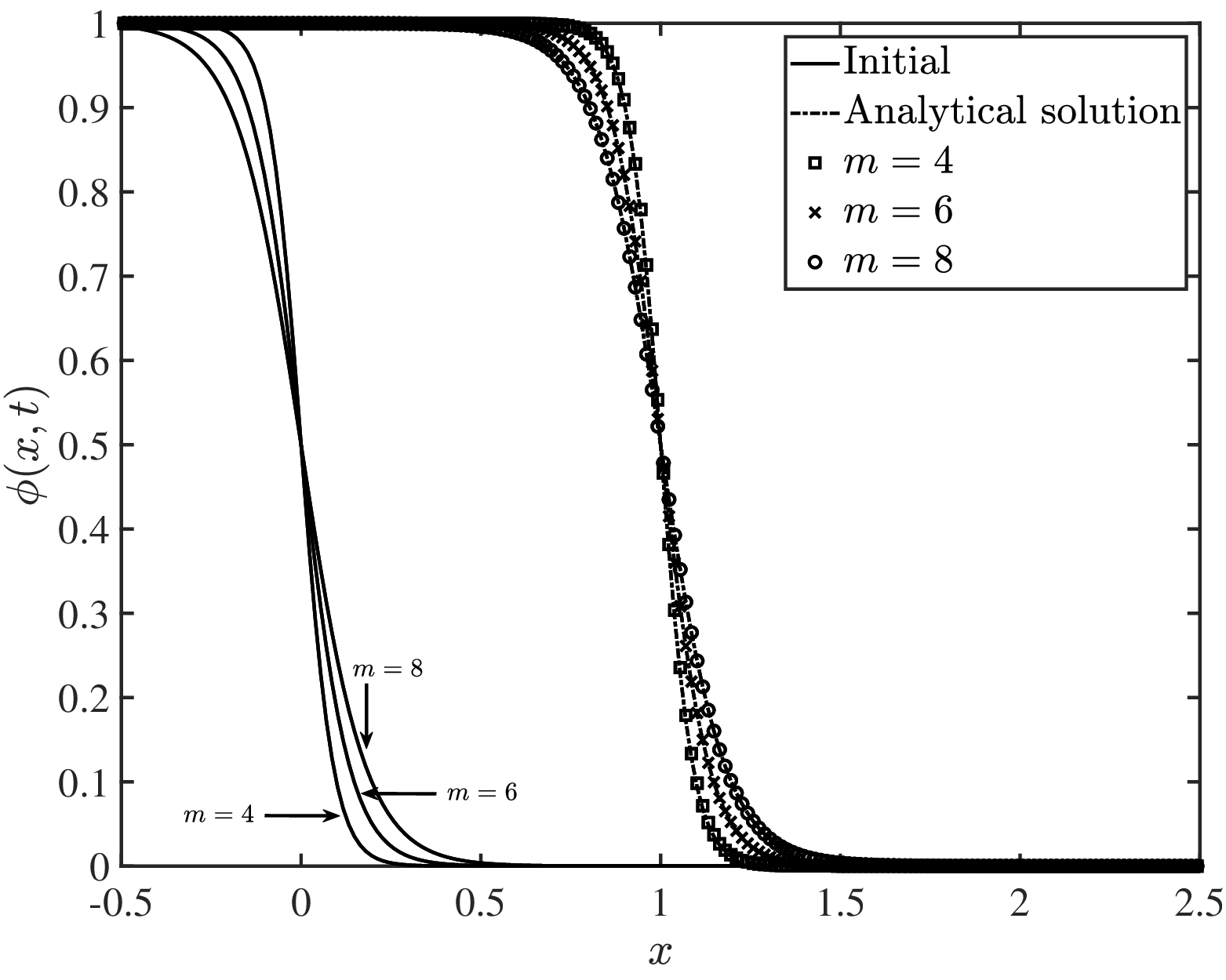} 
		\caption{ The initial, numerical and analytical solutions under different values of interfacial parameter  $\varepsilon$ (symbols: numerical solutions).}   
		\label{Fig-Ex1-1} 
	\end{center}  
\end{figure}

\begin{figure} [H]
	\vspace{-0.4cm}  
	\setlength{\belowcaptionskip}{-1cm}   
	\begin{center} 
		\subfloat[$\Delta t=5\Delta x^2$]
		{
			\includegraphics[width=0.33\textwidth]{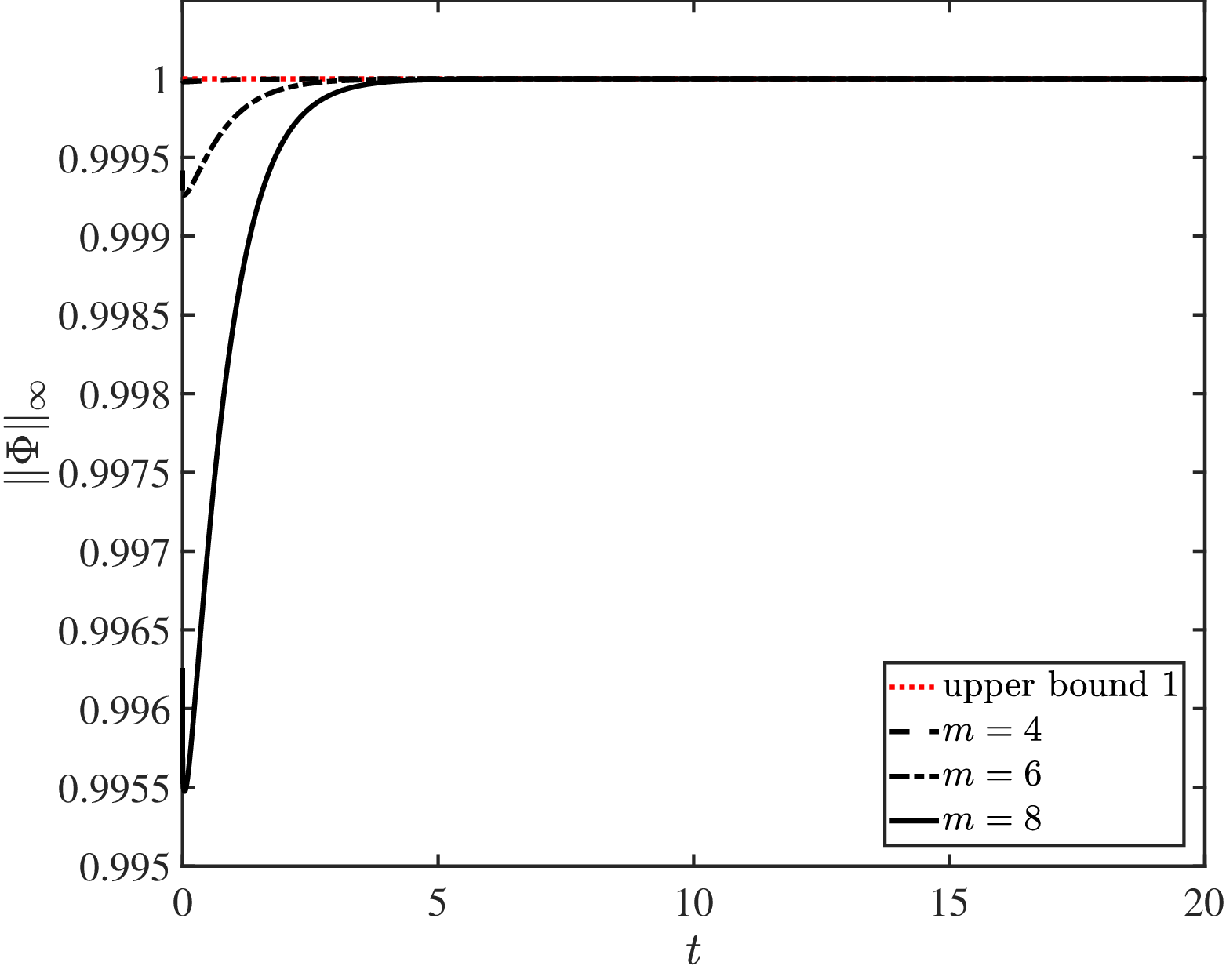}
		} 
		\subfloat[$\Delta t=1/5$]
		{
			\includegraphics[width=0.33\textwidth]{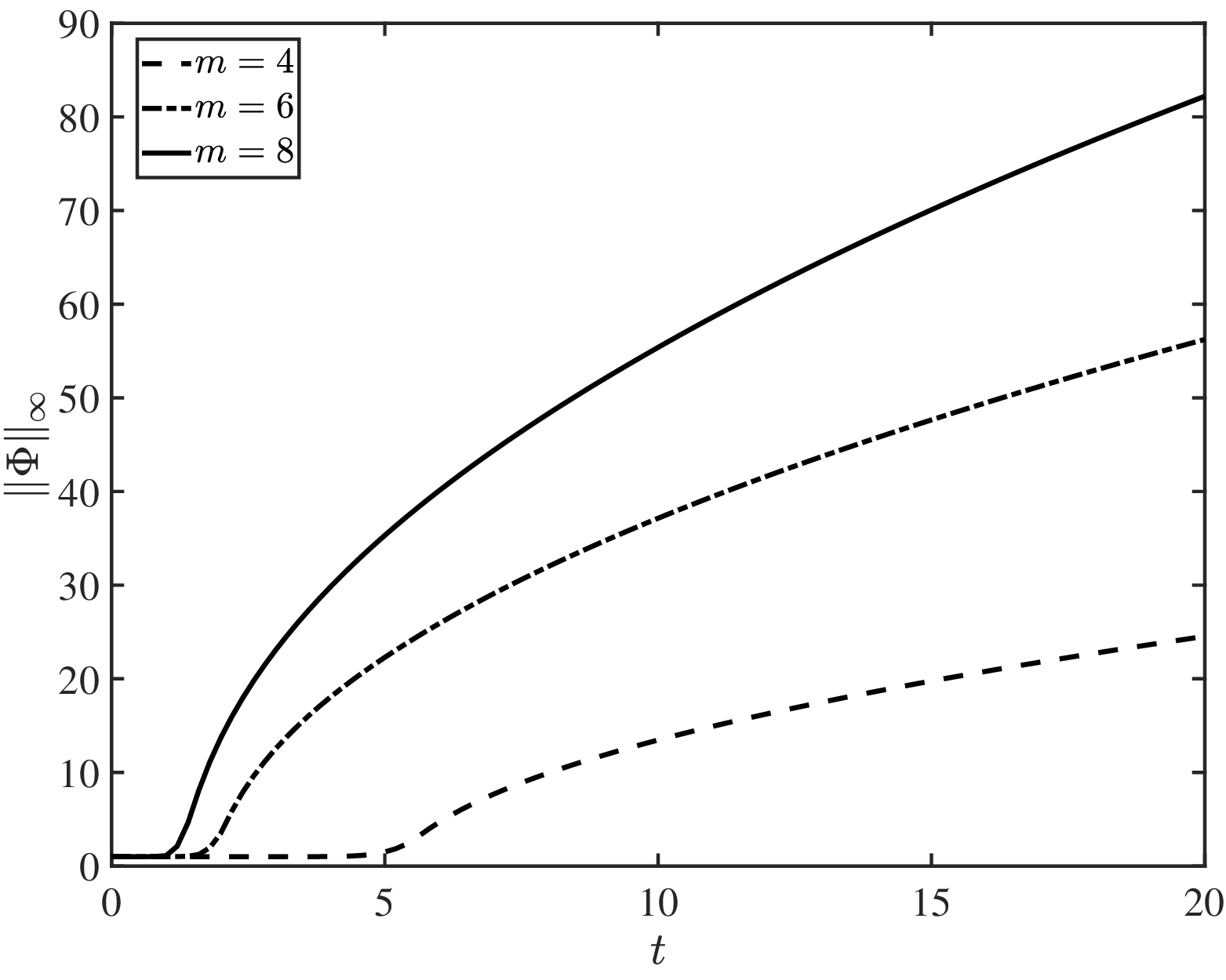}
		}
		
		\subfloat[$\Delta t=5\Delta x^2$]
		{
			\includegraphics[width=0.33\textwidth]{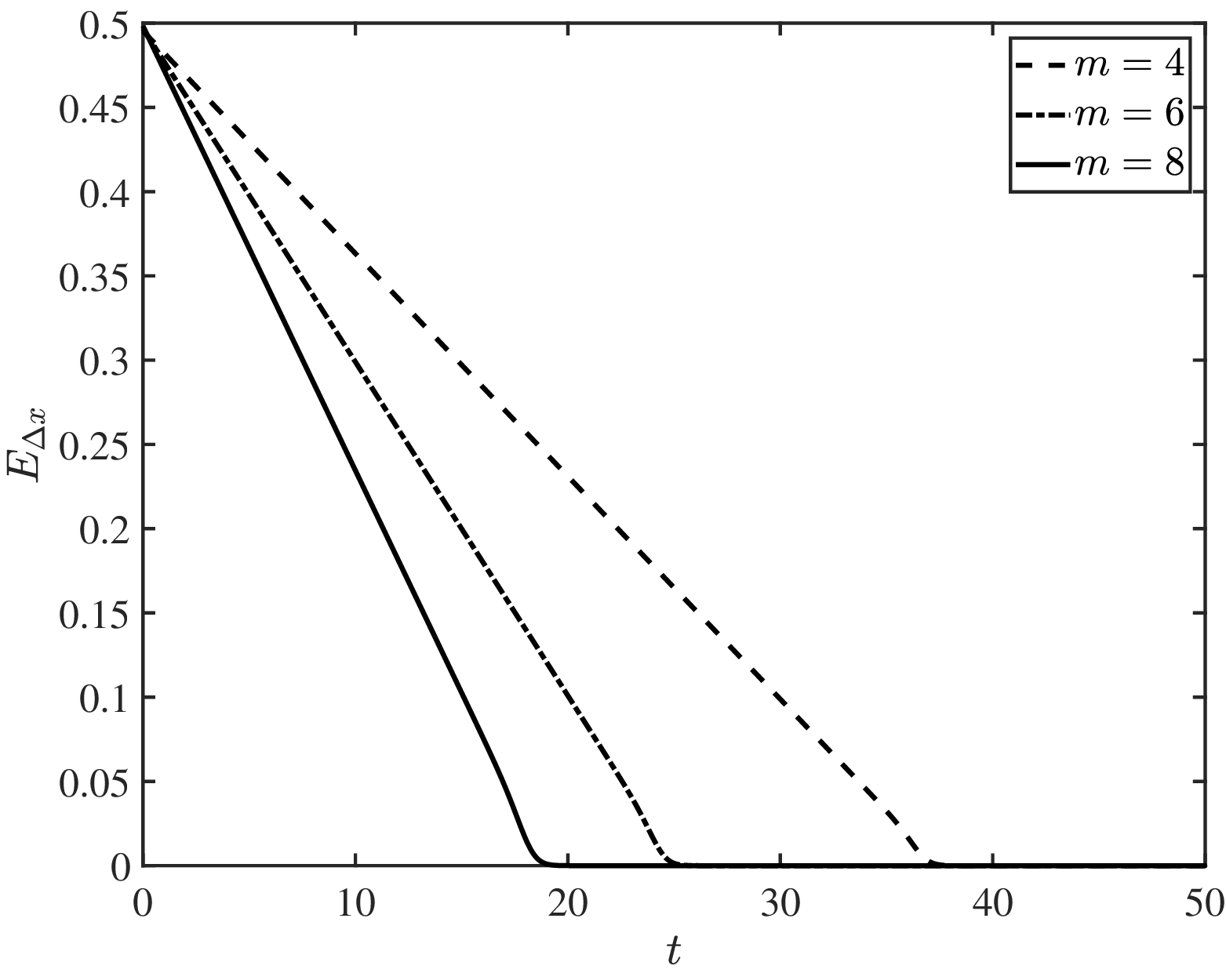}
		}		
		\subfloat[$\Delta t=1/5$ ]
		{
			\includegraphics[width=0.33\textwidth]{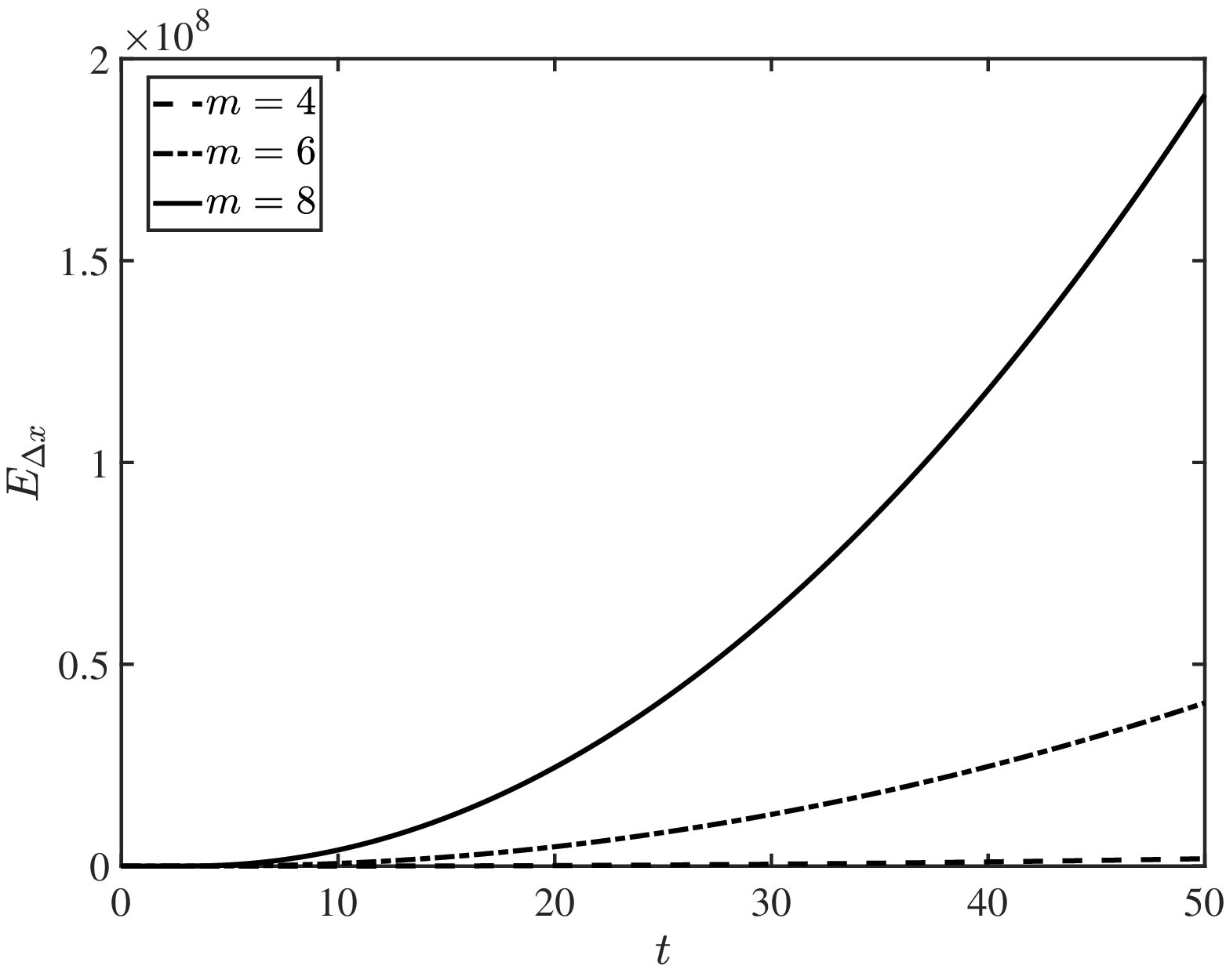}
		}  			 
		\caption{ The discrete maximum norms (a-b) and energies (c-d) of the numerical solutions under different values of interfacial parameter  $\varepsilon$.}  
		\label{Fig-Ex1}  
	\end{center}  
\end{figure}

\begin{figure} [H]
	\vspace{-0.4cm}  
	\setlength{\belowcaptionskip}{-1cm}   
	\begin{center} 
		\subfloat[$\Delta t=1/5$]
		{
			\includegraphics[width=0.33\textwidth]{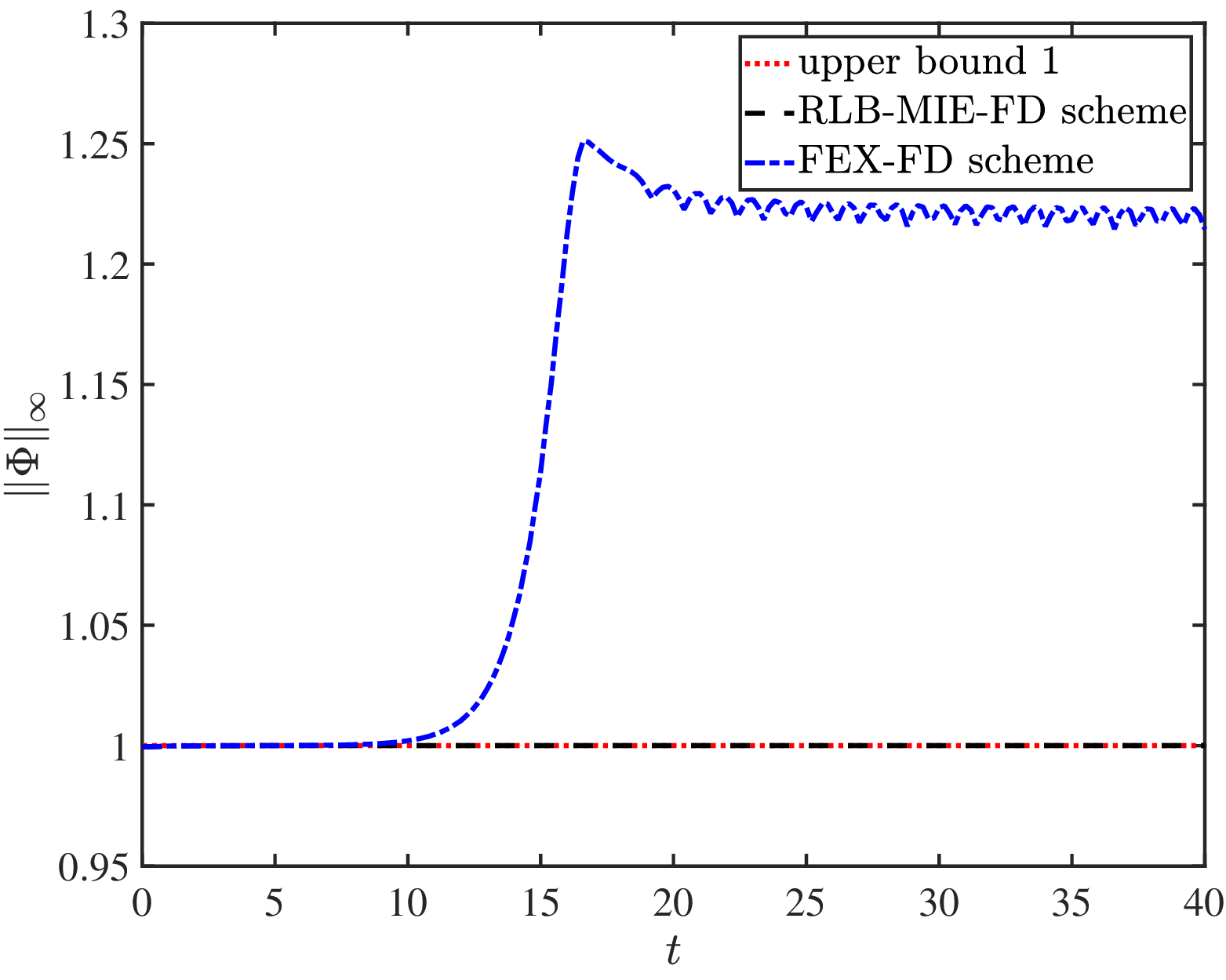}
		} 
		\subfloat[$\Delta t=1/10$]
		{
			\includegraphics[width=0.33\textwidth]{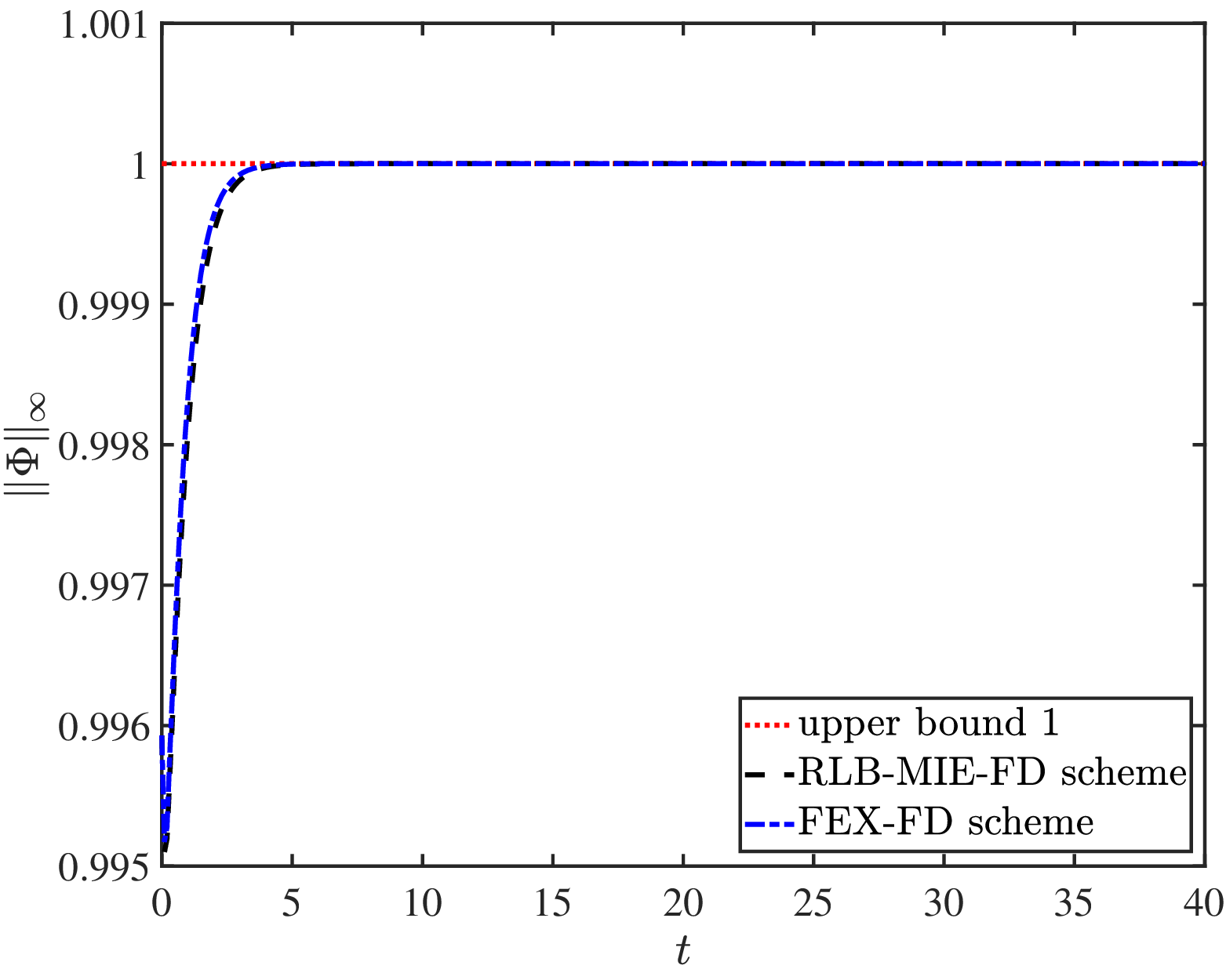}
		}   	
		
		\subfloat[$\Delta t=1/5$]
		{
			\includegraphics[width=0.33\textwidth]{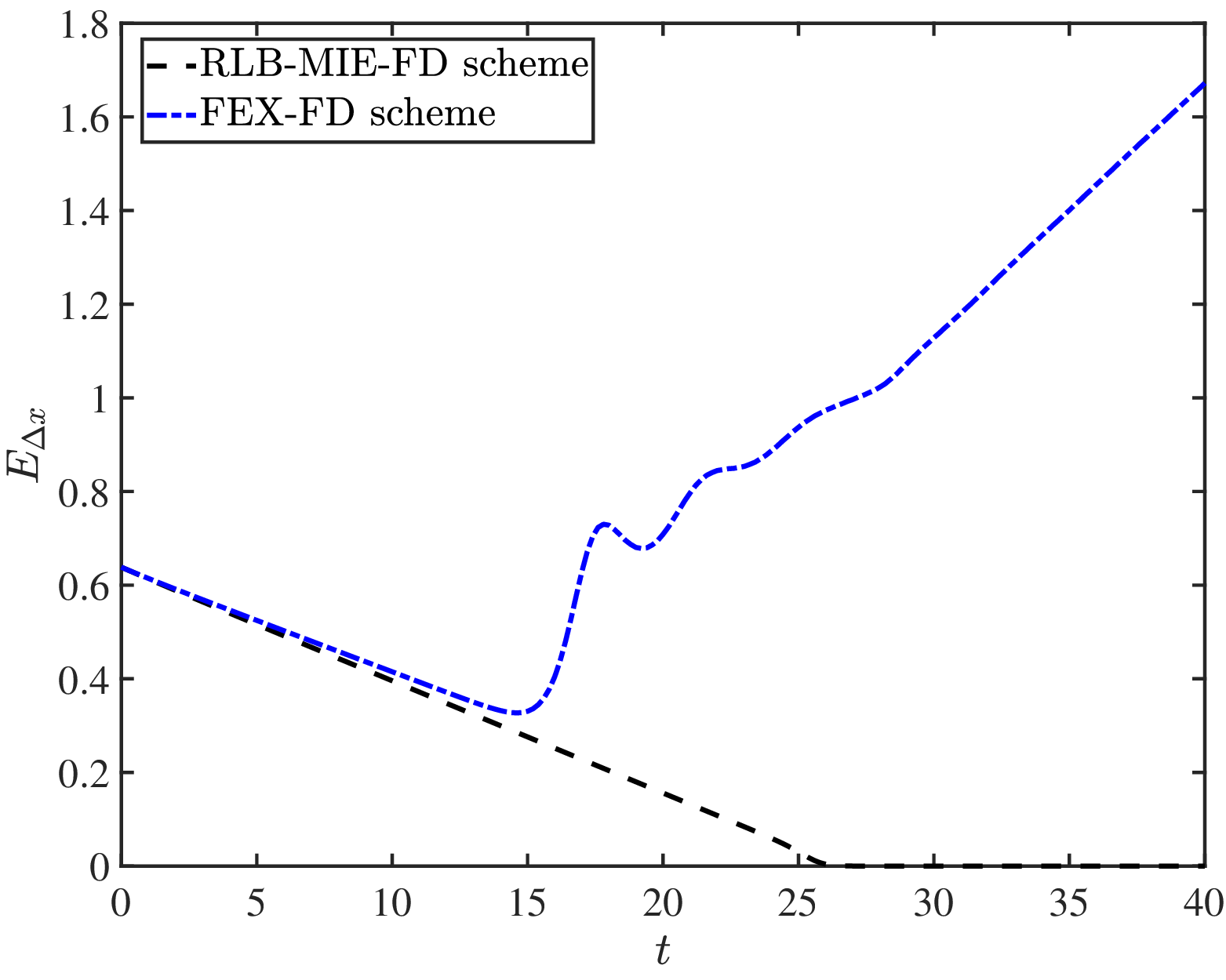}
		} 
		\subfloat[$\Delta t=1/10$]
		{
			\includegraphics[width=0.33\textwidth]{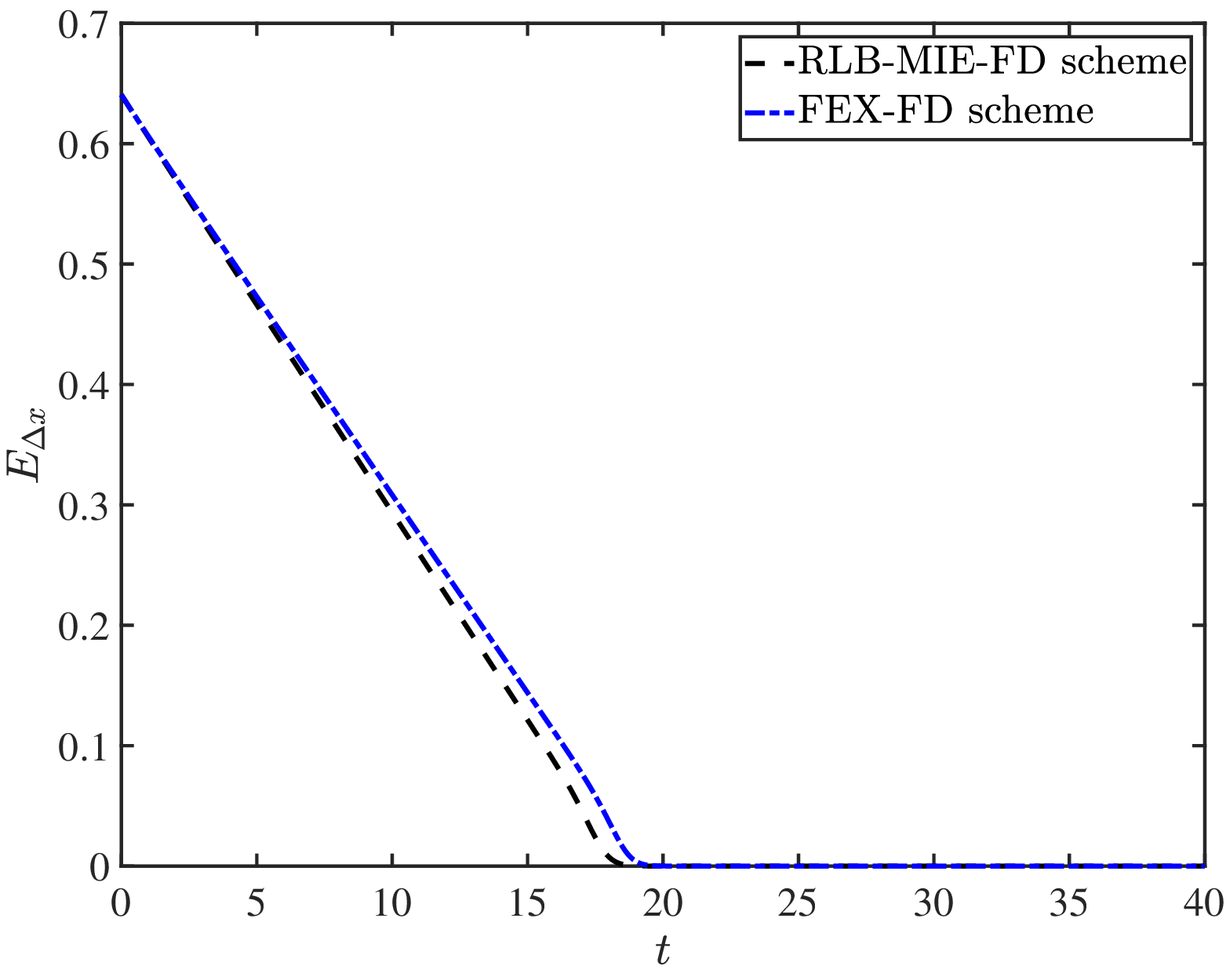}
		}   			 
		\caption{The discrete maximum norms (a-b) and energies (c-d) of the numerical solutions of the RLB-MIE-FD and FEX-FD schemes.}  
		\label{Fig-Ex1-2}  
	\end{center}  
\end{figure}

We now test the performance of the RLB-MIE-FD scheme, and consider three different values of interfacial parameter $\varepsilon=m\tanh(0.9)\Delta x/\sqrt{2}$ ($m=4,6$ and $8$) with the lattice spacing $\Delta x=1/64$ and time step $\Delta t=5\Delta x^2$. We plot the numerical and analytical results at time $T_t=1/s_{\varepsilon}$ in  Fig. \ref{Fig-Ex1-1}, and from this figure, one can find that the numerical solutions of the RLB-MIE-FD scheme are in good agreement with the analytical solutions. 
 Additionally, the evolutions of the corresponding discrete maximum norms and energies of the numerical solutions at different values of the time step $\Delta t$ are shown in Figs. \ref{Fig-Ex1}(a-b) and \ref{Fig-Ex1}(c-d), respectively. As seen from these figures, the  maximum principle and energy dissipation law are well preserved numerically 
when the time step is given by $\Delta t=5\Delta x^2$ $\big[\epsilon(m):=0.0050(4),0.0113(6),0.0200(8)\big]$, while not when $\Delta t=1/5$ $\big[\epsilon(m): =0.8209(4),1.8471(6),3.2837(8)\big]$, these results are consistent with the theoretical results in Theorems 2 and 3. 
To compare the numerical stability of the RLB-MIE-FD and FEX-FD schemes, we plot the evolutions of the discrete maximum norms and energies of the numerical solutions of the two schemes at different values of time step $\Delta t$ in Fig. \ref{Fig-Ex1-2}, where the parameter $\epsilon=0.45$ and lattice spacing $\Delta x=1/32$. Also, the corresponding numerical errors in the maximum and $l^2$ norms of the two schemes are given in Table \ref{Table-Ex1-0}. From Figs. \ref{Fig-Ex1-2}(a) and \ref{Fig-Ex1-2}(c), one can find that the RLB-MIE-FD scheme preserves the discrete maximum principle and energy dissipation law well while the FEX-FD scheme does not, which is consistent with Remark 2. In particular, the RLB-MIE-FD scheme is more accurate than the FEX-FD scheme (see the last row of Table \ref{Table-Ex1-0}). Furthermore, it should be noted that from Table \ref{Table-Ex1-0}, the FEX-FD scheme does not have a convergence behavior while the RLB-MIE-FD scheme does, although Figs. \ref{Fig-Ex1-2}(b) and \ref{Fig-Ex1-2}(d) demonstrate that the discrete maximum principle and energy dissipation law of the FEX-FD scheme can be preserved. From these results, we can conclude that the present RLB-MIE-FD scheme is more  stable than the FEX-FD scheme developed in Ref. \cite{ham2023stability}.
\begin{table}[H]
	\begin{center} 
		\caption{The numerical errors   of the RLB-MIE-FD and FEX-FD schemes at different values of time step $\Delta t$.} 	
		\label{Table-Ex1-0}
		\begin{tabular}{ccccccccccccc }\hline\hline
		\multirow{2}{*}{$\Delta t$}& \multicolumn{2}{c}{RLB-MIE-FD scheme}&&\multicolumn{2}{c}{FEX-FD scheme}\\
		\cline{2-3}\cline{5-6}
				&Err$(\bm{\Phi}^{\star},\infty)$& Err$(\bm{\Phi}^{\star},2)$&&Err$(\bm{\Phi}^{\star},\infty)$& Err$(\bm{\Phi}^{\star},2)$&  \\
			\hline
			$1/5$&  7.0684$\times 10^{-10}$ &5.9314$\times 10^{-10}$&&59.7505&6.4623  \\
			$1/10$&1.5710$\times 10^{-13}$&2.8125$\times 10^{-14}$&&2.0346$\times 10^{-5}$ &1.7359$\times 10^{-6}$\\ 
			\hline\hline
		\end{tabular}
	\end{center}  
\end{table} 
  Finally, we consider the interfacial parameter $\varepsilon=\tanh(0.9)/\big(16\sqrt{2}\big)$ and test the CR in space at time $T_t=1/s_{\varepsilon}$. The RLB-MIE-FD scheme is adopted to calculate the
   numerical solutions of the ACE (\ref{ACE}) at different lattice spacings $\Delta x=2^{-k}/64$ and time steps $\Delta t=5/64^2\times 2^{-2k}$  ($k=0,1,2,3$). As shown in Table \ref{Table-Ex1-1} where the maximum and $l^2$ norms of the numerical errors and CRs are presented, the RLB-MIE-FD scheme has a second-order CR in space.
\begin{table}[H]
	\begin{center} 
		\caption{The numerical errors in the maximum and $l^2$ norms and CRs of the RLB-MIE-FD scheme.} 	\label{Table-Ex1-1}
		\begin{tabular}{ccccccccccc }\hline\hline
			$\Delta x$& Err$(\bm{\Phi}^{\star},\infty)$&CR$(\bm{\Phi}^{\star},\infty)$& Err$(\bm{\Phi}^{\star},2)$&CR$(\bm{\Phi}^{\star},2)$  \\
			\hline
			$1/64$&9.8980 $\times 10^{-3}$&- &2.0477$\times 10^{-3}$&-        \\
			$1/128$&2.4823$\times 10^{-3}$&1.9954&5.1281$\times 10^{-4}$&1.9975 \\
			$1/256$&6.1411$\times 10^{-4}$&2.0151&1.2601$\times 10^{-4}$&2.0147    \\
			$1/512$&1.5323$\times 10^{-4}$&2.0027&3.1678$\times 10^{-5}$&2.0022 \\ 
			\hline\hline
		\end{tabular}
	\end{center} 
	\vspace{-1.2cm}   
\end{table} 

\subsubsection{The one-dimensional ACE with HD-BC}  
\begin{figure}[H]
	\vspace{-0.2cm}  
	\setlength{\belowcaptionskip}{-1cm}   
	\begin{center}           
		\subfloat[]    
		{
			\includegraphics[width=0.33\textwidth]{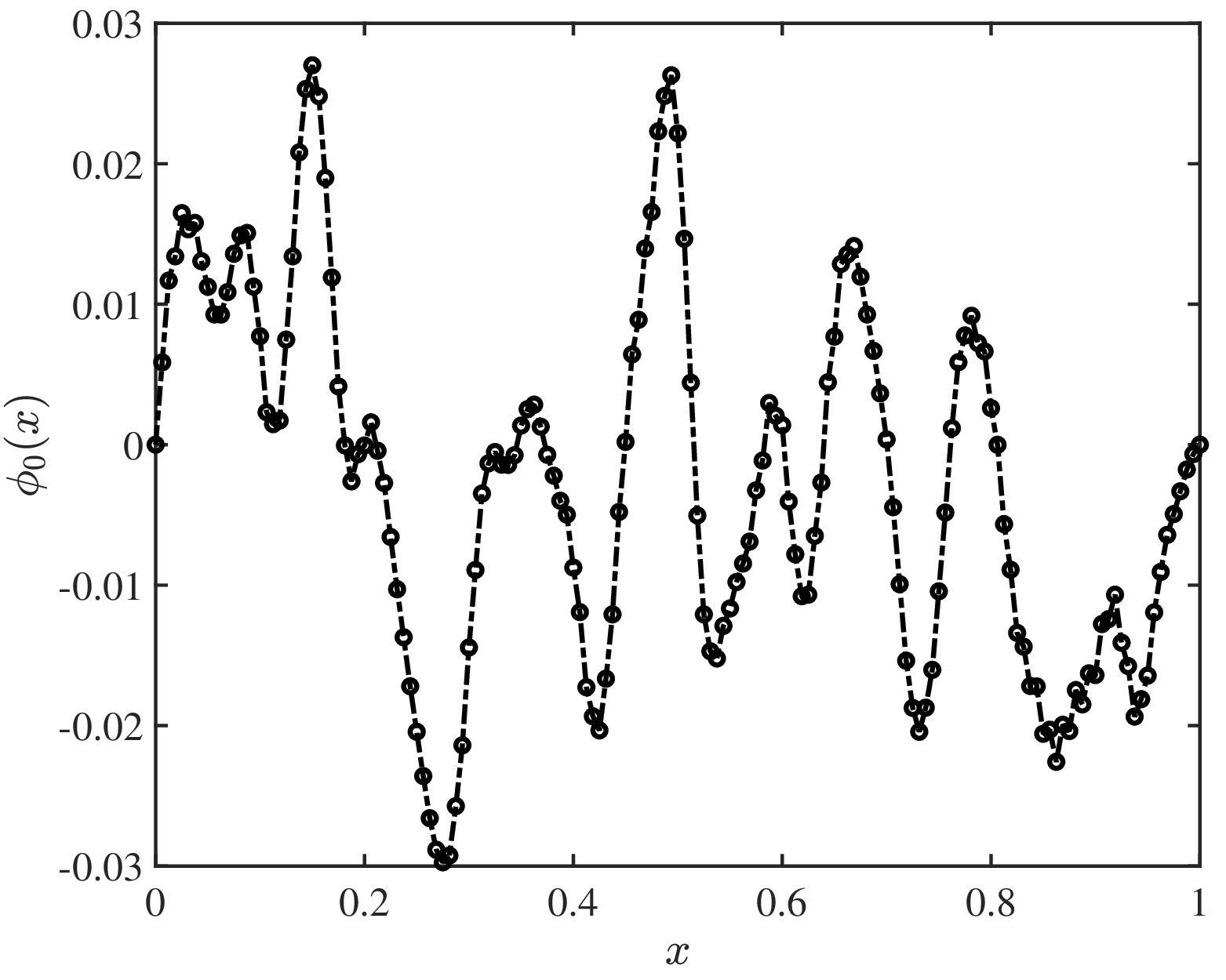}
		}
		\subfloat[   ]
		{
			\includegraphics[width=0.33\textwidth]{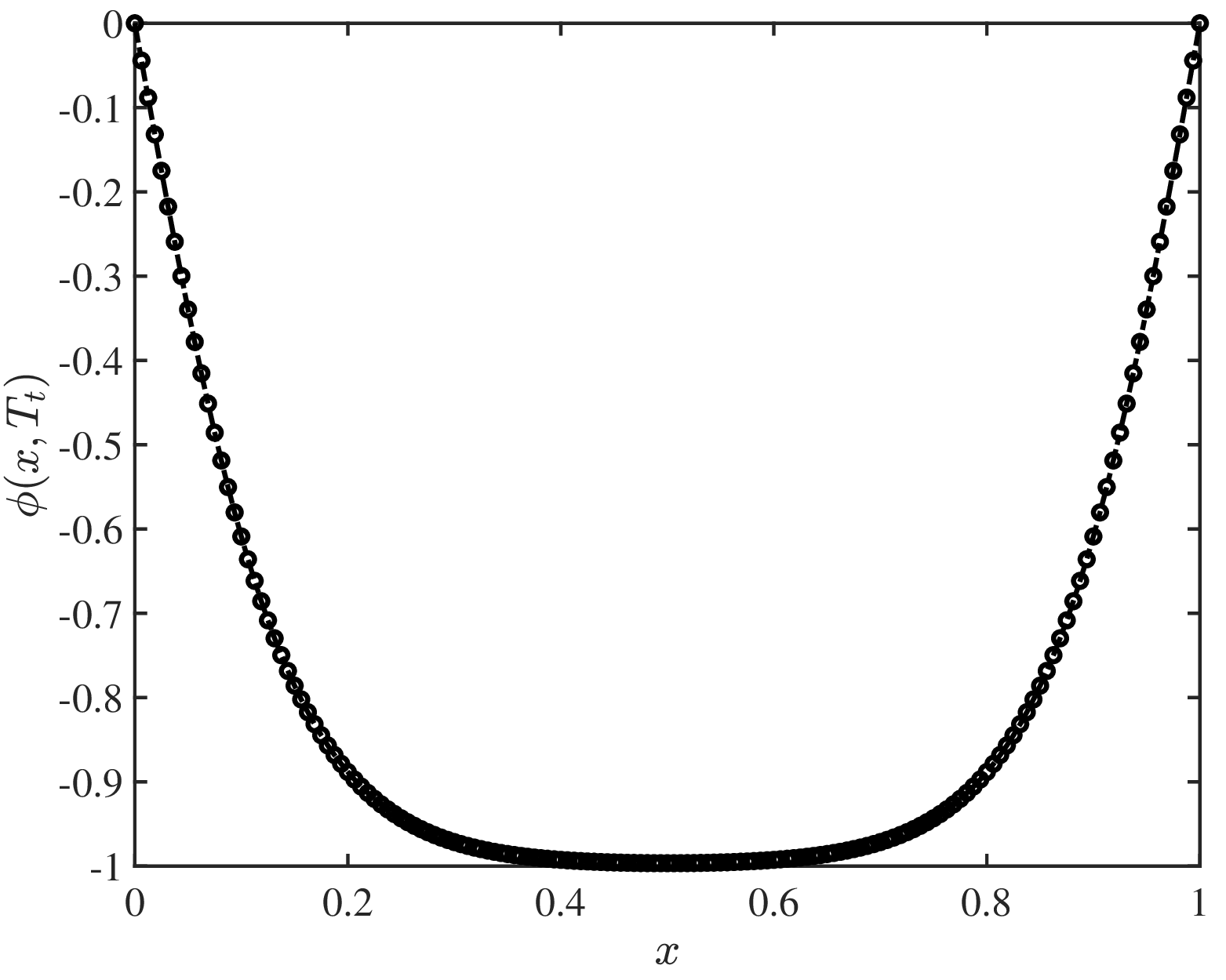}
		}  
		\caption{The initial solution at time $T_t=0$ (a) and numerical solution at time $T_t=40$ (b).} 
		\label{Fig-Ex2} 
	\end{center}       
\end{figure}

\begin{figure}[H] 
	\vspace{-0.2cm}  
	\setlength{\belowcaptionskip}{-1cm}   
	\begin{center}           
		\subfloat[]    
		{
			\includegraphics[width=0.33\textwidth]{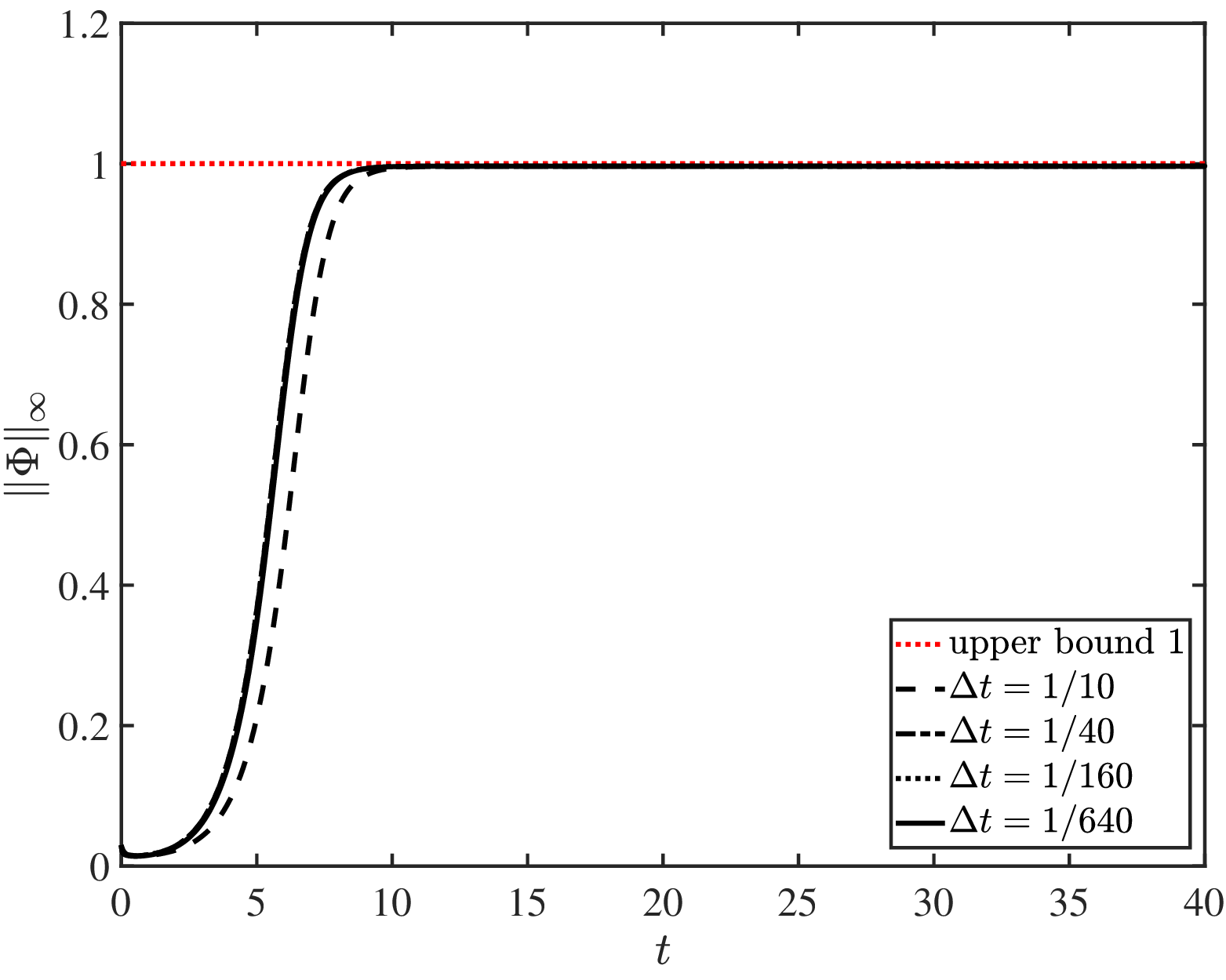}
		}
		\subfloat[   ]
		{
			\includegraphics[width=0.33\textwidth]{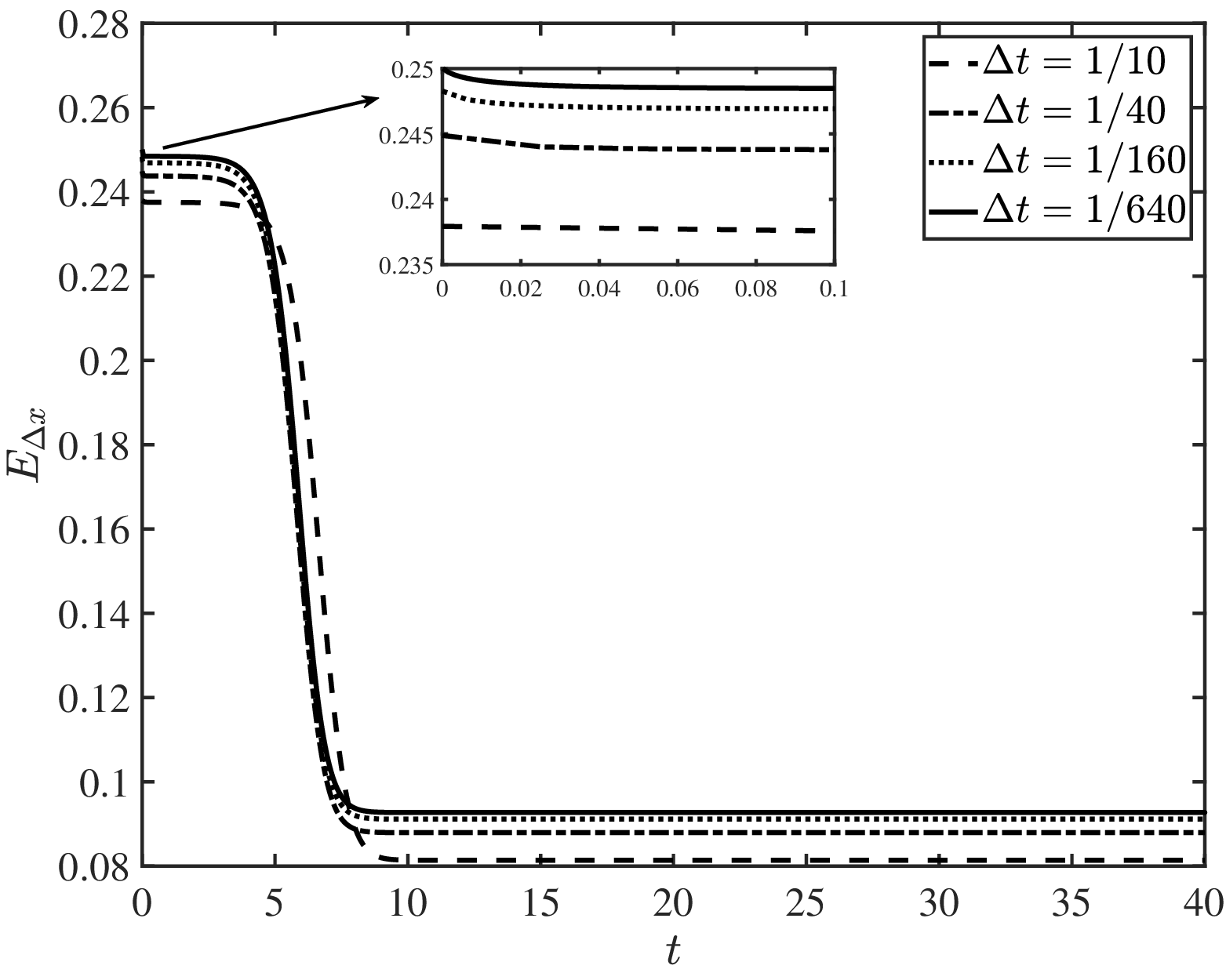}
		}  
		\caption{The discrete maximum norms (a) and energies (b) at different values of time step $\Delta t$.} 
		\label{Fig-Ex2-2}  
	\end{center}      
\end{figure}

\begin{figure}[H] 
	\vspace{-0.2cm}  
	\setlength{\belowcaptionskip}{-1cm}   
	\begin{center}           
		\subfloat[]    
		{
			\includegraphics[width=0.33\textwidth]{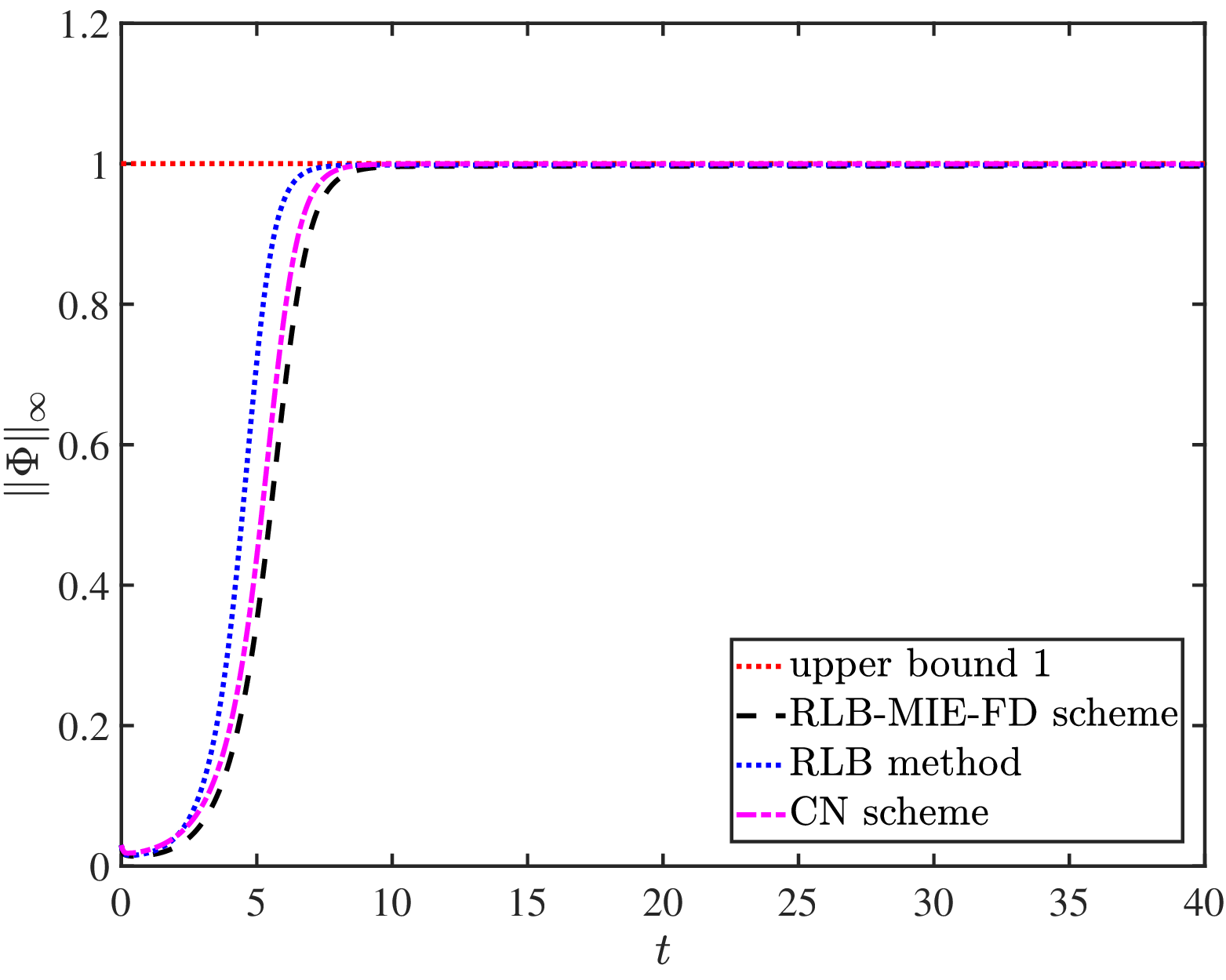}
		}
		\subfloat[   ]
		{
			\includegraphics[width=0.33\textwidth]{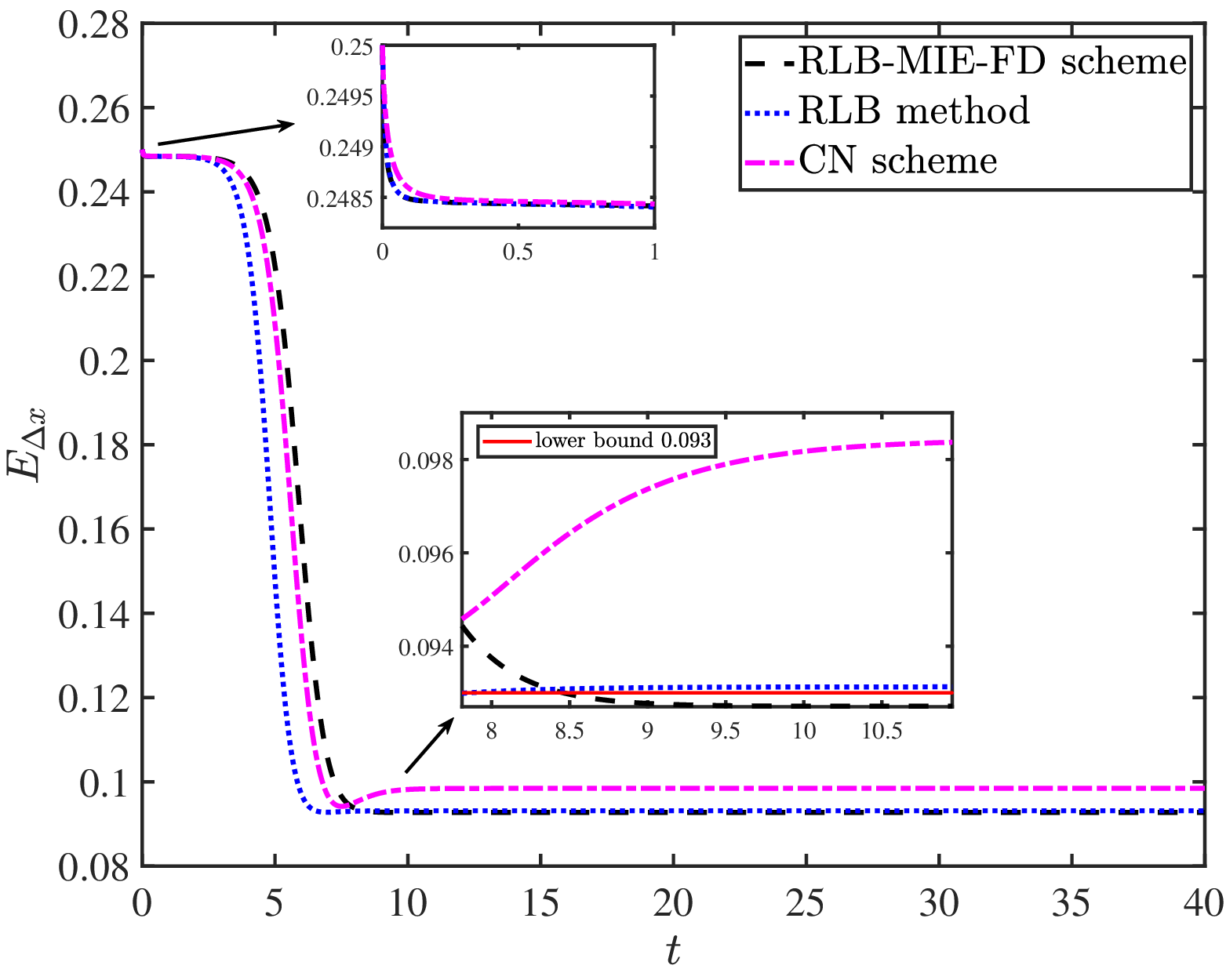}
		}  
		\caption{The discrete maximum norms (a) and energies (b) of the numerical solutions of the RLB-MIE-FD scheme,  RLB method and  CN scheme.} 
		\label{Fig-Ex2-3}  
	\end{center}      
\end{figure}

The second example is the ACE (\ref{ACE}) with the initial value $\phi_0(x)=0.1\times \makebox{rand}(x)-0.05$, $x\in(0,1)$ and the HD-BC [see Fig. \ref{Fig-Ex2}(a)]. We first present the numerical result at time $T_t=40$ in Fig. \ref{Fig-Ex2}(b) where $\epsilon=0.4$, the lattice spacing $\Delta x=1/80$ and time step $\Delta t=1/640$.  From this figure, one can clearly see that the numerical solution is indeed located in the range $[-1, 1]$.
We then further test the maximum values of the numerical solutions and discrete energies in Figs. \ref{Fig-Ex2-2}(a) and \ref{Fig-Ex2-2}(b)  where the time step $\Delta t$ satisfies Eq. (\ref{maximum-condition}) under the condition of  $\Delta x^2/\Delta t=1/10$. As seen from these figures, the discrete maximum principle and energy dissipation law can be preserved, which is in agreement with our theoretical results for the case of the HD-BC. 
Finally, we illustrate the evolutions of the discrete maximum norms and energies of the numerical solutions of the RLB-MIE-FD scheme, RLB method \cite{wang2015regularized}, and CN scheme in Fig. \ref{Fig-Ex2-3} where the time step $\Delta t=1/640$ and lattice spacing $\Delta x=1/160$. As shown in this figure,  all three methods perfectly preserve the maximum principle. However, compared to the RLB-MIE-FD scheme, the discrete energies of the RLB method and CN scheme are not strictly non-increasing [see Fig. \ref{Fig-Ex2-3}(b)], which indicates that to preserve the maximum principle and energy dissipation law, 
 both the RLB method and CN scheme require much stricter conditions than the present RLB-MIE-FD scheme. Additionally, it is worth noting that  compared to the RLB method, the RLB-MIE-FD scheme significantly reduces the memory occupancy (see Ref. \cite{Liu2023rlb}), and it is also more efficient than the CN scheme since it does not need to solve a nonlinear system at each step (see Remark 1).

\subsection{The two- and three-dimensional ACEs}
\subsubsection{The two-dimensional ACE with periodic boundary condition}
In this part, we   consider the two-dimensional ACE (\ref{ACE}) with  the P-BC and smooth initial data $u_0(x,y)=0.05\sin x\sin y$, $(x,y)\in(0,2\pi)\times (0,2\pi)$ \cite{du2019maximum}. 

\begin{figure}[H]
	\vspace{-0.6cm}  
	\setlength{\belowcaptionskip}{-1.2cm}   
	\begin{center}  
		\subfloat[$T_t=0$]{
			\begin{minipage}[t]{0.22\linewidth}
				\centering
				\includegraphics[width=1.6in]{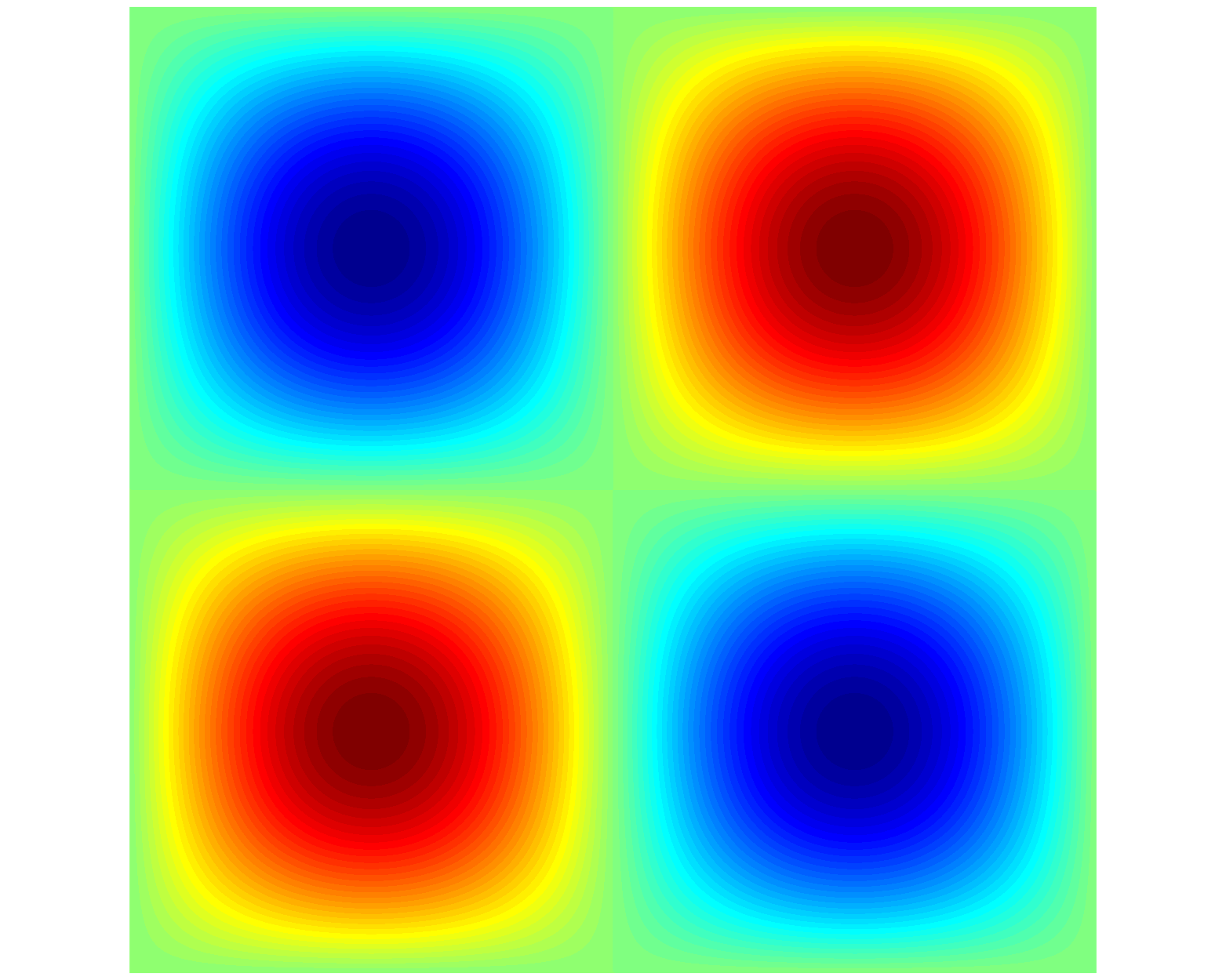} 
			\end{minipage} 
		}
		\subfloat[$T_t=4$]{
			\begin{minipage}[t]{0.22\linewidth}
				\includegraphics[width=1.6in]{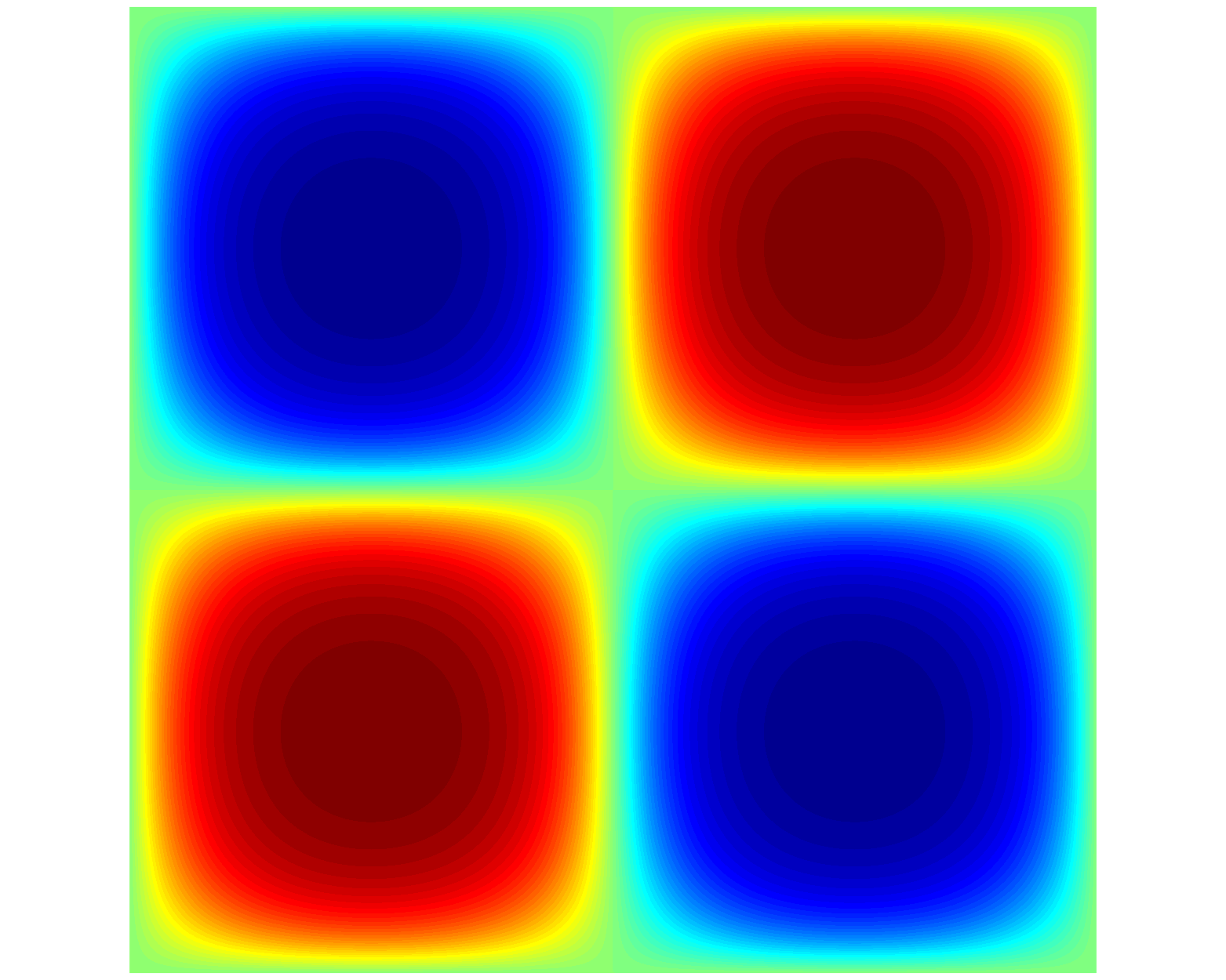} 
			\end{minipage} 
		}
		\subfloat[$T_t=6$]{
			\begin{minipage}[t]{0.22\linewidth}
				\centering
				\includegraphics[width=1.6in]{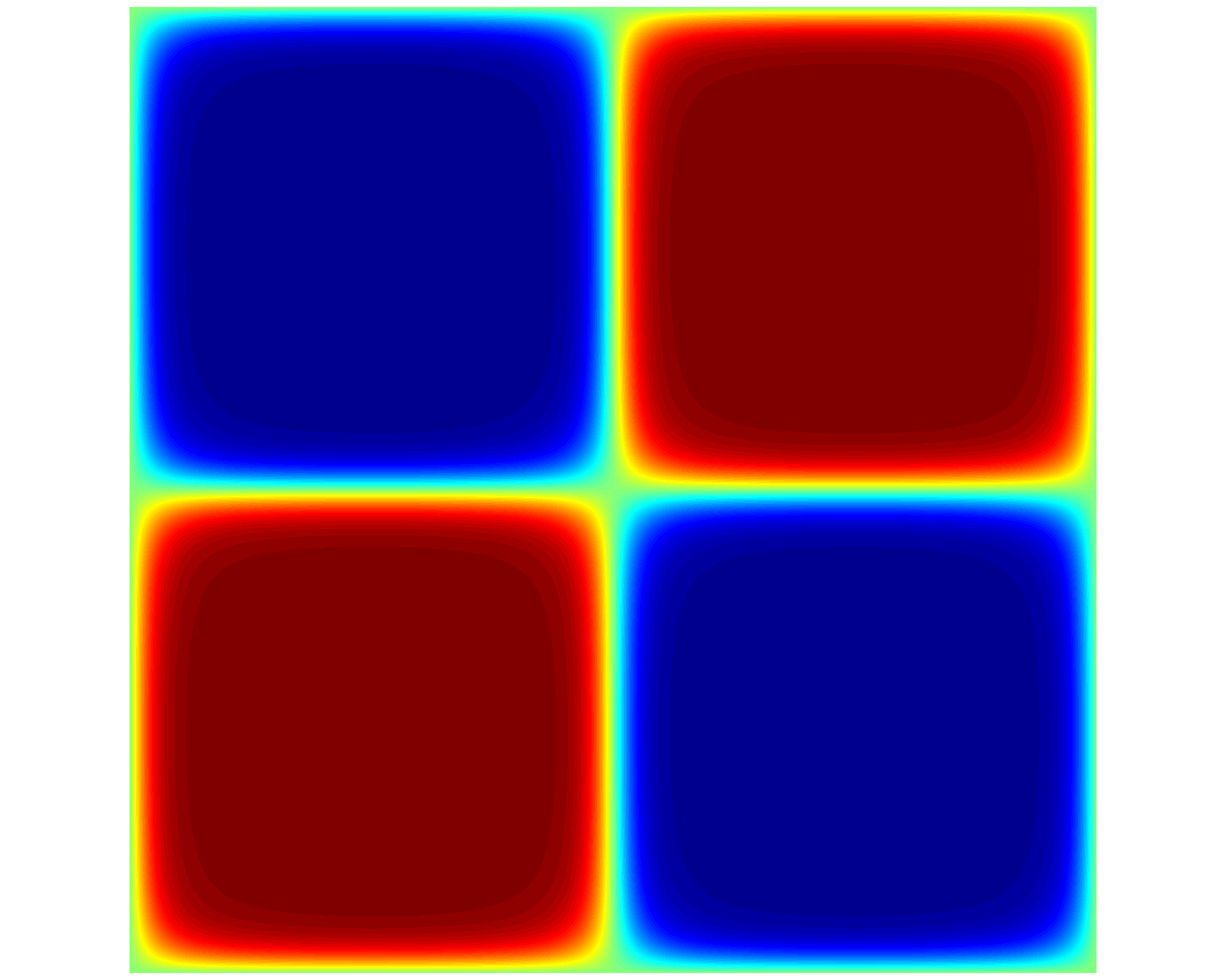} 
			\end{minipage} 
		}
		\subfloat[$T_t=15$]{
			\begin{minipage}[t]{0.22\linewidth}
				\centering
				\includegraphics[width=1.6in]{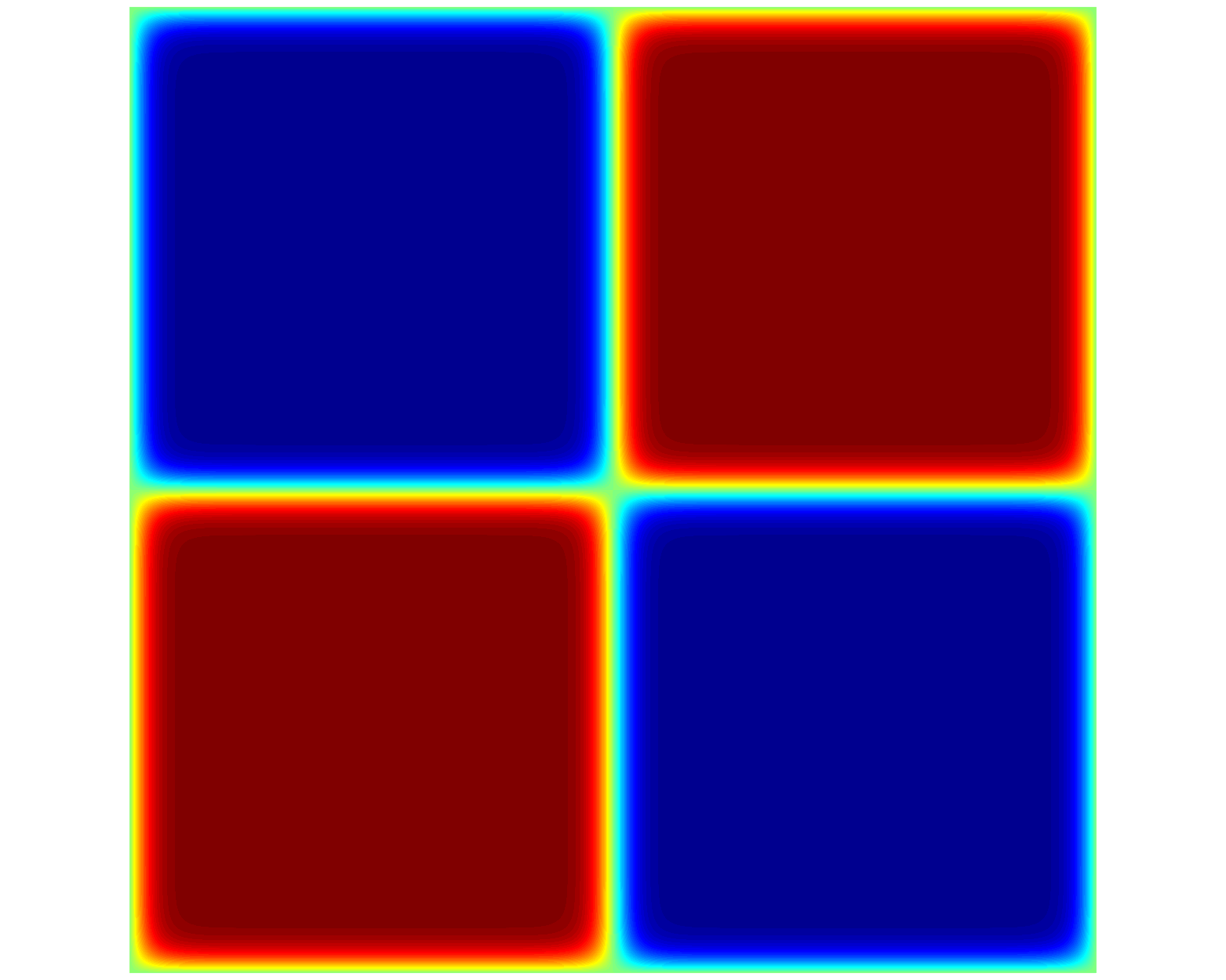} 
			\end{minipage} 
		}  
		\caption{Temporal evolution of the numerical solution.}   	\label{Fig-Ex3-1}    
	\end{center} 
\end{figure}

\begin{figure}[H]
	\vspace{-0.45cm}  
	\setlength{\belowcaptionskip}{-1.2cm}   
	\begin{center}          
		\subfloat[]    
		{
			\includegraphics[width=0.33\textwidth]{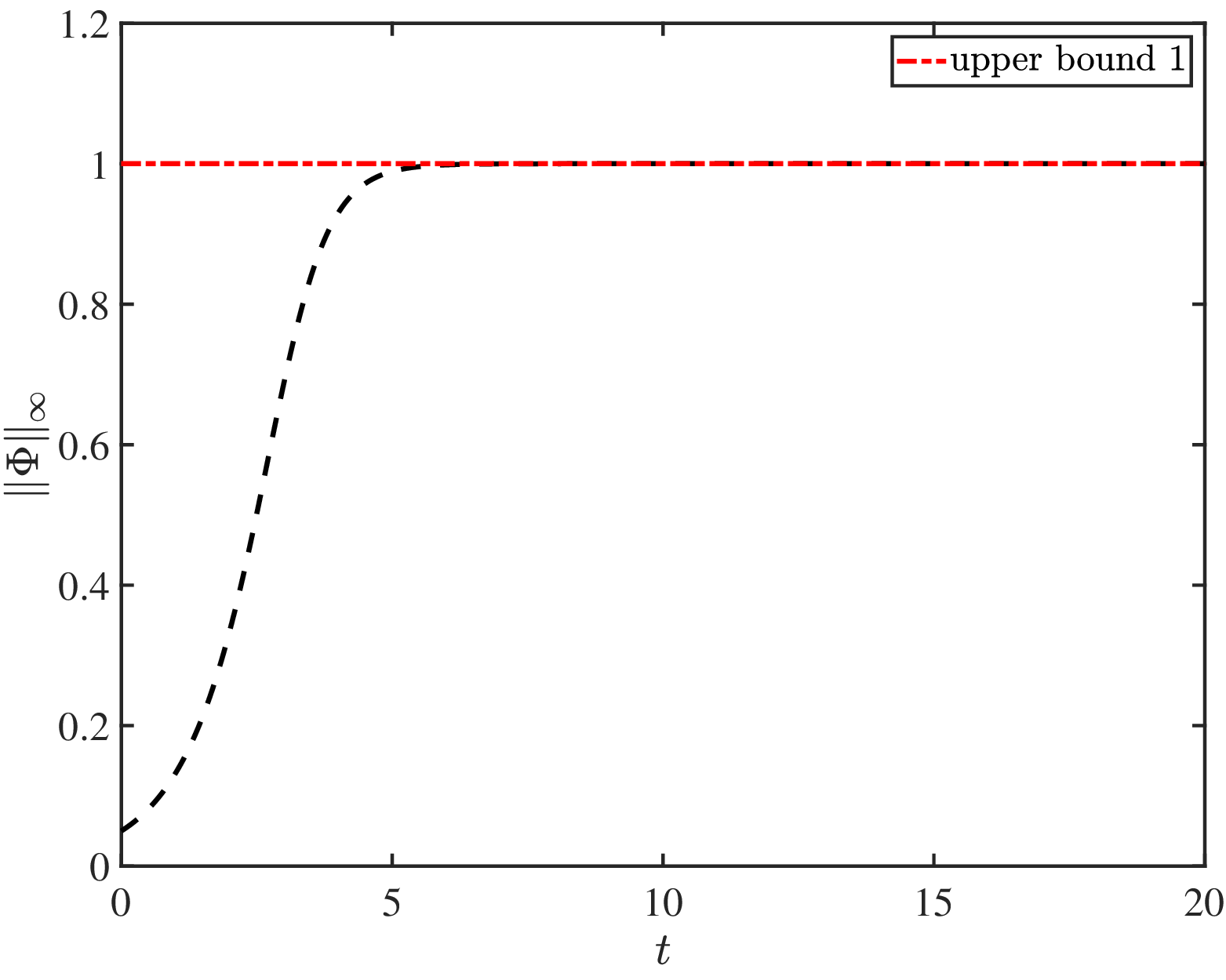}
		}
		\subfloat[   ]
		{
			\includegraphics[width=0.33\textwidth]{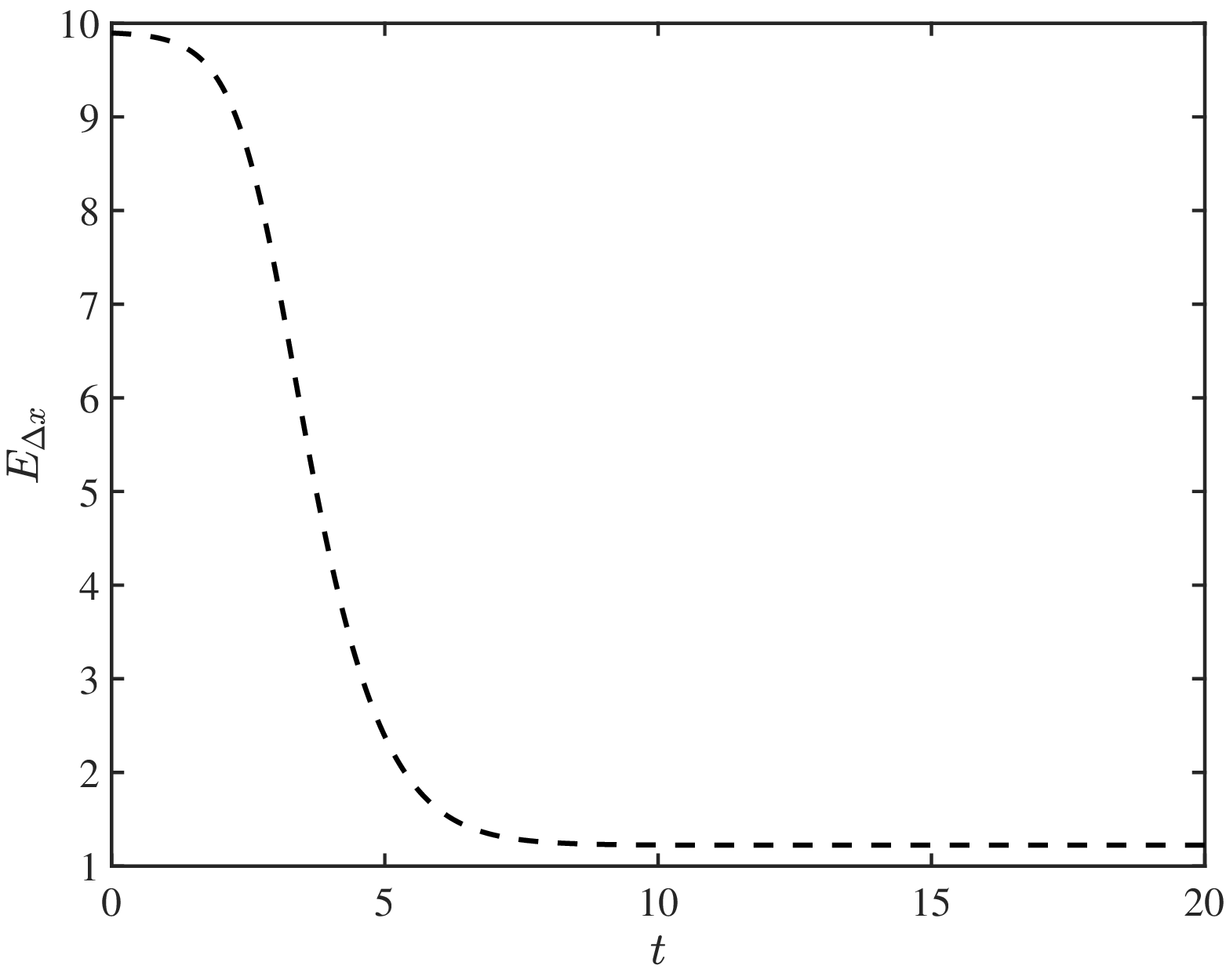}
		} 
		\caption{The discrete maximum norm (a) and energy (b) of the numerical solution.}		 	\label{Fig-Ex3-2}  
	\end{center} 
\end{figure}

We first take the interfacial parameter $\varepsilon= 0.01$ and plot the evolutions of scalar variable $\phi$, discrete maximum norm and energy of numerical solution  in Figs. \ref{Fig-Ex3-1} and  \ref{Fig-Ex3-2}  where the lattice spacing and time step are set to be $\Delta x=\pi/256$ and $\Delta t=1/800$, respectively. 	From these two figures, one can see that   the numerical solution can preserve the discrete maximum principle and energy dissipation law well. 

  Finally, we carry out some simulations to test the CR  of the RLB-MIE-FD scheme in space  under different values of the interfacial parameter $\varepsilon$ ($=0.01,0.2$). The numerical solution $\overline{\Phi}$ of the ACE (\ref{ACE}) with the finest lattice spacing $\overline{\Delta x}$ and time step $\overline{\Delta t}$ is treated as a reference data to calculate the errors of the numerical solutions at different values of lattice spacing $\Delta x=\pi/(10\times k)$ ($k=1, 2,4,5$), and the time step $\Delta t$ is determined by $\Delta x^2/\Delta t=\pi^2/2$. We present the numerical errors with respect to the maximum and $l^2$ norms in Table \ref{Table-Ex3-1} where the time $T_t=1$. From this table, one can find that the RLB-MIE-FD scheme   has a second-order CR in space, which is the same as the theoretical analysis.
\begin{table}  [H]		 	
	\vspace{-0.6cm}  
	\begin{center}		 
		\caption{The CRs of the RLB-MIE-FD scheme under two different values of interfacial parameter $\varepsilon$.}  \label{Table-Ex3-1}
		\begin{tabular}{ccccccccccc }\hline\hline
			$\varepsilon$&$(\overline{\Delta x},\overline{\Delta t})$ &$\Delta x$ &Err$(\bm{\overline{\Phi}},\infty)$&CR$(\bm{\overline{\Phi}},\infty)$&Err$(\bm{\overline{\Phi}},2)$&CR$(\bm{\overline{\Phi}},2)$\\\hline
			\multirow{4}{*}{0.01}&\multirow{4}{*}{($\pi/400,1/4000$)}& $\pi/10$ &4.8590 $\times 10^{-5}$&-  &2.4170$\times 10^{-5}$&-      \\
			
			&	&$\pi/20$ &1.1245$\times 10^{-5}$&2.1114&5.7335$\times 10^{-6}$&2.0757\\
			
			&	    &$\pi/40$ &2.5260$\times 10^{-6}$&2.1543&1.3056$\times 10^{-6}$&2.1347   \\
			
			&		&$\pi/50$ &1.4964$\times 10^{-6}$&2.3464&7.7615$\times 10^{-7}$&2.3306\\ 
			\cline{1-7}
			\multirow{4}{*}{0.2}&\multirow{4}{*}{($\pi/200,1/10000$)}& $\pi/10$ &5.0667$\times 10^{-4}$&- &2.4771$\times 10^{-4}$&-       \\
			
			&		&$\pi/20$ &1.2509$\times 10^{-4}$&2.0181&6.1876$\times 10^{-5}$&1.9836\\
			
			&		&$\pi/40$ &2.9904$\times 10^{-5}$&2.0645 &1.4976$\times 19^{-5}$&2.0468 \\
			
			&		&$\pi/50$ &1.8488$\times 10^{-5}$&2.1550&9.2819$\times 10^{-6}$&2.1437\\ 
			\hline\hline
		\end{tabular}
	\end{center} 
	\vspace{-1.2cm}  
\end{table} 
\subsubsection{Motion by mean curvature}
\begin{figure}[H]
	\vspace{-0.2cm}  
	\begin{center}
		\subfloat[$T=0$]{
			\begin{minipage}[t]{0.22\linewidth}
				\centering
				\includegraphics[width=1.6in]{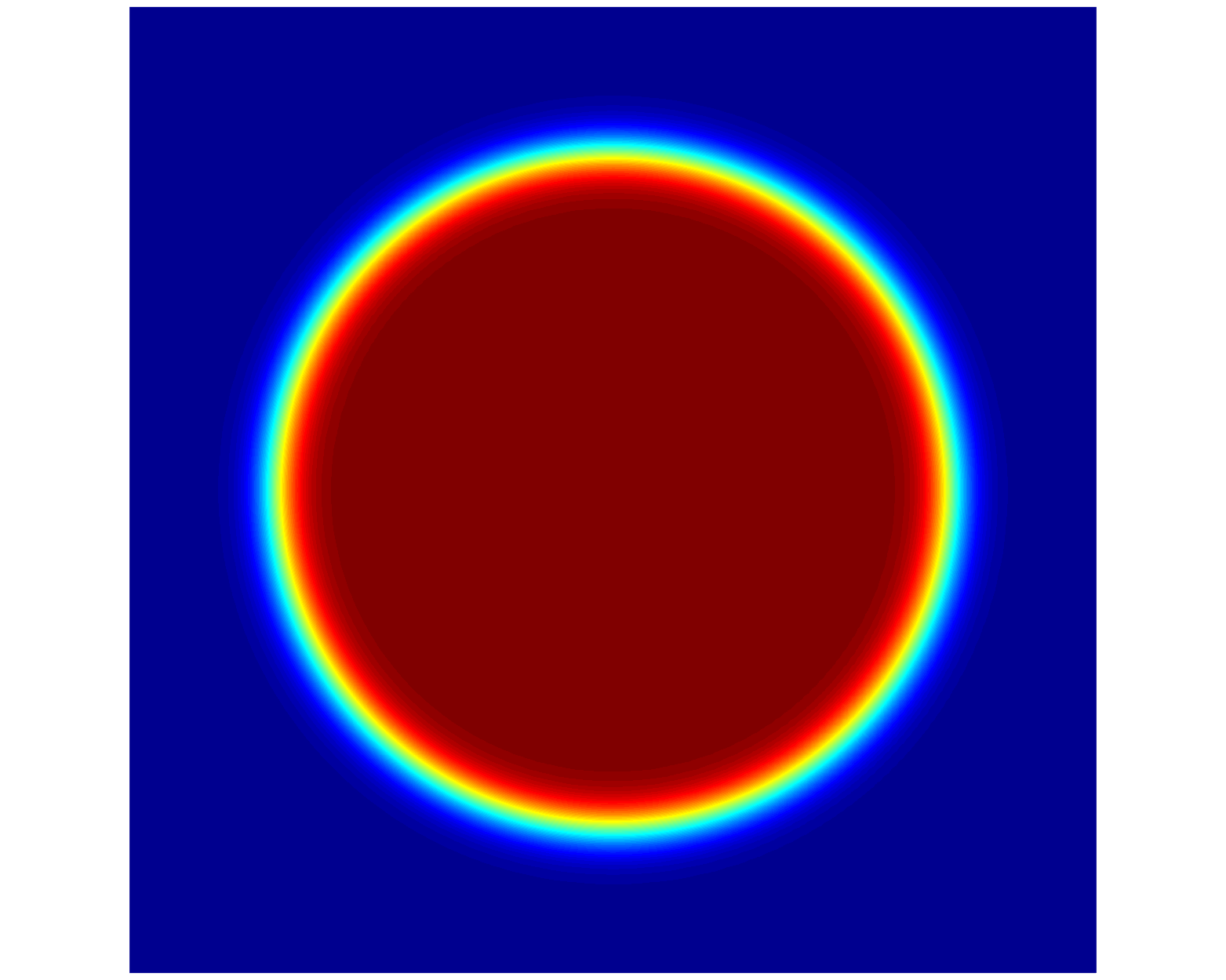} 
			\end{minipage} 
		}
		\subfloat[$T=1$]{
			\begin{minipage}[t]{0.22\linewidth}
				\includegraphics[width=1.6in]{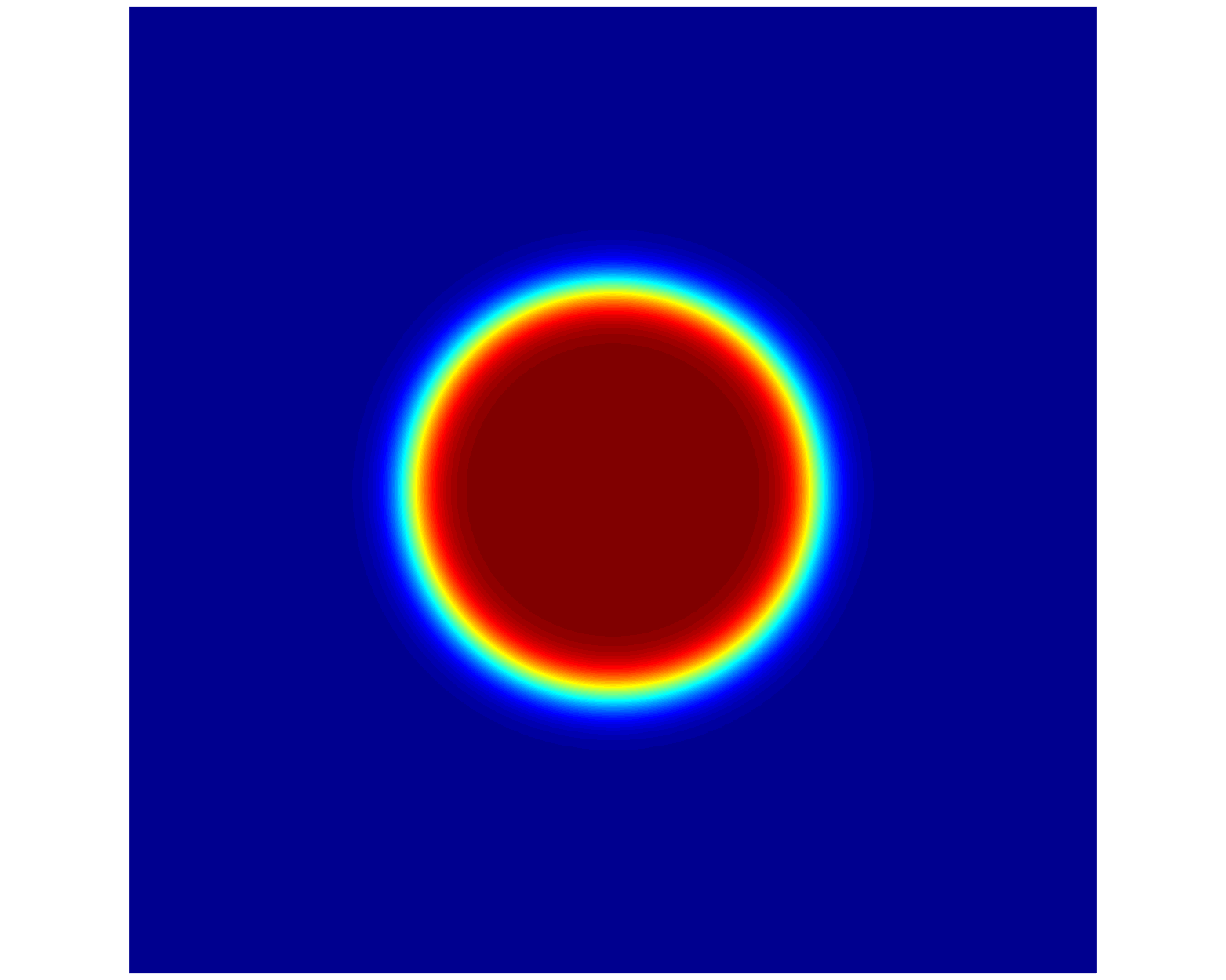}
			\end{minipage} 
		}
		\subfloat[$T=6$]{
			\begin{minipage}[t]{0.22\linewidth}
				\centering
				\includegraphics[width=1.6in]{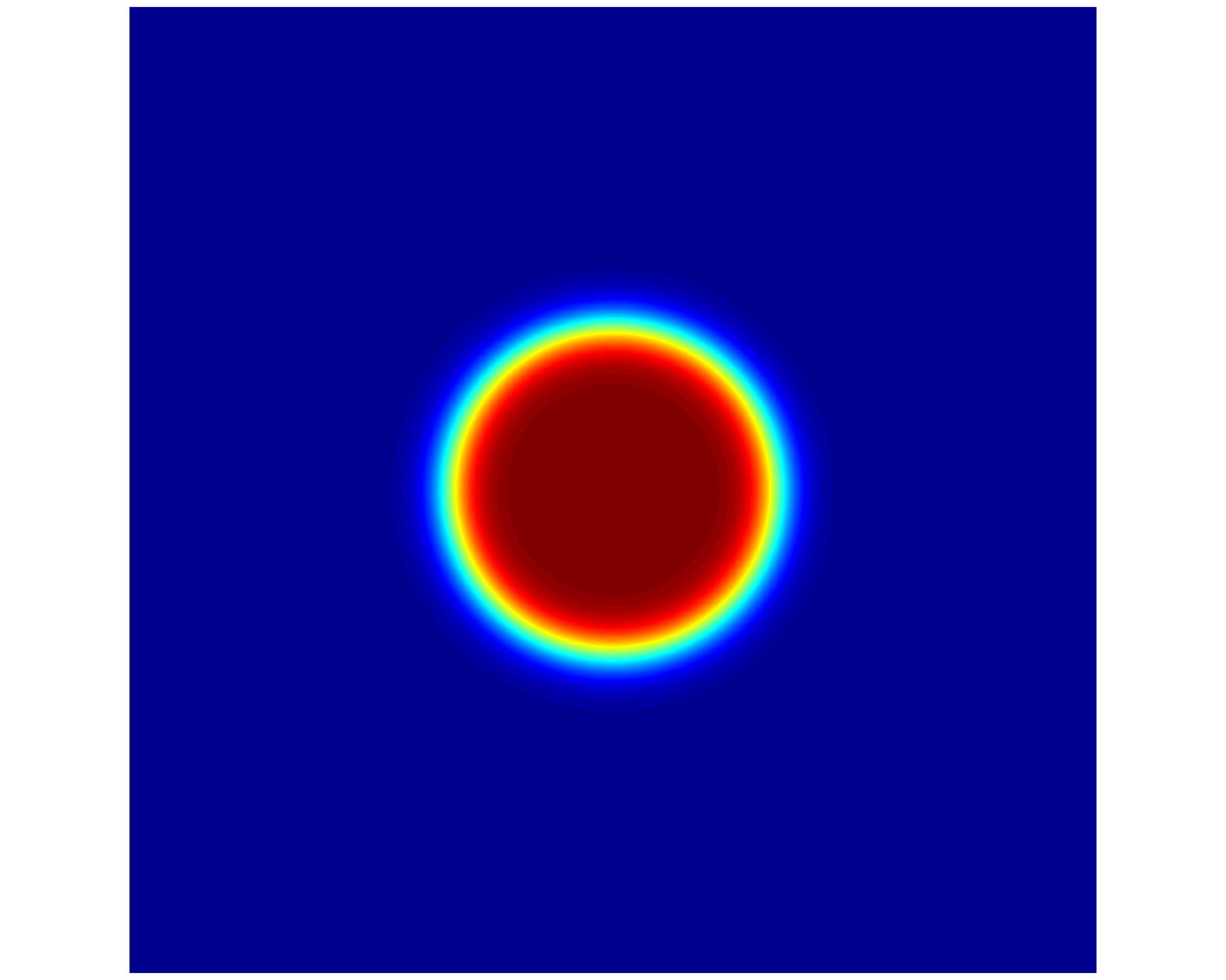}
			\end{minipage} 
		}
		\subfloat[$T=15$]{
			\begin{minipage}[t]{0.22\linewidth}
				\centering
				\includegraphics[width=1.6in]{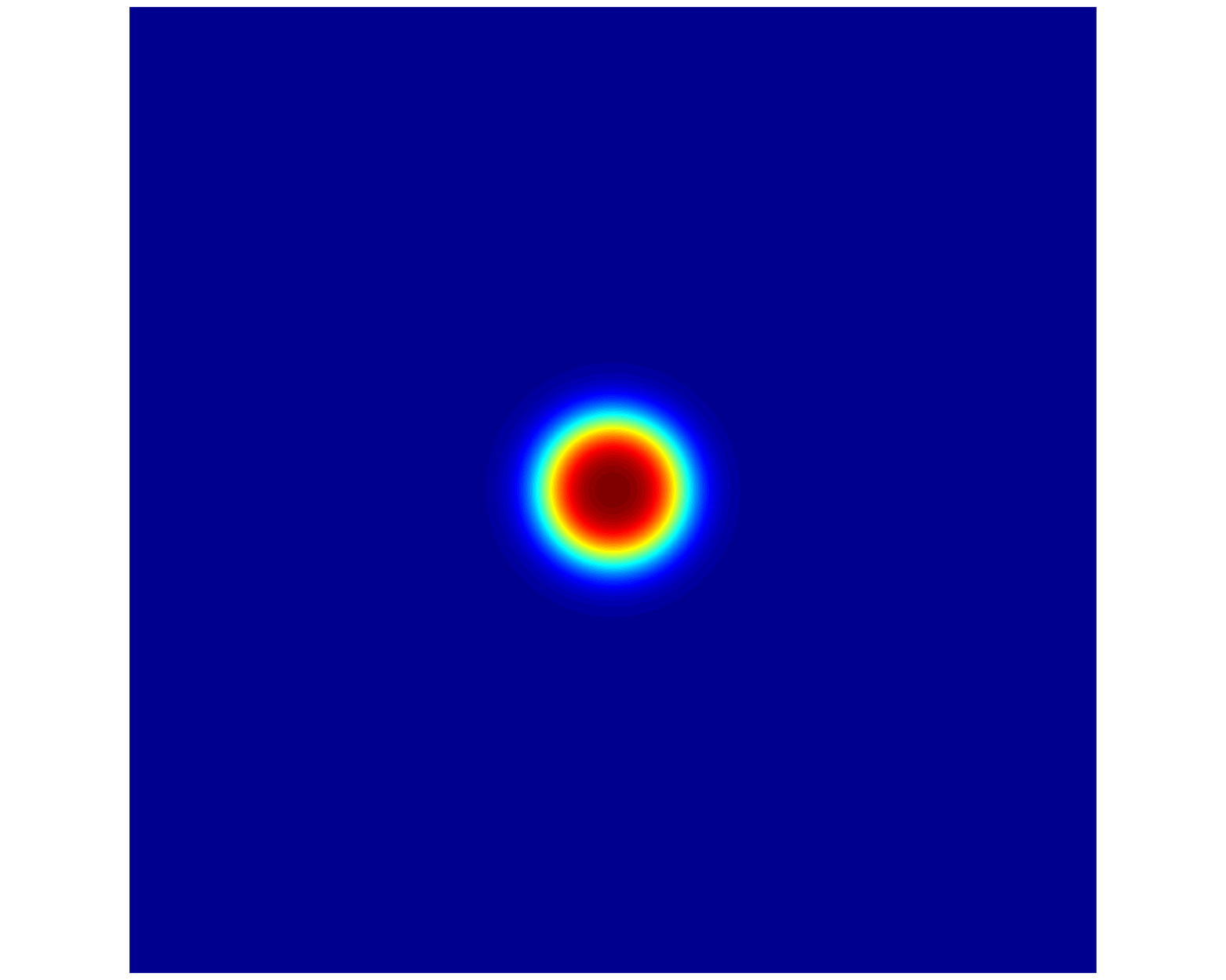}
			\end{minipage} 
		}
		
		\subfloat[   ]
		{
			\includegraphics[width=0.3\textwidth]{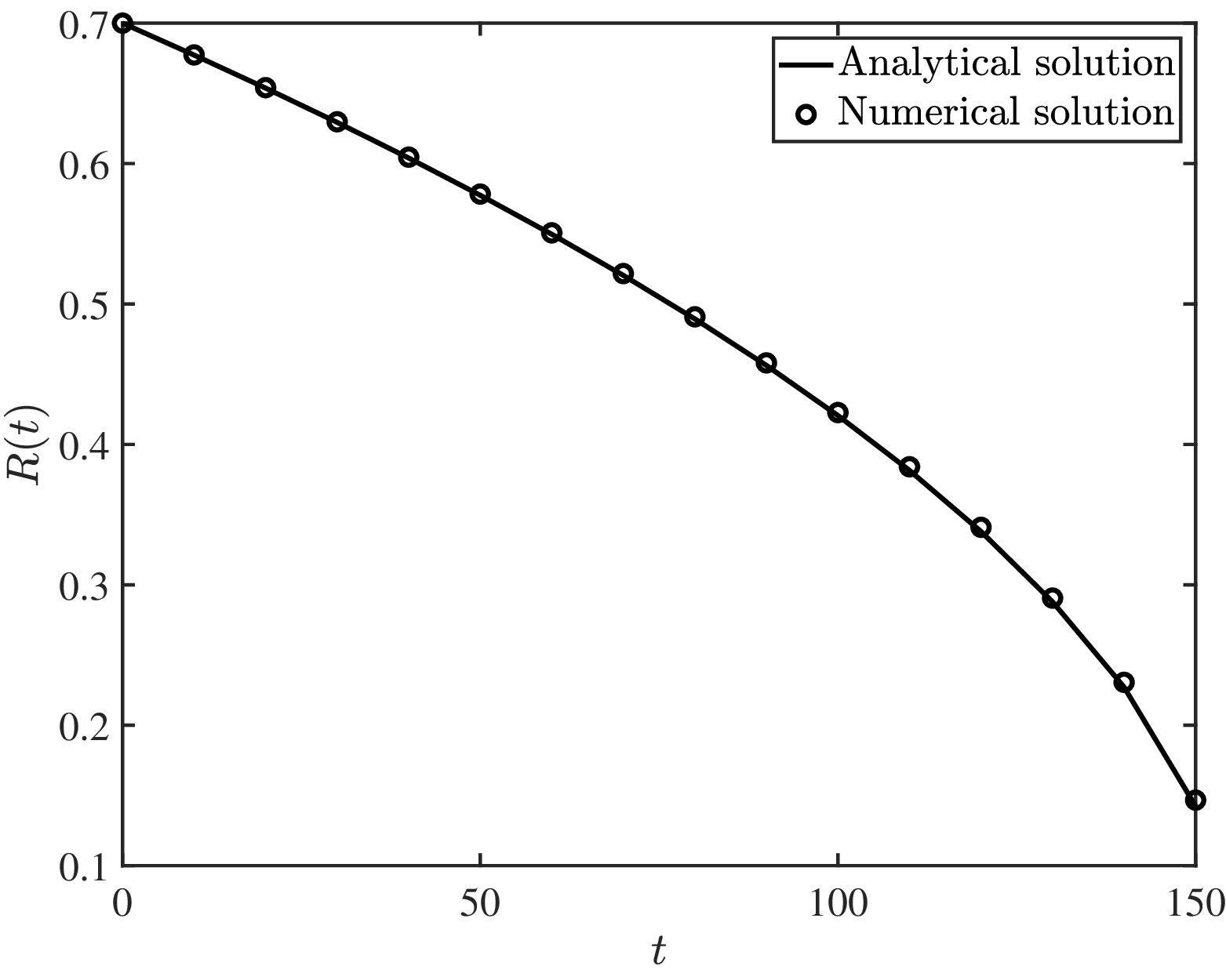}
		}	 
		\caption{The evolution of the numerical solution  of the two-dimensional circle in time.}  	\label{Fig-Ex4-1}  
	\end{center} 
\end{figure}

\begin{figure} [H]
	\vspace{-0.8cm}  
	\setlength{\belowcaptionskip}{-1cm}   
	\begin{center}     
		\subfloat[]    
		{
			\includegraphics[width=0.3\textwidth]{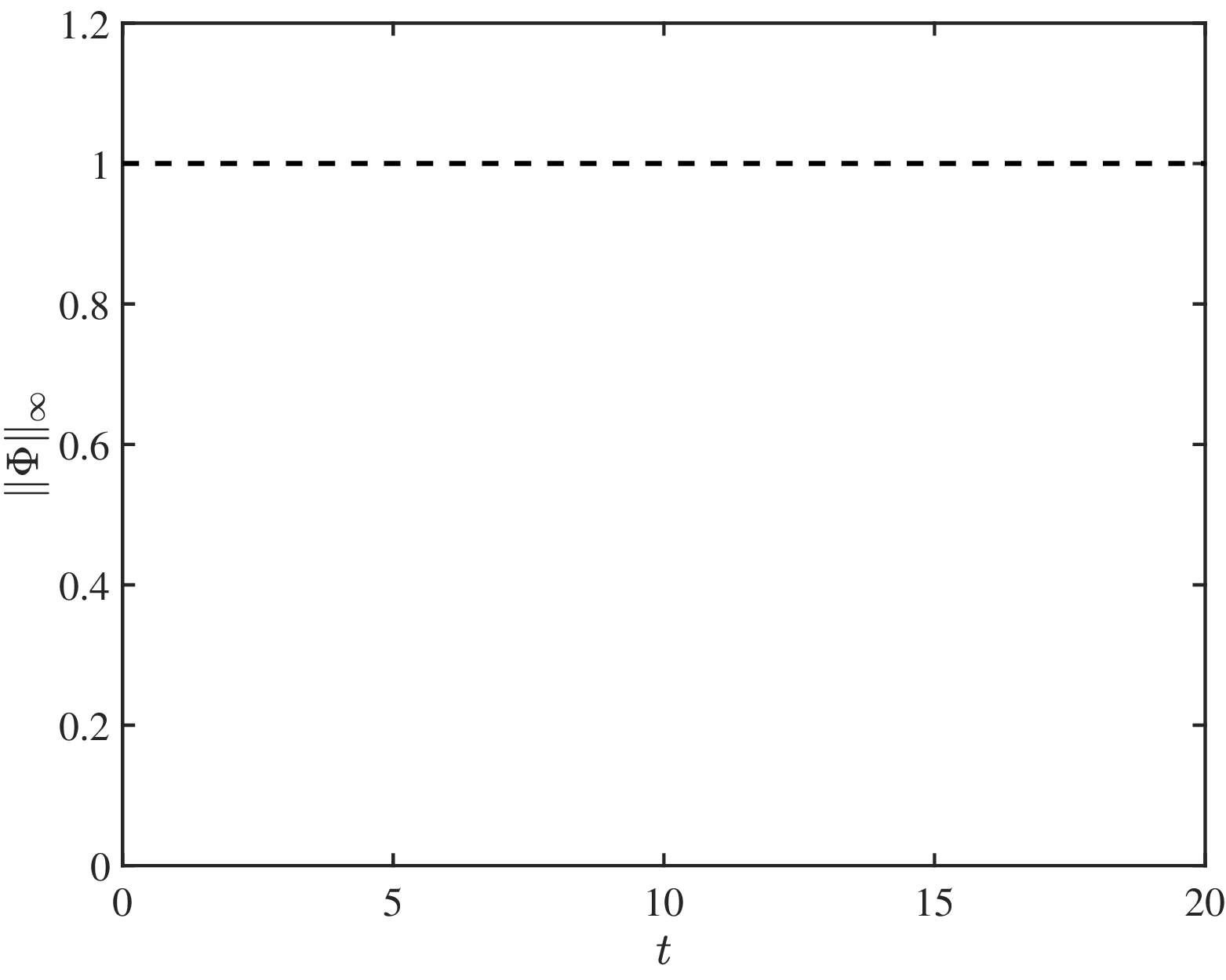}
		}
		\subfloat[   ]
		{
			\includegraphics[width=0.3\textwidth]{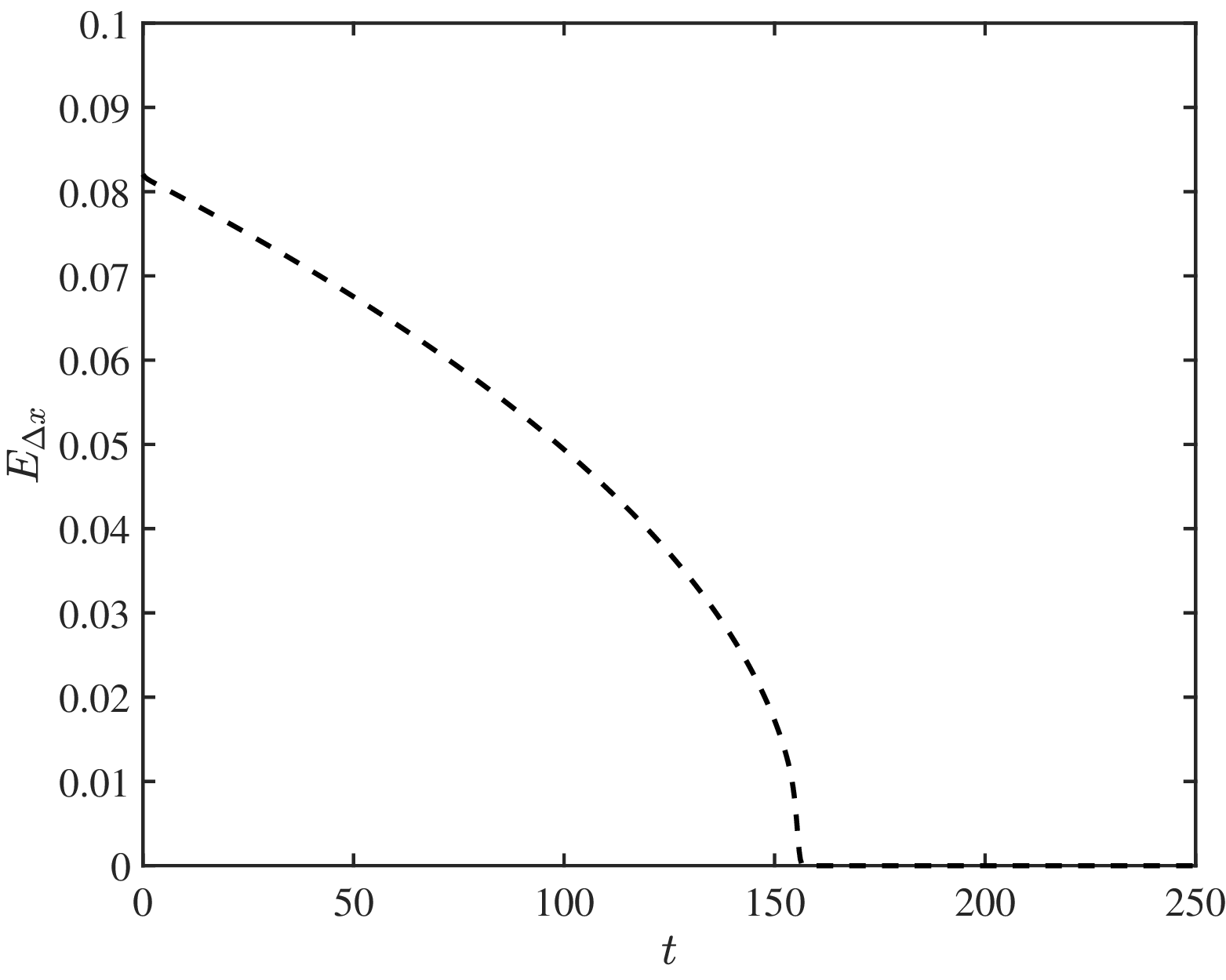}
		}	 
		\caption{The evolutions of the discrete maximum norm (a) and energy (b) of the numerical solution in time.}  \label{Fig-Ex4-2}       
	\end{center}      
	
\end{figure}

\begin{figure}[H]
	\vspace{-0.2cm}  
	\begin{center}
		\subfloat[$T=0$]{
			\begin{minipage}[t]{0.22\linewidth}
				\centering
				\includegraphics[width=1.6in]{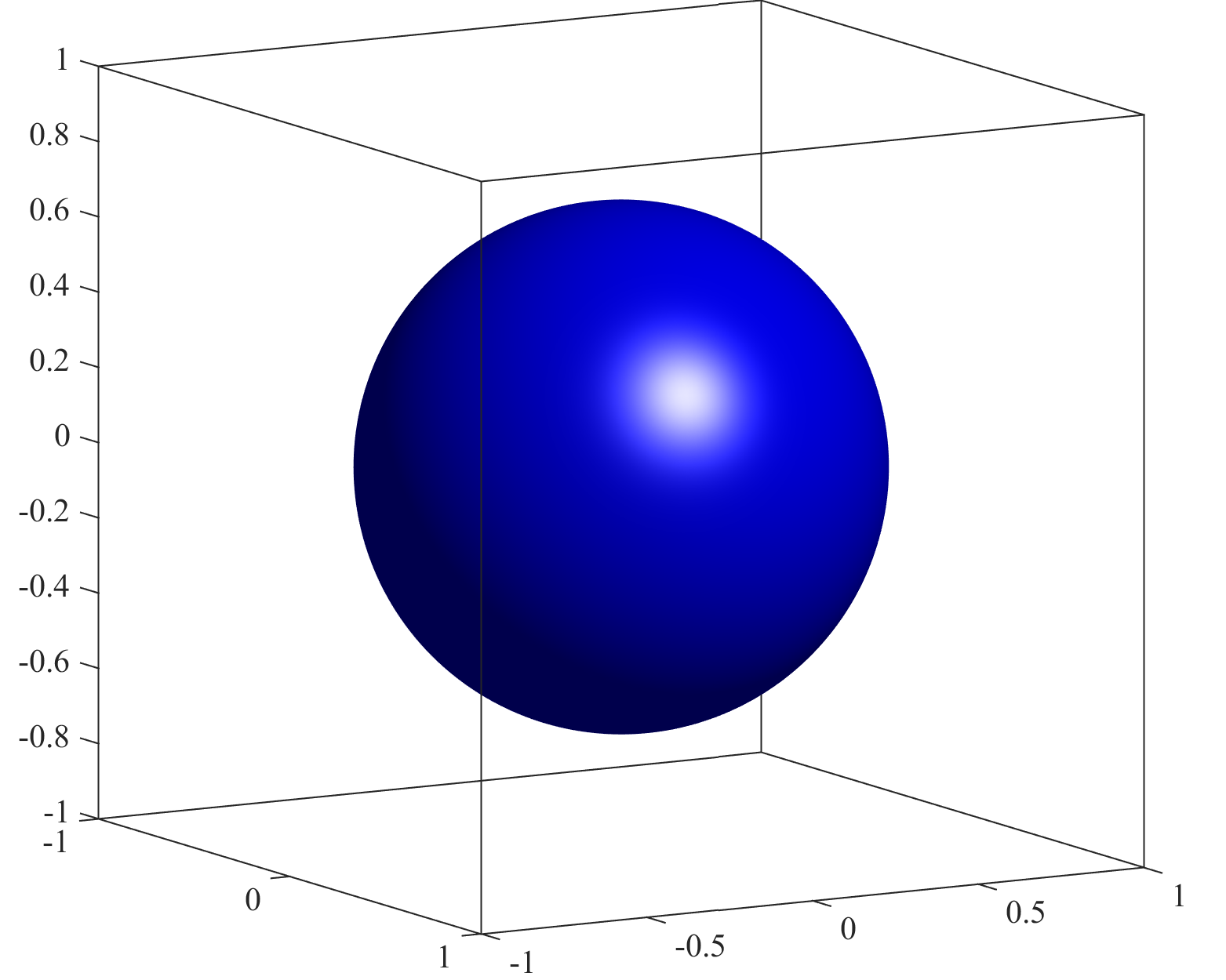} 
			\end{minipage} 
		}
		\subfloat[$T=28$]{
			\begin{minipage}[t]{0.22\linewidth}
				\includegraphics[width=1.6in]{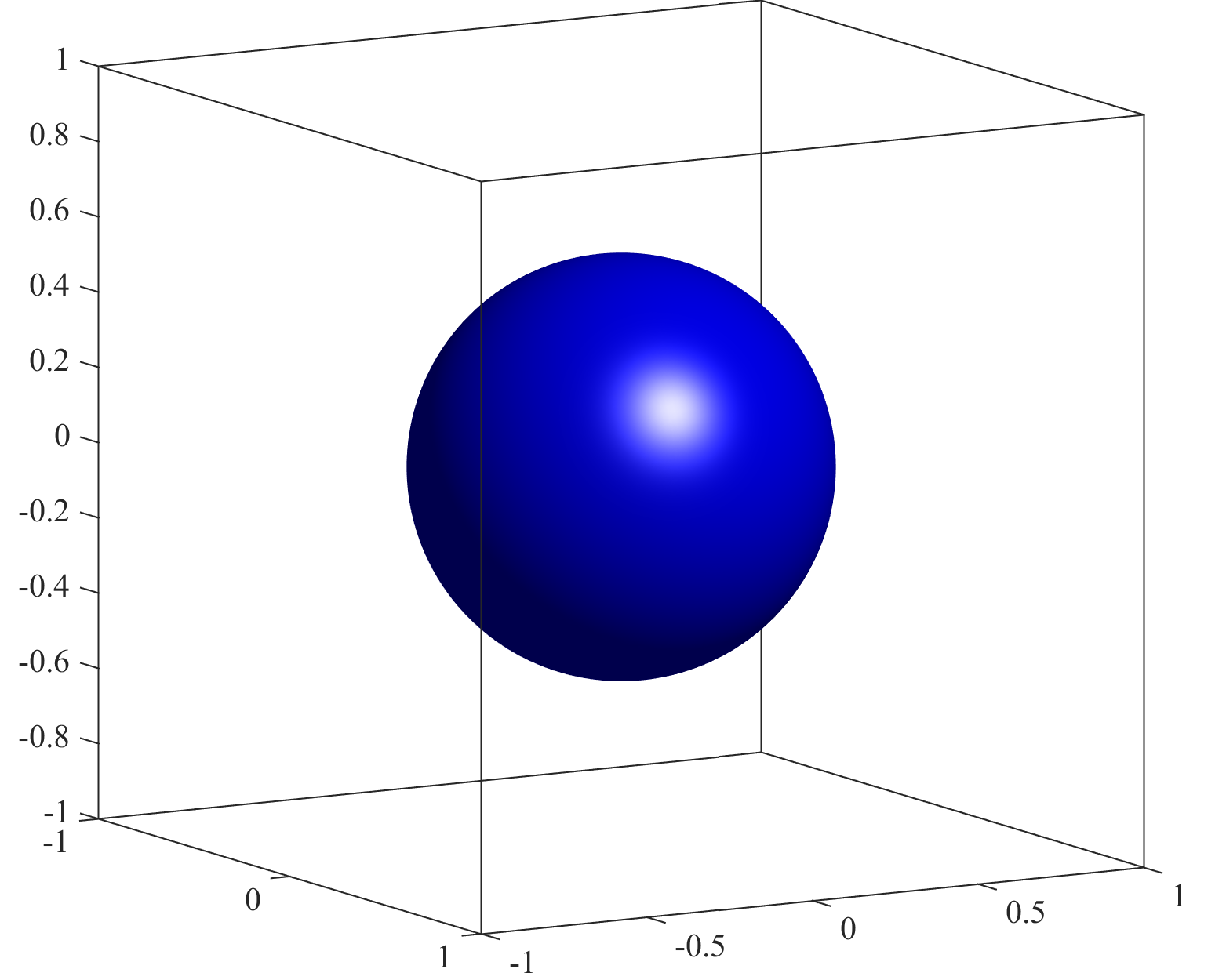} 
			\end{minipage} 
		}
		\subfloat[$T=56$]{
			\begin{minipage}[t]{0.22\linewidth}
				\centering
				\includegraphics[width=1.6in]{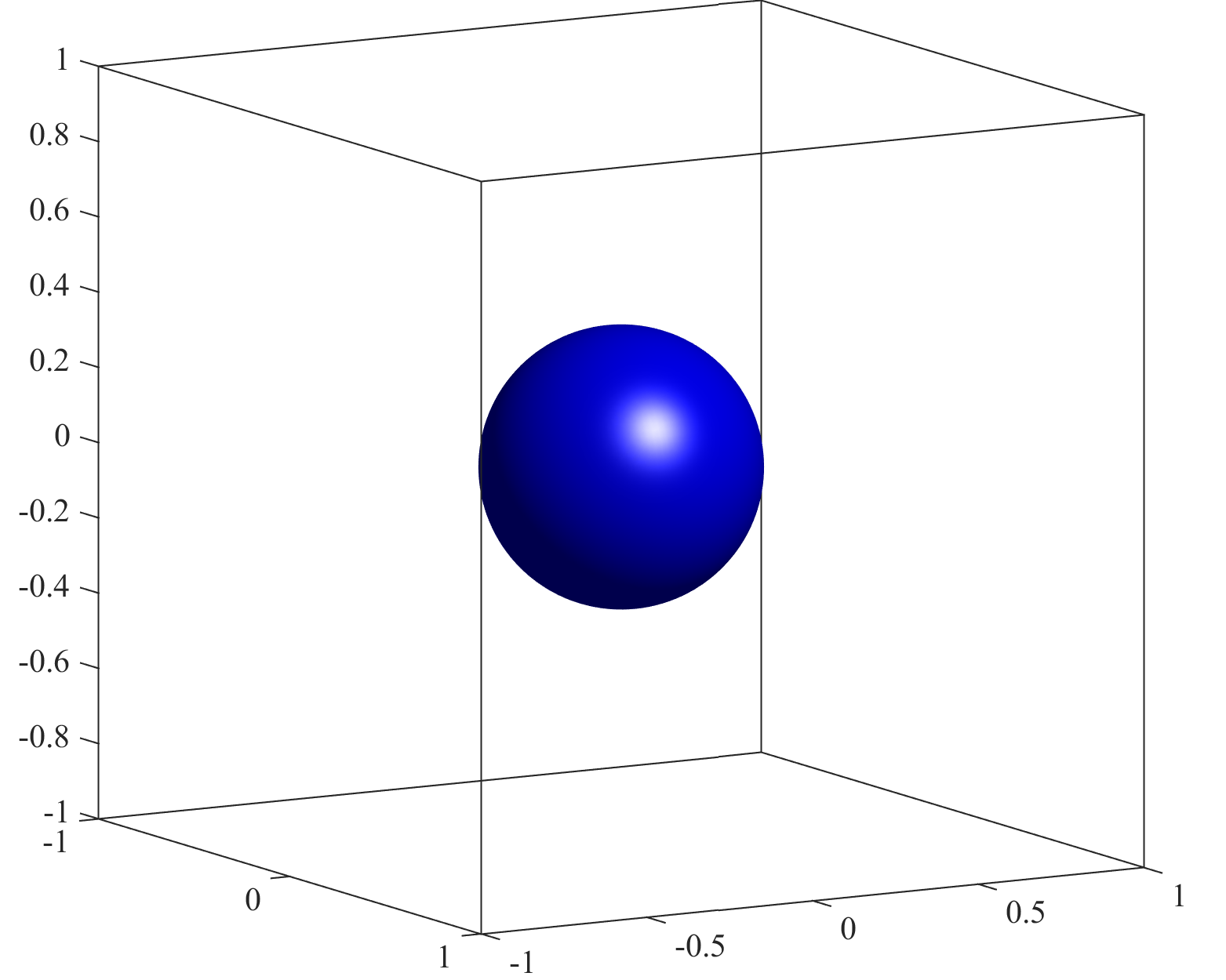} 
			\end{minipage} 
		}
		\subfloat[$T=70$]{
			\begin{minipage}[t]{0.22\linewidth}
				\centering
				\includegraphics[width=1.6in]{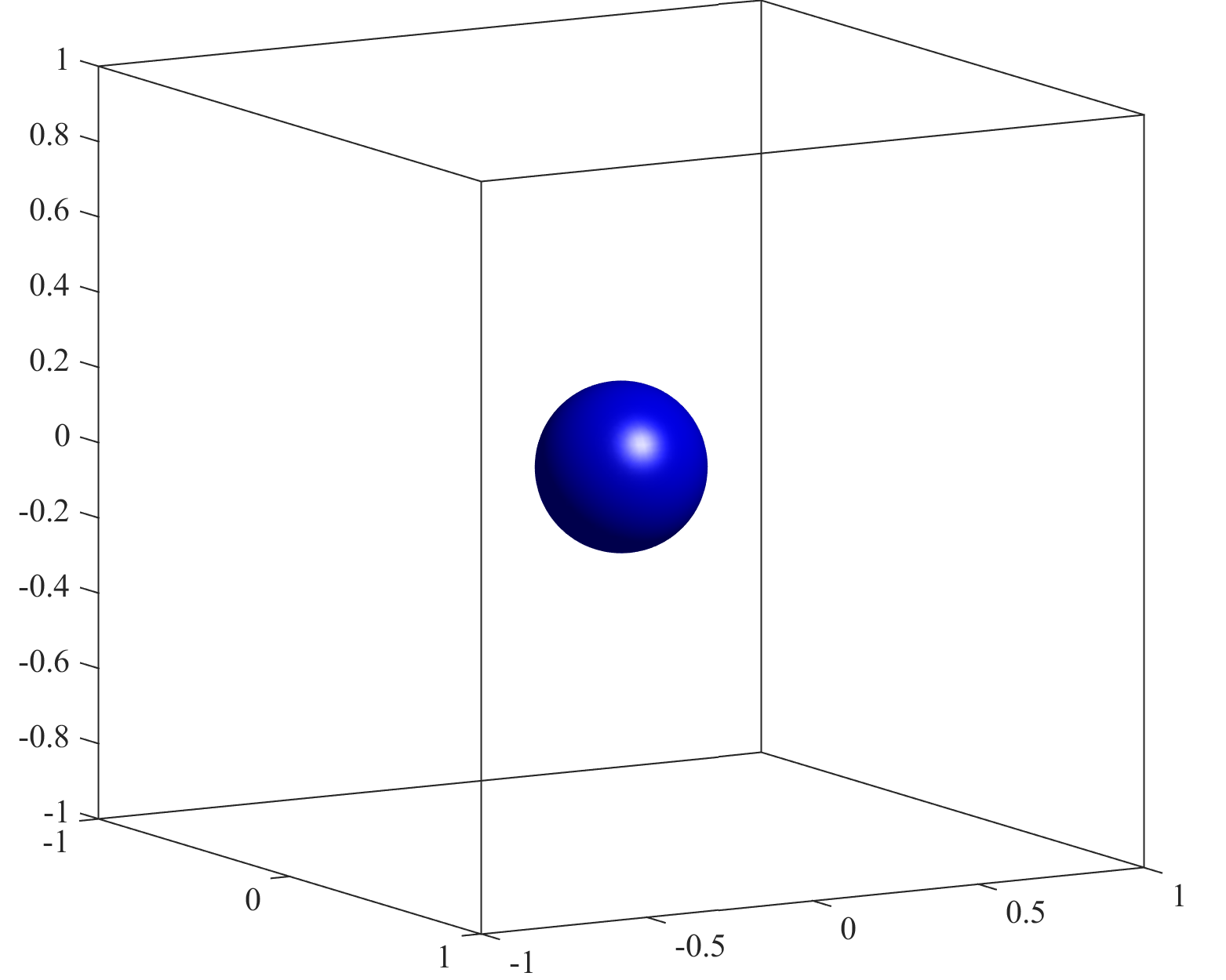} 
			\end{minipage} 
		}
		
		\subfloat[]
		{
			\includegraphics[width=0.3\textwidth]{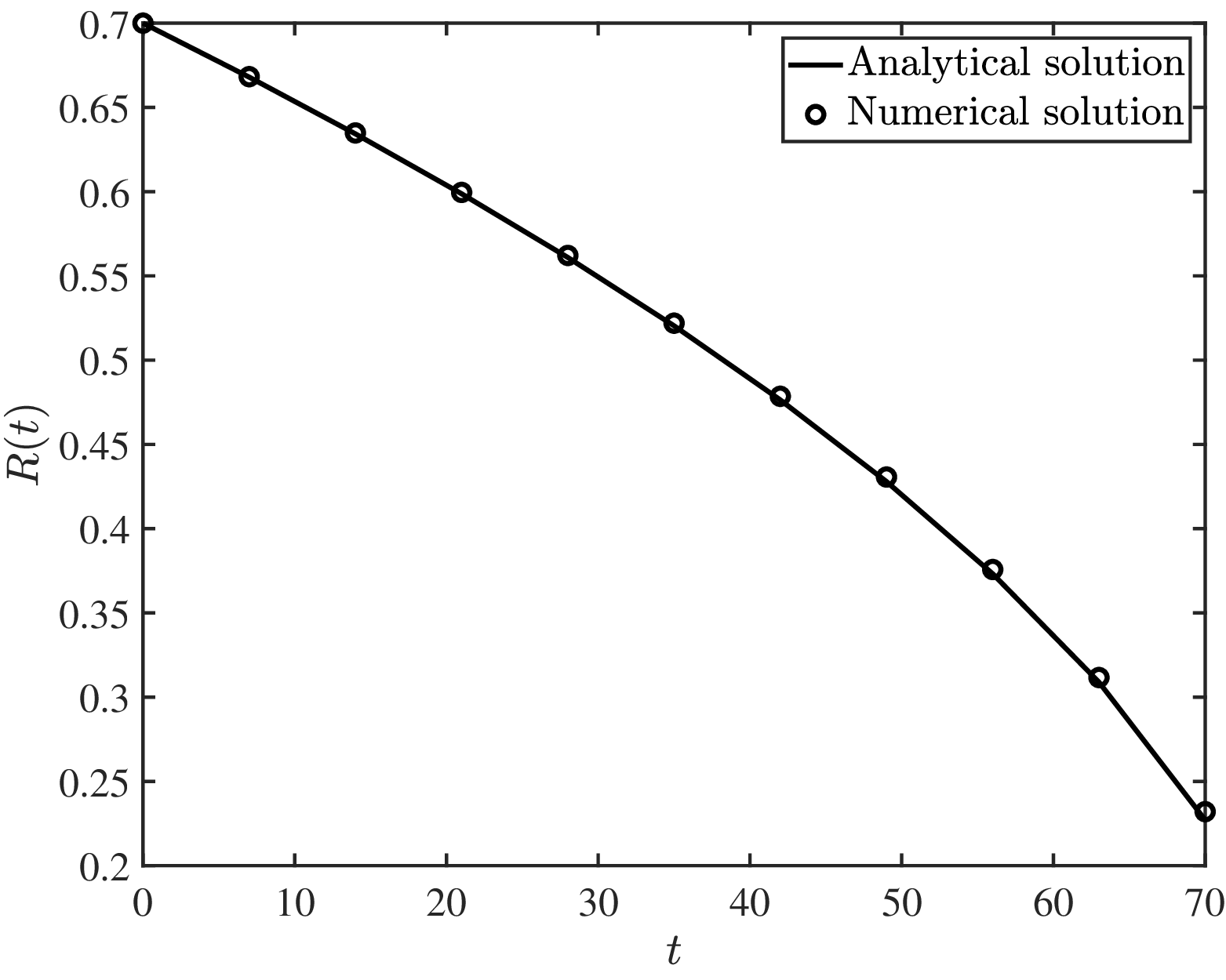}
		} 
		\caption{The evolution of the numerical solution  of the three-dimensional sphere in time.}   	\label{Fig-Ex4-3}   
	\end{center} 
\end{figure}

\begin{figure} [H]
	\vspace{-1.2cm}  
	\setlength{\belowcaptionskip}{-1cm}   
	\begin{center}             
		\subfloat[]    
		{
			\includegraphics[width=0.3\textwidth]{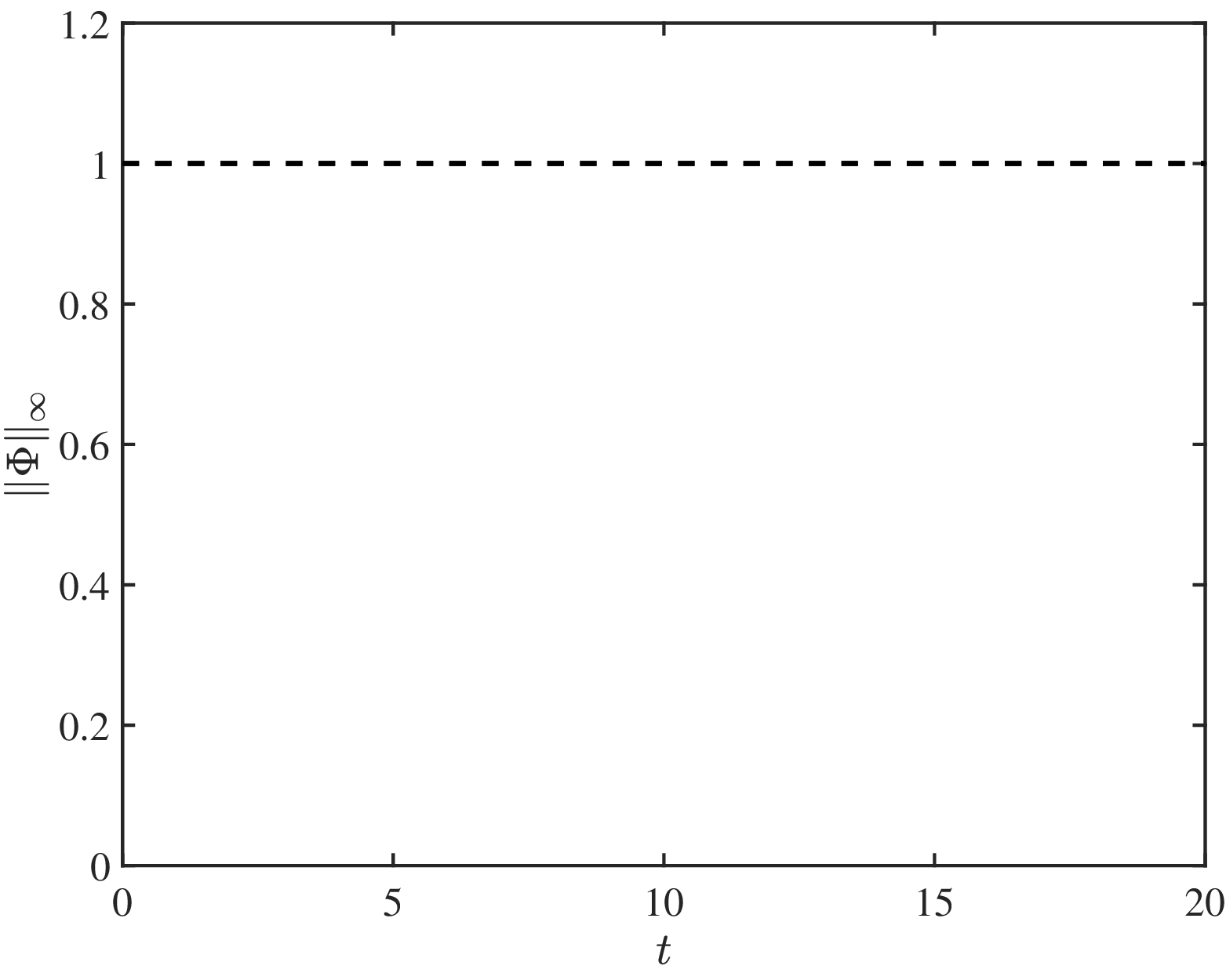}
		} 
		\subfloat[]
		{
			\includegraphics[width=0.3\textwidth]{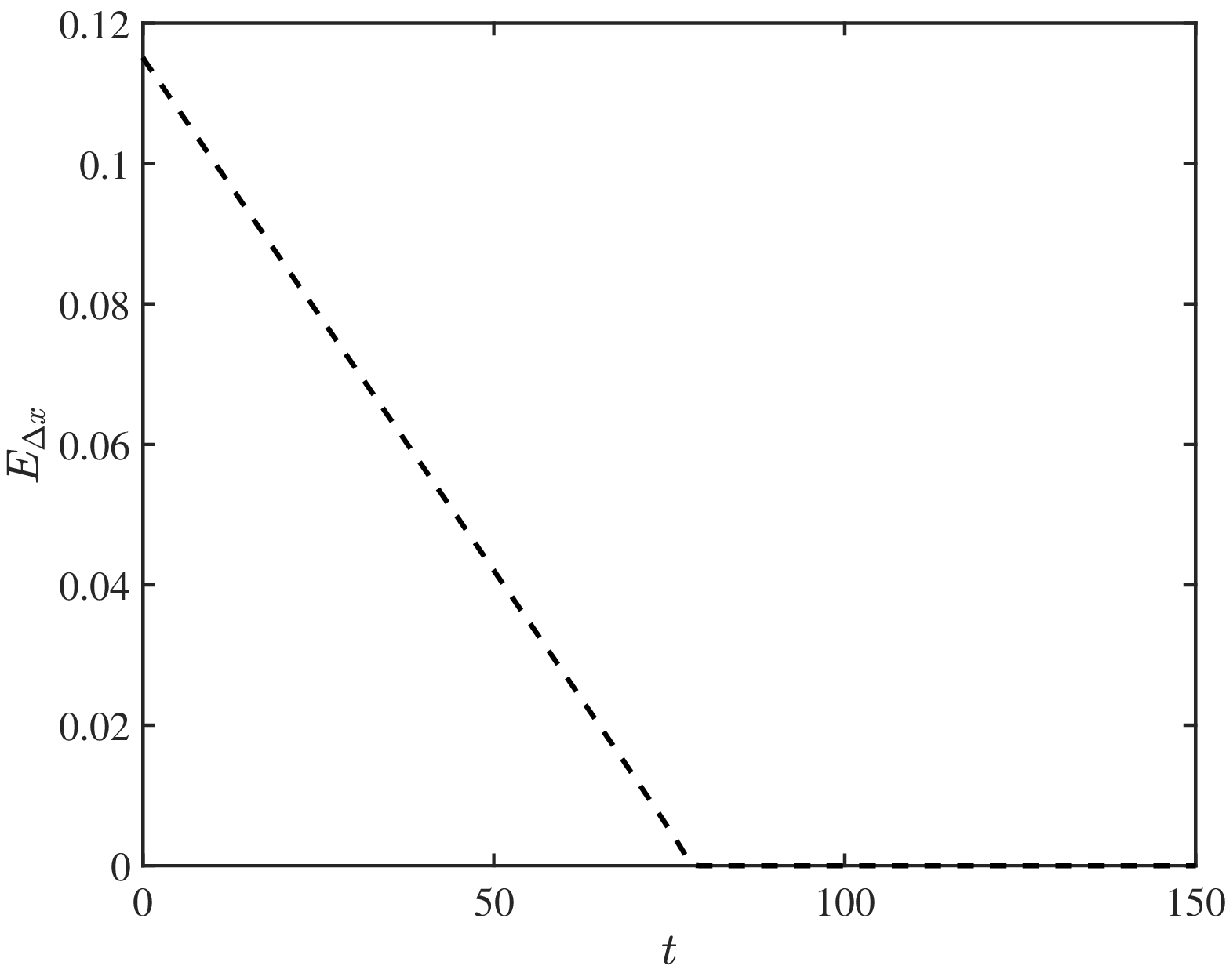}
		} 
		\caption{The evolutions of the discrete maximum norm (a) and energy (b) of the numerical solution in time.}   \label{Fig-Ex4-4}  
	\end{center}     
\end{figure}

To simulate the motion of the curvature flow and to demonstrate that the RLB-MIE-FD scheme is stable under the conditions of Eq. (\ref{maximum-condition}) in Theorem 2, we focus on the temporal evolution of the numerical solution and   the two intrinsic properties of the ACE (\ref{ACE}), i.e., the maximum principle and energy dissipation law. We conduct the numerical simulation with the interfacial parameter $\varepsilon=5\Delta x\tanh(0.9)/\sqrt{2}$, and the lattice spacing and time step are set to be $\Delta x=1/64$ and $\Delta t=1/100$, respectively.
In the two-dimensional case, we consider the computational domain $\Omega=(-1,1)\times (-1,1)$, and the initial condition is given by  a circle with the center ($0,0$) and radius $R_0=0.7$, i.e.,
\begin{align}
	\phi(x,y,0)=\tanh\frac{R_0-\sqrt{x^2+y^2}}{\sqrt{2}\varepsilon}.
\end{align}
Under the HN-BC, one can obtain the analytical solution of radius   $R(t)=\sqrt{R_0^2-2\varepsilon^2t}$ \cite{jeong2018explicit}.  	We present the evolution of the radius of the circle in Fig. \ref{Fig-Ex4-1}, and the numerical result is in good agreement with the analytical solution. In addtion, we also show the evolutions of the discrete maximum norm and energy of the numerical solution in Fig. \ref{Fig-Ex4-2},  and these results are also consistent with our theoretical analysis.
Then we further consider a three-dimensional case with the computational domain $\Omega=(-1,1)\times (-1,1)\times (-1,1)$. The initial condition is given as
\begin{align} 
	\phi(x,y,z,0)=\tanh\frac{R_0-\sqrt{x^2+y^2+z^2}}{\sqrt{2}\varepsilon},
\end{align} 
and the analytical radius can be given by $R(t)=\sqrt{R_0^2-4\varepsilon^2t}$ \cite{jeong2018explicit}. The zero level isosurfaces of the numerical solutions  at several times are plotted in Fig. \ref{Fig-Ex4-3}. From this figure, one can observe that the numerical solution is in good agreement with the analytical solution [see Fig. \ref{Fig-Ex4-4}(e)]. Figs. \ref{Fig-Ex4-4}(a) and \ref{Fig-Ex4-4}(b) show the evolutions of the discrete maximum norm and energy of the numerical solution, respectively,  which agree well with the theoretical analysis.

From the results in Figs. \ref{Fig-Ex4-1}, \ref{Fig-Ex4-2}, \ref{Fig-Ex4-3} and \ref{Fig-Ex4-4}, one can find that the RLB-MIE-FD scheme is not only stable, but also accurate for the mean curvature  flow problems.

\section{Conclusion}\label{conclusion}
In this paper, we first developed a RLB-MIE-FD scheme for the ACE (\ref{ACE}), and theoretically demonstrated that it has the first-order accuracy in time and second-order accuracy in space. Then, we also proved that the developed RLB-MIE-FD scheme can preserve the maximum principle and energy dissipation law at the discrete level under the conditions of Eq. (\ref{maximum-condition}) in Theorem 2. Finally, we carried out some simulations, and the numerical results show that the RLB-MIE-FD scheme not only has a second-order CR in space, but also preserves the discrete maximum principle and energy dissipation law well, which are in good agreement with our theoretical analysis. What is more, it is worth noting that the RLB-MIE-FD scheme is more numerically stable than the FEX-FD scheme, and has the advantage in terms of memory occupancy and efficiency compared to the RLB method and CN scheme.
\section*{Acknowledgements}
The computation is completed in the HPC Platform of Huazhong University of Science and Technology. This work was financially supported by the National Natural Science Foundation of China (Grants No. 12072127 and No. 51836003), the Interdisciplinary Research Program of Hust (2023JCJY002), and the Fundamental Research Funds for the Central Universities, Hust (No. 2023JYCXJJ046).

\bibliographystyle{elsarticle-num} 
\bibliography{reference}

\begin{thebibliography}{10}
\expandafter\ifx\csname url\endcsname\relax
  \def\url#1{\texttt{#1}}\fi
\expandafter\ifx\csname urlprefix\endcsname\relax\def\urlprefix{URL }\fi
\expandafter\ifx\csname href\endcsname\relax
  \def\href#1#2{#2} \def\path#1{#1}\fi

\bibitem{ilmanen1993convergence}
T.~Ilmanen, Convergence of the {A}llen-{C}ahn equation to brakke's motion by
  mean curvature, J. Differ. Geom. 38~(2) (1993) 417--461.

\bibitem{feng2003numerical}
X.~Feng, A.~Prohl, Numerical analysis of the {A}llen-{C}ahn equation and
  approximation for mean curvature flows, Numer. Math. 94 (2003) 33--65.

\bibitem{lee2015mean}
D.~Lee, J.~Kim, Mean curvature flow by the {A}llen-{C}ahn equation, Eur. J.
  Appl. Math. 26~(4) (2015) 535--559.

\bibitem{li2023second}
Y.~Li, K.~Qin, Q.~Xia, J.~Kim, A second-order unconditionally stable method for
  the anisotropic dendritic crystal growth model with an orientation-field,
  Appl. Numer. Math. 184 (2023) 512--526.

\bibitem{benevs2004geometrical}
M.~Bene{\v{s}}, V.~Chalupeck{\`y}, K.~Mikula, Geometrical image segmentation by
  the {A}llen-{C}ahn equation, Appl. Numer. Math. 51~(2-3) (2004) 187--205.

\bibitem{kim2012phase}
J.~Kim, Phase-field models for multi-component fluid flows, Commun. Comput.
  Phys. 12~(3) (2012) 613--661.

\bibitem{Wang2019}
H.~Wang, X.~Yuan, H.~Liang, Z.~Chai, B.~Shi, A brief review of the
  phase-field-based lattice {B}oltzmann method for multiphase flows,
  Capillarity 2~(3) (2019) 33--52.

\bibitem{yang2018uniform}
J.~Yang, Q.~Du, W.~Zhang, Uniform ${L}^p$-bound of the {A}llen-{C}ahn equation
  and its numerical discretization, Int. J. Numer. Anal. Mod. 15 (2018).

\bibitem{tang2016implicit}
T.~Tang, J.~Yang, Implicit-explicit scheme for the {A}llen-{C}ahn equation
  preserves the maximum principle, J. Comput. Math. (2016) 451--461.

\bibitem{shen2016maximum}
J.~Shen, T.~Tang, J.~Yang, On the maximum principle preserving schemes for the
  generalized {A}llen-{C}ahn equation, Commun. Math. Sci. 14~(6) (2016)
  1517--1534.

\bibitem{ham2023stability}
S.~Ham, J.~Kim, Stability analysis for a maximum principle preserving explicit
  scheme of the {A}llen-{C}ahn equation, Math. Comput. Simulat. 207 (2023)
  453--465.

\bibitem{hou2017numerical}
T.~Hou, T.~Tang, J.~Yang, Numerical analysis of fully discretized
  {C}rank-{N}icolson scheme for fractional-in-space {A}llen-{C}ahn equations,
  J. Sci. Comput. 72 (2017) 1214--1231.

\bibitem{hou2020numerical}
T.~Hou, H.~Leng, Numerical analysis of a stabilized
  {C}rank-{N}icolson/{A}dams-{B}ashforth finite difference scheme for
  {A}llen-{C}ahn equations, Appl. Math. Lett. 102 (2020) 106150.

\bibitem{feng2021maximum}
J.~Feng, Y.~Zhou, T.~Hou, A maximum-principle preserving and unconditionally
  energy-stable linear second-order finite difference scheme for
  {A}llen-{C}ahnn equations, Appl. Math. Lett. 118 (2021) 107179.

\bibitem{liao2020energy}
H.-L. Liao, T.~Tang, T.~Zhou, On energy stable, maximum-principle preserving,
  second-order {BDF} scheme with variable steps for the {A}llen-{C}ahn
  equation, SIAM J. Numer. Anal. 58~(4) (2020) 2294--2314.

\bibitem{du2019maximum}
Q.~Du, L.~Ju, X.~Li, Z.~Qiao, Maximum principle preserving exponential time
  differencing schemes for the nonlocal {A}llen-{C}ahn equation, SIAM J. Numer.
  Anal. 57~(2) (2019) 875--898.

\bibitem{fu2022energy1}
L.~Ju, X.~Li, Z.~Qiao, Generalized {SAV}-exponential integrator schemes for
  {A}llen-{C}ahn type gradient flows, SIAM J. Numer. Anal. 60~(4) (2022)
  1905--1931.

\bibitem{fu2022energy2}
L.~Ju, X.~Li, Z.~Qiao, Stabilized exponential-{SAV} schemes preserving energy
  dissipation law and maximum bound principle for the {A}llen-{C}ahn type
  equations, J. Sci. Comput. 92~(2) (2022) 66.

\bibitem{li2020arbitrarily}
B.~Li, J.~Yang, Z.~Zhou, Arbitrarily high-order exponential cut-off methods for
  preserving maximum principle of parabolic equations, SIAM J. Sci. Comput.
  42~(6) (2020) A3957--A3978.

\bibitem{cheng2022new1}
Q.~Cheng, J.~Shen, A new lagrange multiplier approach for constructing
  structure preserving schemes, {I}. {P}ositivity preserving, Comput. Methods
  Appl. Mech. Eng. 391 (2022) 114585.

\bibitem{cheng2022new2}
Q.~Cheng, J.~Shen, A new lagrange multiplier approach for constructing
  structure preserving schemes, {I}{I}. {B}ound preserving, SIAM J. Numer.
  Anal. 60~(3) (2022) 970--998.

\bibitem{yang2022arbitrarily}
J.~Yang, Z.~Yuan, Z.~Zhou, Arbitrarily high-order maximum bound preserving
  schemes with cut-off postprocessing for {A}llen-{C}ahn equations, J. Sci.
  Comput. 90~(2) (2022) 76.

\bibitem{akrivis2019energy}
G.~Akrivis, B.~Li, D.~Li, Energy-decaying extrapolated {RK}-{SAV} methods for
  the {A}llen-{C}ahn and {C}ahn-{H}illiard equations, SIAM J. Sci. Comput.
  41~(6) (2019) A3703--A3727.

\bibitem{shen2018scalar}
J.~Shen, J.~Xu, J.~Yang, The scalar auxiliary variable ({SAV}) approach for
  gradient flows, J. Comput. Phys. 353 (2018) 407--416.

\bibitem{zhang2021numerical}
H.~Zhang, J.~Yan, X.~Qian, S.~Song, Numerical analysis and applications of
  explicit high order maximum principle preserving integrating factor
  {R}unge-{K}utta schemes for {A}llen-{C}ahn equation, Appl. Numer. Math. 161
  (2021) 372--390.

\bibitem{zhang2021preserving}
H.~Zhang, J.~Yan, X.~Qian, X.~Gu, S.~Song, On the preserving of the maximum
  principle and energy stability of high-order implicit-explicit runge-kutta
  schemes for the space-fractional {A}llen-{C}ahn equation, Numer. Algor. 88
  (2021) 1309--1336.

\bibitem{li2021stabilized}
J.~Li, X.~Li, L.~Ju, X.~Feng, Stabilized integrating factor {R}unge-{K}utta
  method and unconditional preservation of maximum bound principle, SIAM J.
  Sci. Comput. 43~(3) (2021) A1780--A1802.

\bibitem{Guo2013}
Z.~Guo, C.~Shu, Lattice {B}oltzmann method and its application in engineering,
  World Scientific, 2013.

\bibitem{Succi2008}
S.~Succi, The lattice {B}oltzmann equation: for fluid dynamics and beyond,
  Oxford university press, 2001.

\bibitem{Kruger2017}
T.~Kr{\"u}ger, H.~Kusumaatmaja, A.~Kuzmin, O.~Shardt, G.~Silva, E.~Viggen, The
  {L}attice {B}oltzmann {M}ethod: {P}rinciples and {P}ractice, Oxford
  university press, 2017.

\bibitem{liu2022diffuse}
X.~Liu, Z.~Chai, C.~Zhan, B.~Shi, W.~Zhang, A diffuse-domain phase-field
  lattice boltzmann method for two-phase flows in complex geometries,
  Multiscale Model. Sim. 20~(4) (2022) 1411--1436.

\bibitem{liu2023improved}
X.~Liu, Z.~Chai, B.~Shi, Improved hybrid {A}llen-{C}ahn phase-field-based
  lattice {B}oltzmann method for incompressible two-phase flows, Phys. Rev. E
  107~(3) (2023) 035308.

\bibitem{xiao2022second}
X.~Xiao, X.~Feng, A second-order maximum bound principle preserving operator
  splitting method for the {A}llen-{C}ahn equation with applications in
  multi-phase systems, Math. Comput. Simulat. 202 (2022) 36--58.

\bibitem{simonis2023lattice}
S.~Simonis, Lattice {B}oltzmann {Me}thods for {P}artial {D}ifferential
  {E}quations, Ph.D. thesis, Dissertation, Karlsruhe, Karlsruher Institut
  f{\"u}r Technologie (KIT) (2023).

\bibitem{Suga2010}
S.~Suga, An accurate multi-level finite difference scheme for 1{D} diffusion
  equations derived from the lattice {B}oltzmann method, J. Stat. Phys. 140
  (2010) 494--503.

\bibitem{Lin2022}
Y.~Lin, N.~Hong, B.~Shi, Z.~Chai, Multiple-relaxation-time lattice {B}oltzmann
  model-based four-level finite-difference scheme for one-dimensional diffusion
  equations, Phys. Rev. E 104~(1) (2021) 015312.

\bibitem{silva2023discrete}
G.~Silva, Discrete effects on the source term for the lattice {B}oltzmann
  modelling of one-dimensional reaction-diffusion equations, Comput. Fluids 251
  (2023) 105735.

\bibitem{Chen2023}
Y.~Chen, Z.~Chai, B.~Shi, A general fourth-order mesoscopic
  multiple-relaxation-time lattice {B}oltzmann model and equivalent macroscopic
  finite-difference scheme for two-dimensional diffusion equations, arXiv
  preprint arXiv:2308.05280 (2023).

\bibitem{straka2020accuracy}
R.~Straka, K.~Sharma, An accuracy analysis of the cascaded lattice {B}oltzmann
  method for the 1{D} advection-diffusion equation, Comput. Methods Mater. Sci.
  20~(4) (2020) 174.

\bibitem{chen2023fourth}
Y.~Chen, Z.~Chai, B.~Shi, Fourth-order multiple-relaxation-time lattice
  {B}oltzmann model and equivalent finite-difference scheme for one-dimensional
  convection-diffusion equations, Phys. Rev. E 107~(5) (2023) 055305.

\bibitem{li2012multilevel}
Q.~Li, Z.~Zheng, S.~Wang, J.~Liu, A multilevel finite difference scheme for
  one-dimensional {B}urgers equation derived from the lattice {B}oltzmann
  method, J. Appl. Math. 2012 (2012).

\bibitem{chen2023cole}
Y.~Chen, X.~Liu, Z.~Chai, B.~Shi, A {C}ole-{H}opf transformation based
  fourth-order multiple-relaxation-time lattice {B}oltzmann model for the
  coupled {B}urgers' equations, arXiv preprint arXiv:2309.02825 (2023).

\bibitem{junk2001finite}
M.~Junk, A finite difference interpretation of the lattice {B}oltzmann method,
  Numer. Methods Part. Diff. Equ. 17~(4) (2001) 383--402.

\bibitem{Be2023}
T.~Bellotti, B.~Graille, M.~Massot, Finite {D}ifference formulation of any
  lattice {B}oltzmann scheme, Numer. Math. 152~(1) (2022) 1--40.

\bibitem{Chen2023-GP}
Y.~Chen, X.~Liu, Z.~Chai, B.~Shi, The macroscopic finite-difference scheme and
  modified equations of the general propagation multiple-relaxation-time
  lattice {B}oltzmann model, arXiv preprint arXiv:2308.11882 (2023).

\bibitem{Liu2023rlb}
X.~Liu, Y.~Chen, Z.~Chai, B.~Shi, Efficient macroscopic finite-difference
  scheme of regularized lattice {B}oltzmann method, Sumbitted (2023).

\bibitem{Chai2020square}
Z.~Chai, B.~Shi, Multiple-relaxation-time lattice {B}oltzmann method for the
  {N}avier-{S}tokes and nonlinear convection-diffusion equations: {M}odeling,
  analysis, and elements, Phys. Rev. E 102~(2) (2020) 023306.

\bibitem{wang2015regularized}
L.~Wang, B.~Shi, Z.~Chai, Regularized lattice {B}oltzmann model for a class of
  convection-diffusion equations, Phys. Rev. E 92~(4) (2015) 043311.

\bibitem{lee2020numerical}
D.~Lee, The numerical solutions for the energy-dissipative and
  mass-conservative {A}llen-{C}ahn equation, Comput. Math. Appl. 80~(1) (2020)
  263--284.

\bibitem{varga1976m}
R.~S. Varga, M-matrix theory and recent results in numerical linear algebra,
  in: Sparse matrix computations, Elsevier, 1976, pp. 375--387.

\bibitem{jeong2018explicit}
D.~Jeong, J.~Kim, An explicit hybrid finite difference scheme for the
  {A}llen-{C}ahn equation, J. Comput. Appl. Math. 340 (2018) 247--255.

\end{thebibliography}

\end{document}